\documentclass[12pt,reqno]{amsart}

\usepackage{amssymb}

\def\II{j}
\def\JJ{i}

\usepackage{color}
\def\Rot{\color{red}} 
\def\Blau{\color{blue}} 
\def\Gruen{\color{green}} 
 
\font\Riesig=cmr10 at 80pt
\font\riesig=cmr10 at 50pt
\font\grosz=cmr10 at 20pt

\input texdraw
\def\Ringerl(#1 #2){\move(#1 #2)\fcir f:0 r:.15}
\def\ringerl(#1 #2){\move(#1 #2)\fcir f:0 r:.1}
\def\Mark(#1 #2){\move(#1 #2)\fcir f:0 r:.2}

\newtheorem{theorem}{Theorem}

\newtheorem{lemma}[theorem]{Lemma}
\newtheorem{corollary}[theorem]{Corollary}
\newtheorem*{conjectureFM}{Conjecture FM}

\theoremstyle{remark}
\newtheorem*{remark}{Remark}

\numberwithin{equation}{section}

\def\({\left(}
\def\){\right)}
\def\coef#1{\left\langle#1\right\rangle}

\def\rk{\operatorname{rk}}
\def\Nar{\operatorname{Nar}}

\def\de{\delta}
\def\om{\omega}
\def\ph{\varphi}
\def\si{\sigma}
\def\ep{\varepsilon}
\def\al{\alpha}

\def\ga{\gamma}
\def\la{\lambda}

\def\ka{\kappa}
\def\ta{\tau}
\def\De{\Delta}

\def\Na{\bigtriangledown}

\begin{document}

\title[Decomposition numbers for finite Coxeter groups]{%
Decomposition numbers for finite Coxeter groups
and generalised non-crossing partitions}

\author[C. Krattenthaler and T. W. M\"uller]%
{C. Krattenthaler$^\dagger$ and T. W. M\"uller}

\address{Fakult\"at f\"ur Mathematik, Universit\"at Wien,
Nordbergstra{\ss}e~15, A-1090 Vienna, Austria.
WWW: \tt http://www.mat.univie.ac.at/\~{}kratt.}

\address{School of Mathematical Sciences, Queen Mary
\& Westfield College, University of London,
Mile End Road, London E1 4NS, United Kingdom.\newline
WWW: \tt http://www.maths.qmw.ac.uk/\~{}twm/.}

\thanks{$^\dagger$Research partially supported 
by the Austrian Science Foundation FWF, grant S9607-N13,
in the framework of the National Research Network
``Analytic Combinatorics and Probabilistic Number Theory"}


\begin{abstract}
Given a finite irreducible Coxeter group $W$, a positive integer $d$, 
and types
$T_1,T_2,\dots,T_d$ (in the sense of the classification of finite
Coxeter groups), we compute the number of
decompositions $c=\si_1\si_2\cdots\si_d$ of a Coxeter element $c$ of
$W$, such that $\si_i$ is a Coxeter element in a subgroup of type
$T_i$ in $W$, $i=1,2,\dots,d$, and such that the factorisation is 
``minimal" in the sense that the sum of the ranks of the $T_i$'s,
$i=1,2,\dots,d$, equals the rank of $W$. For the exceptional types,
these decomposition numbers have been computed by the first author in
[{\it``Topics in Discrete Mathematics,''} M.~Klazar et al. (eds.), 
Springer--Verlag, Berlin, New York, 2006, pp.~93--126]
and [{\it S\'e\-mi\-naire Lotharingien Combin.} {\bf 54} (2006),
Article~B54l]. The type $A_n$ decomposition
numbers have been computed by Goulden and Jackson in [{\it Europ.\ J.
Combin.} {\bf 13} (1992), 357--365], albeit using a somewhat different
language. We explain how to extract the type $B_n$
decomposition numbers from results of
B\'ona, Bousquet, Labelle and Leroux [Adv.\ Appl.\ Math.\ {\bf 24} 
(2000), 22--56] on map enumeration. Our formula for the type $D_n$ 
decomposition numbers is new. These results are then used to
determine, for a fixed positive integer $l$ and fixed integers
$r_1\le r_2\le \dots\le r_l$, the number of multi-chains 
$\pi_1\le \pi_2\le \dots\le \pi_l$ in Armstrong's generalised non-crossing 
partitions poset, where the poset rank of $\pi_i$ equals
$r_i$, and where the ``block structure" of $\pi_1$ is prescribed.
We demonstrate that this result implies all known enumerative results on
ordinary and generalised non-crossing partitions via appropriate
summations. Surprisingly, this result on multi-chain enumeration
is new even for the original
non-crossing partitions of Kreweras. Moreover, the result allows one
to solve the problem of rank-selected chain enumeration
in the type $D_n$ generalised non-crossing partitions poset, which, in
turn, leads to a proof of Armstrong's $F=M$ Conjecture in type
$D_n$, thus completing a computational proof of the $F=M$ Conjecture
for all types. It also allows to address another conjecture of Armstrong
on maximal intervals containing a random multichain in the 
generalised non-crossing partitions poset.
\end{abstract}

\maketitle

\section{Introduction}
\label{sec:0} 

The introduction of {\it non-crossing partitions for finite reflection 
groups} (finite Coxeter groups) by Bessis \cite{BesDAA} and 
Brady and Watt \cite{BRWaAA}
marks the creation of a new, exciting subject of combinatorial
theory, namely the study of these new combinatorial objects which
possess numerous beautiful properties, and seem to relate to
several other objects of combinatorics and algebra, most notably to the
cluster complex of Fomin and Zelevinsky \cite{FOZeAB} (cf.\
\cite{AthaAE,AthaAI,AtBWAA,AtReAA,BesDAA,BeCoAA,BraTAA,BRWaAA,BRWaAB,ChaFAA,FoReAB}). 
They reduce to
the classical non-crossing partitions of Kreweras \cite{KrewAC} for
the irreducible reflection groups of type $A_{n}$ (i.e., the
symmetric groups), and to Reiner's \cite{ReivAG} type $B_n$
non-crossing partitions for the irreducible reflections groups of
type $B_n$. (They differ, however, from the type $D_n$ non-crossing
partitions of \cite{ReivAG}.) The subject 
has been enriched by Armstrong through introduction of his
{\it generalised non-crossing partitions for reflection groups} in
\cite{ArmDAA}. In the symmetric group case, these
reduce to the $m$-divisible non-crossing partitions of
Edelman \cite{EdelAA}, while they produce new combinatorial objects
already for the reflection groups of type $B_n$. Again, these
generalised non-crossing partitions
possess numerous beautiful properties, and seem to relate to
several other objects of combinatorics and algebra, most
notably to the generalised cluster complex of Fomin and Reading
\cite{FoReAA} (cf.\
\cite{ArmDAA,AtTzAA,AtTzAB,FoReAB,KratCB,KratCF,KratCG,TzanAA,TzanAC,TzanAB}). 

From a technical point of view, the {\it main subject matter\/} of the present paper is the
computation of the number of certain factorisations of the Coxeter
element of a reflection group. These {\it decomposition numbers}, as
we shall call them from now on (see Section~\ref{sec:1} for the
precise definition),
arose in \cite{KratCB,KratCF}, where it was shown that they
play a crucial role in the computation of enumerative invariants of
(generalised) non-crossing partitions. Moreover, in these two papers
the decomposition numbers for the exceptional reflection groups have
been computed, and it was pointed out that the decomposition numbers
in type $A_n$ (i.e., the decomposition numbers for the symmetric
groups) had been earlier computed by Goulden and Jackson in
\cite{GoJaAS}. Here, we explain how the decomposition numbers in
type $B_n$ can be extracted from results of B\'ona, Bousquet, Labelle
and Leroux \cite{BoBLAA} on the enumeration of certain planar maps, 
and we find formulae for the decomposition numbers in type $D_n$, 
thus completing the project of computing the decomposition numbers
for all the irreducible reflection groups.

The {\it main goal\/} of the present paper, however, is to access the
enumerative theory of the generalised non-crossing
partitions of Armstrong via these decomposition numbers.
Indeed, one finds numerous enumerative results on ordinary and
generalised non-crossing
partitions in the literature (cf.\ 
\cite{ArmDAA,AthaAE,AtReAA,BesDAA,BeCoAA,EdelAA,KrewAC,ReivAG,TzanAC}): 
results on the total number of (generalised) non-crossing partitions of a given size, 
of those with a fixed number of blocks, of those with a given block
structure, results on the number of (multi-)chains of a given length in a given poset
of (generalised) non-crossing partitions, results on rank-selected
chain enumeration (that is, results on the number of chains in
which the ranks of the elements of the chains have been fixed), etc.
We show that not only can {\it all\/} these results be rederived from
our decomposition numbers, we are also able to find several new enumerative
results. In this regard, the most general type of result that we find is
formulae for the number of (multi-)chains $\pi_1\le \pi_2\le \dots\le
\pi_{l-1}$ in the poset of non-crossing partitions of type $A_n$,
$B_n$, respectively $D_n$, in which the block structure of $\pi_1$ is
fixed as well as the ranks of $\pi_2,\dots,\pi_{l-1}$. Even the corresponding
result in type $A_n$, for the non-crossing partitions of
Kreweras, is new. Furthermore, from the result in type $D_n$, 
by a suitable summation, we are able to find 
a formula for the rank-selected chain enumeration in the poset of 
generalised non-crossing partitions of type $D_n$, thus generalising
the earlier formula of Athanasiadis and Reiner \cite{AtReAA} for the
rank-selected chain enumeration of ``ordinary" 
non-crossing partitions of type $D_n$. In conjunction with the
results from \cite{KratCB,KratCF}, this generalisation in turn
allows us to complete a
computational case-by-case proof of Armstrong's ``$F=M$ Conjecture" 
\cite[Conjecture~5.3.2]{ArmDAA} predicting 
a surprising relationship between a certain face count in the
generalised cluster complex of Fomin and Reading and the M\"obius
function in the poset of generalised non-crossing partitions of
Armstrong. (A case-free proof had been found earlier by Tzanaki in
\cite{TzanAB}.) 
Our results allow us also to address another conjecture of Armstrong 
\cite[Conj.~3.5.13]{ArmDAA} on maximal intervals containing a random
multichain in the poset of generalised
non-crossing partitions. We show that the conjecture is indeed true
for types $A_n$ and $B_n$, but that it fails for type $D_n$ (and 
we suspect that it will also fail for most of the exceptional types).

We remark that a
totally different approach to the enumerative theory of (generalised)
non-crossing partitions is proposed in \cite{KratCG}. This approach
is, however, completely combinatorial and avoids, in particular, 
reflection groups. It is, therefore, not capable of
computing our decomposition numbers nor anything else which is
intrinsic to the combinatorics of reflection groups. 
A similar remark applies to \cite[Theorem~4.1]{ReadAA}, 
where a remarkable
uniform recurrence is found for rank-selected chain enumeration
in the generalised non-crossing partitions of {\it any} type.
It could be used, for example, for verifying our result 
in Corollary~\ref{cor:6} on the
rank-selected chain enumeration in the generalised non-crossing
partitions of type $D_n$, but it is not capable of
computing our decomposition numbers nor of verifying results with restrictions
on block structure. 

Our paper is organised as follows. In the next section we define the
decomposition numbers for finite reflection groups from \cite{KratCB,KratCF}, 
the central objects in our paper, together with a combinatorial
variant, which depends on combinatorial realisations of non-crossing
partitions, which we also explain in the same section.
This is followed by an intermediate
section in which we collect together some auxiliary results that will be needed 
later on. In Section~\ref{sec:2}, we recall Goulden and Jackson's
formula \cite{GoJaAS} for the full rank decomposition numbers of type
$A_n$, together with the formula from \cite[Theorem~10]{KratCF} that
it implies for the decomposition numbers of type $A_n$ of arbitrary rank.
The purpose of Section~\ref{sec:3} is to explain how formulae for
the decomposition numbers of type $B_n$ can be extracted from results
of B\'ona, Bousquet, Labelle and Leroux in \cite{BoBLAA}. The type
$D_n$ decomposition numbers are computed in Section~\ref{sec:4}. The
approach that we follow is, essentially, the approach of Goulden and
Jackson in \cite{GoJaAS}: we translate the counting problem into the
problem of enumerating certain maps. This problem is then solved by a
combinatorial decomposition of these maps, translating the
decomposition into a system of equations for corresponding generating
functions, and finally solving this system with the help of the
multidimensional Lagrange inversion formula of Good. 
Sections~\ref{sec:6}--\ref{sec:9} form the ``applications part" of the paper.
In the preparatory Section~\ref{sec:6}, 
we recall the definition of the generalised non-crossing
partitions of Armstrong, and we explain the combinatorial
realisations of the generalised non-crossing partitions for the
types $A_n$, $B_n$, and $D_n$ from \cite{ArmDAA} and
\cite{KratCG}. The bulk of the applications is contained in
Section~\ref{sec:7}, where we present three theorems,
Theorems~\ref{thm:4}, \ref{thm:5}, and \ref{thm:6}, on the
number of factorisations of a Coxeter element of type
$A_n$, $B_n$, respectively $D_n$, with less stringent restrictions on
the factors than for the decomposition numbers. These theorems result
from our formulae for the (combinatorial) decomposition numbers upon
appropriate summations. Subsequently, it is shown that the corresponding
formulae imply all known enumeration results on non-crossing
partitions and generalised non-crossing partitions, plus several new
ones, see Corollaries~\ref{cor:2}, \ref{cor:3}, \ref{cor:4}--\ref{cor:6}
and the accompanying remarks. Section~\ref{sec:8} presents the
announced computational proof of the $F=M$ (ex-)Conjecture for type
$D_n$, based on our formula in Corollary~\ref{cor:6} for the
rank-selected chain enumeration in the poset of generalised
non-crossing partitions of type $D_n$, while Section~\ref{sec:8a}
addresses Conjecture~3.5.13 from \cite{ArmDAA}, showing that it does
not hold in general since it fails in type $D_n$.
In the final
Section~\ref{sec:9} we point out that the decomposition numbers do
not only allow one to derive enumerative results for the generalised
non-crossing partitions of the classical types, they also provide all the
means for doing this for the exceptional types. For the convenience
of the reader, we list the values of the decomposition numbers for
the exceptional types that have been computed in \cite{KratCB,KratCF}
in an appendix.

In concluding the introduction, we want to attract the reader's
attention to the fact that many of the formulae presented here are
very combinatorial in nature (see Sections~\ref{sec:2}, \ref{sec:3},
\ref{sec:7}). This raises the natural question as to whether it is
possible to find combinatorial proofs for them. Indeed, a
combinatorial (and, in fact, almost bijective) proof of the formula
of Goulden and Jackson, presented here in Theorem~\ref{thm:1}, has
been given by Bousquet, Chauve and Schaeffer in \cite{BoCSAA}.
Moreover, most of the proofs for the known enumeration results on
(generalised) non-crossing partitions presented in 
\cite{ArmDAA,AthaAE,AtReAA,EdelAA,ReivAG} are combinatorial. On the other
hand, to our knowledge so far nobody has given a combinatorial proof 
for Theorem~\ref{thm:2}, the formula for the decomposition numbers of
type $B_n$, essentially due to B\'ona, Bousquet, Labelle and Leroux
\cite{BoBLAA}, although we believe that this should be possible by
modifying the ideas from \cite{BoCSAA}. There are also other formulae
in our paper (see e.g.\ Corollaries~\ref{cor:2} and \ref{cor:3},
Eqs.~\eqref{eq:6} and \eqref{eq:45a}) which seem amenable to
combinatorial proofs.
However, to find combinatorial proofs
for our type $D_n$ results (cf.\ in particular Theorem~\ref{thm:3}.(ii)
and Corollaries~\ref{cor:4}--\ref{cor:6})
seems rather hopeless to us.

\section{Decomposition numbers for finite Coxeter groups}
\label{sec:1} 

In this section, we introduce the decomposition numbers from
\cite{KratCB,KratCF}, which are (Coxeter) group-theoretical in nature, 
plus combinatorial variants for Coxeter groups of types $B_n$ and
$D_n$, which will be important in combinatorial applications. 
These variants depend on the combinatorial realisation of
these Coxeter groups, which we also explain here.

Let $\Phi$ be a finite root system of rank $n$. (We refer the reader to
\cite{HumpAC} for all terminology on root systems.)
For an element $\al\in\Phi$, let $t_\al$
denote the corresponding reflection in the central hyperplane perpendicular to
$\al$. Let $W=W(\Phi)$ be the group generated by these reflections. 
As is well-known (cf.\ e.g.\ \cite[Sec.~6.4]{HumpAC}), 
any such reflection group is at the same time a
finite Coxeter group, and all finite Coxeter groups can be realised
in this way. By
definition, any element $w$ of $W$ can be represented as a product 
$w=t_1t_2\cdots t_\ell$, where the $t_i$'s are reflections. We call
the minimal number of reflections which is needed for such a
product representation the {\it absolute length\/} of $w$, and we
denote it by $\ell_T(w)$. We then define the {\it absolute order} on
$W$, denoted by $\le_T$, via
$$u\le_T w\quad \text{if and only if}\quad
\ell_T(w)=\ell_T(u)+\ell_T(u^{-1}w).$$
As is well-known and easy to see, 
this is equivalent to the statement that every
shortest representation of $u$ by reflections
occurs as an initial segment in some shortest product representation
of $w$ by reflections. 

Now, for a finite root system $\Phi$ of rank $n$, types
$T_1,T_2,\dots,T_d$ (in the sense of the classification of finite
Coxeter groups), and a {\it Coxeter element\/} $c$,
the {\it decomposition number} $N_\Phi(T_1,T_2,\dots,T_d)$ is defined
as the number of ``minimal" products $c_1c_2\cdots c_d$
less than or equal to $c$ in absolute order, 
``minimal" meaning that 
$\ell_T(c_1)+\ell_T(c_2)+\dots+\ell_T(c_d)=\ell_T(c_1c_2\cdots c_d)$, 
such that, for $i=1,2,\dots,d$, the
type of $c_i$ as a parabolic Coxeter element is $T_i$.
(Here, the term ``parabolic Coxeter element" means a Coxeter element
in some parabolic subgroup.
The reader should recall that it follows from \cite[Lemma~1.4.3]{BesDAA}
that any element $c_i$ is indeed a Coxeter element in a parabolic subgroup
of $W=W(\Phi)$. By definition, the type of $c_i$ is the type of this
parabolic subgroup. The reader should also note that, because of the
rewriting
\begin{equation} \label{eq:Umord}
c_1c_2\cdots c_d=c_i(c_i^{-1}c_1c_i)(c_i^{-1}c_2c_i)\cdots
(c_i^{-1}c_{i-1}c_i)c_{i+1}\cdots c_d,
\end{equation}
any $c_i$ in such a minimal product $c_1c_2\cdots c_d\le_T c$ is
itself $\le_T c$.) It is easy to see that 
the decomposition numbers are independent
of the choice of the Coxeter element $c$.
(This follows from the well-known fact that any two Coxeter elements
are conjugate to each other; cf.\ \cite[Sec.~3.16]{HumpAC}.)

The decomposition numbers satisfy several linear relations between
themselves. First of all,
the number $N_\Phi(T_1,T_2,\dots,T_d)$ is independent
of the order of the types $T_1,T_2,\dots,\break T_d$; that is, we have
\begin{equation} \label{Aa} 
N_\Phi(T_{\si(1)},T_{\si(2)},\dots,T_{\si(d)})=N_\Phi(T_1,T_2,\dots,T_d)
\end{equation}
for every permutation $\si$ of $\{1,2,\dots,d\}$. This is, in fact, a
consequence of the rewriting \eqref{eq:Umord}.
Furthermore, by the definition of these numbers, those of ``lower rank" can be
computed from those of ``full rank." To be precise, we have
\begin{equation} \label{Ab} 
N_\Phi(T_1,T_2,\dots,T_d)=
\sum _{T} ^{}N_\Phi(T_1,T_2,\dots,T_d,T),
\end{equation}
where the sum is taken over all types $T$ of rank $n-\rk T_1-\rk T_2-
\dots -\rk T_d$ (with $\rk T$ denoting the rank of the root
system $\Psi$ of type $T$, and $n$ still denoting the rank of the 
fixed root system $\Phi$; for later use we record that
\begin{equation} \label{eq:ellrk}
\ell_T(w_0)=\rk T_0
\end{equation}
for any parabolic Coxeter element $w_0$ of type $T_0$).

The decomposition numbers for the exceptional types have been
computed in \cite{KratCB,KratCF}. For the benefit of the reader,
we reproduce these numbers in the appendix. The 
decomposition numbers
for type $A_n$ are given in Section~\ref{sec:2}, the ones for type
$B_n$ are computed in Section~\ref{sec:3}, while the ones for type
$D_n$ are computed in Section~\ref{sec:4}. 

\medskip
Next we introduce variants of the above decomposition numbers for
the types $B_n$ and $D_n$, which depend on the combinatorial
realisation of the Coxeter groups of these types. 

As is well-known,
the reflection group $W(A_n)$ can be realised as the
symmetric group $S_{n+1}$ on $\{1,2,\dots,n+1\}$. 
The reflection groups $W(B_n)$
and $W(D_n)$, on the other hand, can be realised as 
subgroups of the symmetric group on $2n$ elements. (See e.g.\
\cite[Sections~8.1 and 8.2]{BjBrAB}.)
Namely, the reflection group $W(B_n)$ can be realised as the subgroup
of the group of all permutations $\pi$ of 
\begin{equation*} 
\{1,2,\dots,n,\bar 1,\bar 2,\dots,\bar n\}
\end{equation*}
satisfying the property 
\begin{equation} \label{eq:-}
\pi(\bar i)=\overline{\pi(i)}.
\end{equation}
(Here, and in what follows, $\bar{\bar i}$ 
is identified with $i$ for all $i$.) In this realisation, there is an
analogue of the disjoint cycle decomposition of permutations. Namely,
every $\pi\in W(B_n)$ can be decomposed as
\begin{equation} \label{eq:cycB}
\pi=\ka_1\ka_2\cdots \ka_s,
\end{equation}
where, for $i=1,2,\dots,s$, 
$\ka_i$ is of one of two possible types of ``cycles": 
a {\it type $A$ cycle},
by which we mean a permutation of the form
\begin{equation} \label{eq:acycle}
((a_1,a_2,\dots,a_k)):=(a_1,a_2,\dots,a_k)\,
(\overline{a_1},\overline{a_2},\dots,\overline{a_k}),
\end{equation}
or a {\it type $B$ cycle}, by which we mean a permutation of the form
\begin{equation} \label{eq:bcycle} 
[a_1,a_2,\dots,a_k]:=(a_1,a_2,\dots,a_k,
\overline{a_1},\overline{a_2},\dots,\overline{a_k}),
\end{equation}
$a_1,a_2,\dots,a_k\in\{1,2,\dots,n,\bar 1,\bar 2,\dots,\bar n\}$.
(Here we adopt notation from \cite{BRWaAA}.)
In both cases, we call $k$ the {\it length} of the ``cycle."
The decomposition \eqref{eq:cycB} is unique up to a reordering of the
$\ka_i$'s.

We call a type $A$ cycle of length $k$ {\it of
combinatorial type} $A_{k-1}$, while we call a type $B$ cycle of length
$k$ {\it of combinatorial type} $B_k$, $k=1,2,\dots$. 
The reader should observe that,
when regarded as a parabolic Coxeter element, for $k\ge2$ a type $A$ cycle 
of length $k$ has type
$A_{k-1}$, while a type $B$ cycle of length $k$ has type
$B_{k}$. However, a type $B$ cycle of length $1$, that
is, a permutation of the form $(i,\bar i)$, has type $A_1$ when
regarded as a parabolic Coxeter element, while we say that it has
combinatorial type $B_1$. (The reader should recall that, in the
classification of finite Coxeter groups, the type $B_1$ does not
occur, respectively, that sometimes $B_1$ is identified with $A_1$. 
Here, when we speak of ``combinatorial type," then we {\it do} distinguish
between $A_1$ and $B_1$. For example, the ``cycles"
$((1,2))=(1,2)\,(\bar 1,\bar 2)$ or 
$((\bar 1,2))=(\bar 1,2)\,(1,\bar 2)$ have combinatorial type $A_1$,
whereas the cycles $[1]=(1,\bar 1)$ or $[2]=(2,\bar 2)$ have
combinatorial type $B_1$.)

As Coxeter element for $W(B_n)$, 
we choose 
$$c=(1,2,\dots,n,\bar1,\bar2,\dots,\bar
n)=[1,2,\dots,n].$$ 
Now, given {\it combinatorial\/} types
$T_1,T_2,\dots,T_d$, each of which being a product of $A_k$'s and
$B_k$'s, $k=1,2,\dots$, 
the {\it combinatorial decomposition number} 
$N_{B_n}^{\text {comb}}(T_1,T_2,\dots,T_d)$ is defined
as the number of minimal products $c_1c_2\cdots c_d$
less than or equal to $c$ in absolute order, where 
``minimal" has the same meaning as above,
such that for $i=1,2,\dots,d$ the {\it combinatorial\/} type of $c_i$ is $T_i$.
Because of \eqref{eq:Umord}, the combinatorial decomposition numbers
$N_{B_n}^{\text {comb}}(T_1,T_2,\dots,T_d)$ satisfy also \eqref{Aa}
and \eqref{Ab}.

The reflection group $W(D_n)$ can be realised as the subgroup
of the group of all permutations $\pi$ of $\{1,2,\dots,n,\bar 1,\bar
2,\dots,\bar n\}$
satisfying \eqref{eq:-} and the property that an even number of 
elements from $\{1,2,\dots,n\}$ is mapped to an element of negative
sign. (Here, the elements $1,2,\dots,n$ are considered to have sign
$+$, while the elements $\bar1,\bar2,\dots,\bar n$ are considered to have
sign $-$.) Since $W(D_n)$ is a subgroup of $W(B_n)$, and since the
above realisation of $W(D_n)$
is contained as a subset in the realisation of $W(B_n)$
that we just described,
any $\pi\in W(D_n)$ can be decomposed as in \eqref{eq:cycB},
where, for $i=1,2,\dots,d$, 
$\ka_i$ is either a type $A$ or a type $B$ cycle. 
Requiring that $\pi$ is in the subgroup $W(D_n)$ of $W(B_n)$ is equivalent
to requiring that there is an even number of type $B$ cycles in the
decomposition \eqref{eq:cycB}.
Again, the decomposition \eqref{eq:cycB} for $\pi\in W(D_n)$
is unique up to a reordering of the
$\ka_i$'s.

As Coxeter element, we choose 
$$c=(1,2,\dots,n-1,\bar1,\bar2,\dots,
\overline{n-1})\,(n,\bar n)=[1,2,\dots,n-1]\,[n].$$
We shall be entirely concerned with elements $\pi$ of $W(D_n)$ which
are less than or equal to $c$. It is not difficult to see (and it is
shown in \cite[Sec.~3]{AtReAA}) that the unique factorisation of 
any such element $\pi$ has either $0$ or $2$ type $B$ cycles, and in
the latter case one of the type $B$ cycles is $[n]=(n,\bar n)$. 
In this latter case, in abuse of terminology, we call the product of
these two type $B$ cycles,
$
[a_1,a_2,\dots,a_{k-1}]\,[n]
$
say, a ``cycle" {\it of combinatorial type} $D_k$. More generally, we
shall say for any product of two disjoint type $B$ cycles of the form 
\begin{equation} \label{eq:dcycle}
[a_1,a_2,\dots,a_{k-1}]\,[a_k]
\end{equation}
that it is a ``cycle" of {\it combinatorial type} $D_k$.
The reader should observe that,
when regarded as parabolic Coxeter element, for $k\ge4$ an element
of the form \eqref{eq:dcycle} has type
$D_{k}$. However, if $k=3$, it has type $A_3$ when
regarded as parabolic Coxeter element, while we say that it has
combinatorial type $D_3$, and, if $k=2$, it has type $A_1^2$ when
regarded as parabolic Coxeter element, while we say that it has
combinatorial type $D_2$. (The reader should recall that, in the
classification of finite Coxeter groups, the types $D_3$ and $D_2$ do not
occur, respectively, that sometimes $D_3$ is identified with $A_3$,
$D_2$ being identified with $A_1^2$. 
Here, when we speak of ``combinatorial type," then we {\it do} distinguish
between $D_3$ and $A_3$, and between $D_2$ and $A_1^2$.)

Now, given {\it combinatorial\/} types
$T_1,T_2,\dots,T_d$, each of which being a product of $A_k$'s and
$D_k$'s, $k=1,2,\dots$, 
the {\it combinatorial decomposition number} 
$N_{B_n}^{\text {comb}}(T_1,T_2,\dots,T_d)$ is defined
as the number of minimal products $c_1c_2\cdots c_d$
less than or equal to $c$ in absolute order, where 
``minimal" has the same meaning as above,
such that for $i=1,2,\dots,d$ the {\it combinatorial\/} type of $c_i$ is $T_i$.
Because of \eqref{eq:Umord}, the combinatorial decomposition numbers
$N_{D_n}^{\text {comb}}(T_1,T_2,\dots,T_d)$ satisfy also \eqref{Aa}
and \eqref{Ab}.

\section{Auxiliary results}
\label{sec:1a} 

In our computations in the proof of Theorem~\ref{thm:3}, leading to the determination
of the decomposition numbers of type $D_n$,
we need to apply
the Lagrange--Good inversion formula \cite{GoodAA}
(see also \cite[Sec.~5]{KratAC} and the references cited therein).
We recall it here for the convenience of the reader. In doing
so, we use standard multi-index notation. Name\-ly, given a positive
integer $d$, and vectors $\mathbf z=(z_1,z_2,\dots,z_d)$ and
$\mathbf n=(n_1,n_2,\dots,n_d)$, we write
$\mathbf z^{\mathbf n}$ for $z_1^{n_1}z_2^{n_2}\cdots z_d^{n_d}$.
Furthermore, in abuse of notation, given a formal power series $f$ in $d$
variables, $f(\mathbf z)$ stands for
$f(z_1,z_2,\dots,z_d)$. Moreover, given $d$ formal power series
$f_1,f_2,\dots,f_d$ in $d$ variables, $\mathbf f^{\mathbf n}(\mathbf
z)$ is short for
$$f_1^{n_1}(z_1,z_2,\dots,z_d)f_2^{n_2}(z_1,z_2,\dots,z_d)
\cdots f_d^{n_d}(z_1,z_2,\dots,z_d).$$
Finally, if $\mathbf
m=(m_1,m_2,\dots,m_d)$ is another vector, then $\mathbf m+\mathbf n$
is short for $(m_1+n_1,m_2+n_2,\dots,m_d+n_d)$. Notation such as $\mathbf
m-\mathbf n$ has to be interpreted in a similar way.

\begin{theorem}[\sc Lagrange--Good inversion] \label{thm:LG}
Let $d$ be a positive integer, and let $f_1(\mathbf z),f_2(\mathbf z),
\dots,f_d(\mathbf z)$ be formal
power series in $\mathbf z=(z_1,z_2,\dots,z_d)$ with the property
that, for all $i$,
$f_i(\mathbf z)$ is of the form $z_i/\ph_i(\mathbf z)$ for some formal
power series $\ph_i(\mathbf z)$ with $\ph_i(0,0,\dots,0)\ne0$. Then, if we
expand a formal power series $g(\mathbf z)$ in terms of powers of the
$f_i(\mathbf z),$
\begin{equation} \label{eq:LG}
g(\mathbf z)=\sum _{\mathbf n} ^{}\ga_{\mathbf
n}\mathbf f^{\mathbf n}(\mathbf z),
\end{equation}
the coefficients $\ga_{\mathbf n}$ are given by
$$
\ga_{\mathbf n}=\coef{\mathbf z^{-\boldsymbol e}}g(\mathbf z)\mathbf
f^{-\mathbf n-\boldsymbol e}(\mathbf z)\det_{1\le i,j\le d}
\(\frac {\partial f_i} {\partial z_j}(\mathbf z)\),$$
where $\mathbf e=(1,1,\dots,1),$ 
where the sum in \eqref{eq:LG} 
runs over all $d$-tuples $\mathbf n$ of non-negative
integers, and where $\coef{\mathbf z^{\mathbf m}}h(\mathbf z)$
denotes the coefficient of $\mathbf z^{\mathbf m}$ in the formal
Laurent series $h(\mathbf z)$.
\end{theorem}

Next, we prove a determinant lemma and a corollary,
both of which will also be used in the proof of
Theorem~\ref{thm:3}.

\begin{lemma} \label{lem:1}
Let $d$ be a positive integer, and let
$X_1,X_2,\dots,X_d,Y_2,Y_3,\dots,Y_d$ be indeterminates. 
Then
\begin{equation} \label{eq:4}
\det_{1\le i,j\le d}\(\left\{\begin{matrix}
\displaystyle 1-\chi(1\ne j)\frac {Y_j} {X_1},&i=1\\
\displaystyle 1-\chi(i\ne j)\frac {Y_i} {X_i},&i\ge2
\end{matrix}\right\}\)=
\frac {\displaystyle 
\(\sum _{i=1} ^{d}X_i-\sum _{i=2} ^{d}Y_i\)Y_2Y_3\cdots
Y_d}
{X_1X_2\cdots X_d},
\end{equation}
where $\chi(\mathcal S)=1$ if $\mathcal S$ is
true and $\chi(\mathcal S)=0$ otherwise.
\end{lemma}
\begin{proof} By using multilinearity in the rows, we rewrite the
determinant on the left-hand side of \eqref{eq:4} as
$$
\frac {1} {X_1X_2\cdots X_d}
\det_{1\le i,j\le d}\(\left\{\begin{matrix}
X_1-\chi(1\ne j) {Y_j} ,&i=1\\
X_i-\chi(i\ne j) {Y_i} ,&i\ge2
\end{matrix}\right\}\).
$$
Next, we subtract the first column from all other columns.
As a result, we obtain the determinant
$$
\frac {1} {X_1X_2\cdots X_d}\det_{1\le i,j\le d}\(
\left\{\begin{matrix}
X_1 ,&i=j=1\\
-Y_j ,&i=1\text { and }j\ge2\\
X_i-Y_i ,&i\ge2\text { and }j=1\\
\chi(i=j)Y_i ,&i,j\ge2
\end{matrix}\right\}\).
$$
Now we add rows $2,3,\dots,d$ to the first row.
After that, our determinant becomes lower triangular, with the entry in the
first row and column equal to $\sum _{i=1} ^{d}X_i-
\sum _{i=2} ^{d}Y_i$, and
the diagonal entry in row $i$, $i\ge2$, equal to $Y_i$. 
Hence, we obtain the claimed result.
\end{proof}

\begin{corollary} \label{cor:1}
Let $d$ and $r$ be positive integers, $1\le r\le d,$ and let
$X_1,X_2,\dots,X_d,$ $Y$ and $Z$ be indeterminates. 
Then, with notation as in Lemma~{\em\ref{lem:1},}
\begin{multline} \label{eq:5}
\det_{1\le i,j\le d}\(\left\{\begin{matrix}
1-\chi(r\ne j)\frac {Z} {X_r},&i=r\\
1-\chi(i\ne j)\frac {Y} {X_i},&i\ne r
\end{matrix}\right\}\)\\=
\frac {Y^{d-2}\big(
Z\sum _{i=1} ^{d}X_i+(Y-Z)X_r-(d-1)YZ\big)}
{X_1X_2\cdots X_d}.
\end{multline}
\end{corollary}
\begin{proof}
We write the diagonal entry in the $r$-th row 
of the determinant in \eqref{eq:5} as 
$$1=\frac {X_r+Y-Z} {X_r}-\frac {Y-Z} {X_r},$$
and then use linearity of the determinant in the $r$-th row to
decompose the determinant as
$$\frac {X_r+Y-Z} {X_r}D_1-\frac {Y-Z} {X_r}D_2,$$
where $D_1$ is the determinant in \eqref{eq:4} with
$X_r$ replaced by $X_r+Y-Z$, and with $Y_i=Y$ for all
$i$, and where $D_2$ is the determinant in \eqref{eq:4} with $d$
replaced by $d-1$, with $Y_i=Y$ for all $i$, and with $X_i$ replaced
by $X_{i-1}$ for $i=r+1,r+2,\dots,d$. Hence, using Lemma~\ref{lem:1},
we deduce that the determinant in \eqref{eq:5} is equal to
\begin{equation*}
\frac {Y^{d-1}\big(\sum _{i=1} ^{d}X_i+Y-Z-(d-1)Y\big)}
{X_1X_2\cdots X_d}
-
\frac {(Y-Z)Y^{d-2}\big(\sum _{i=1} ^{d}X_i-X_r-(d-2)Y\big)}
{X_1X_2\cdots X_d}.
\end{equation*}
Little simplification then leads to \eqref{eq:5}.
\end{proof}

We end this section with a summation lemma, which we shall
need in Sections~\ref{sec:3} and \ref{sec:4} in order to compute the
$B_n$, respectively $D_n$, decomposition numbers of arbitrary rank 
from those of full rank, and in Section~\ref{sec:7} to derive
enumerative results for (generalised) non-crossing partitions 
from our formulae for the decomposition numbers.

\begin{lemma} \label{lem:binsum}
Let $M$ and $r$ be non-negative integers. Then
\begin{equation} \label{eq:binsum}
\sum _{m_1+2m_2+\dots+rm_r=r} ^{}\binom M{m_1,m_2,\dots,m_r}=
\binom {M+r-1}r,
\end{equation}
where the multinomial coefficient
is defined by
$$\binom {M}{m_1,m_2,\dots,m_r}=\frac {M!} {m_1!\,m_2!\cdots
m_r!\,(M-m_1-m_2-\dots-m_r)!}.$$
\end{lemma}
\begin{proof}
The identity results directly by comparing coefficients of 
$z^r$ on both sides of the identity
$$(1+z+z^2+z^3+\cdots)^M=(1-z)^{-M}.$$ 
\end{proof}

\section{Decomposition numbers for type $A$}
\label{sec:2} 

As was already pointed out in \cite[Sec.~10]{KratCF},
the decomposition numbers for type $A_n$ have already been computed
by Goulden and Jackson in \cite[Theorem~3.2]{GoJaAS}, 
albeit using a somewhat different language. (The condition on the sum 
$l(\al_1)+l(\al_2)+\dots+l(\al_m)$ is misstated throughout the latter
paper. It should be replaced by
$l(\al_1)+l(\al_2)+\dots+l(\al_m)=(m-1)n+1$.) In our terminology,
their result reads as follows.

\begin{theorem} \label{thm:1}
Let $T_1,T_2,\dots,T_d$ be types
with $\rk T_1+\rk T_2+\dots+\rk T_d=n,$
where 
$$T_i=A_1^{m_1^{(i)}}*A_2^{m_2^{(i)}}*\dots*A_n^{m_n^{(i)}},\quad 
i=1,2,\dots,d.$$ 
Then
\begin{equation} \label{eq:1}
N_{A_n}(T_1,T_2,\dots,T_d)=(n+1)^{d-1}\prod _{i=1} ^{d}
\frac {1} {n-\rk T_i+1}\binom {n-\rk T_i+1}{m_1^{(i)},m_2^{(i)},\dots,
m_n^{(i)}},
\end{equation}
where the multinomial coefficient is defined as in
Lemma~{\em\ref{lem:binsum}}.
\end{theorem}

Here we have used Stembridge's \cite{StemAZ} notation for the
decomposition of types into a product of irreducibles; for example,
the equation $T=A_2^3*A_5$ means that the root system of type $T$ 
decomposes into the orthogonal product of $3$ copies of root systems
of type $A_2$ and one copy of the root system of type $A_5$.

It was shown in \cite[Theorem~10]{KratCF} that, upon applying the
summation formula in Lemma~\ref{lem:binsum} to the result in
Theorem~\ref{thm:1} in a suitable manner, one obtains a compact formula
for {\it all\/} type $A_n$ decomposition numbers.

\begin{theorem} \label{thm:20}
Let the types $T_1,T_2,\dots,T_d$ be given,
where 
$$T_i=A_1^{m_1^{(i)}}*A_2^{m_2^{(i)}}*\dots*A_n^{m_n^{(i)}},\quad 
i=1,2,\dots,d.$$ 
Then
\begin{multline} \label{eq:100}
N_{A_n}(T_1,T_2,\dots,T_d)=(n+1)^{d-1}
\binom {n+1}{\rk T_1+\rk T_2+\dots+\rk T_d+1}\\
\times
\prod _{i=1} ^{d}
\frac {1} {n-\rk T_i+1}\binom {n-\rk T_i+1}{m_1^{(i)},m_2^{(i)},\dots,
m_n^{(i)}},
\end{multline}
where the multinomial coefficient is defined as in
Lemma~{\em\ref{lem:binsum}}.
All other decomposition numbers $N_{A_n}(T_1,T_2,\dots,T_d)$
are zero.
\end{theorem}

\section{Decomposition numbers for type $B$}
\label{sec:3} 

In this section we compute the decomposition numbers in type $B_n$.
We show that one can extract the corresponding formulae
from results of B\'ona, Bousquet, Labelle and Leroux \cite{BoBLAA} on
the enumeration of certain planar maps, which they call $m$-ary cacti.
While reading the statement of the theorem, the reader should
recall from Section~\ref{sec:1} the distinction between group-theoretic and
combinatorial decomposition numbers.

\begin{theorem} \label{thm:2}
{\em(i)} If $T_1,T_2,\dots,T_d$ are types
with $\rk T_1+\rk T_2+\dots+\rk T_d=n,$
where 
$$T_i=A_1^{m_1^{(i)}}*A_2^{m_2^{(i)}}*\dots*A_n^{m_n^{(i)}},\quad 
i=1,2,\dots,j-1,j+1,\dots,d,$$ 
and
$$T_j=B_\al*A_1^{m_1^{(j)}}*A_2^{m_2^{(j)}}*\dots*A_n^{m_n^{(j)}},
$$ 
for some $\al\ge1,$
then
\begin{equation} \label{eq:2}
N_{B_n}^{\text {comb}}(T_1,T_2,\dots,T_d)=n^{d-1}
\binom {n-\rk T_j}{m_1^{(j)},m_2^{(j)},\dots,
m_n^{(j)}}
\underset{i\ne j}{\prod _{i=1} ^{d}}
\frac {1} {n-\rk T_i}\binom {n-\rk T_i}{m_1^{(i)},m_2^{(i)},\dots,
m_n^{(i)}},
\end{equation}
where the multinomial coefficient is defined as in
Lemma~{\em\ref{lem:binsum}}. 
For $\al\ge2,$ the number 
$N_{B_n}(T_1,T_2,\dots, T_d)$ is given by the same formula.

\smallskip
{\em(ii)} If $T_1,T_2,\dots,T_d$ are types
with $\rk T_1+\rk T_2+\dots+\rk T_d=n,$
where 
$$T_i=A_1^{m_1^{(i)}}*A_2^{m_2^{(i)}}*\dots*A_n^{m_n^{(i)}},\quad 
i=1,2,\dots,d,$$ 
then
\begin{equation} \label{eq:3}
N_{B_n}(T_1,T_2,\dots,T_d)=n^{d-1}\Bigg(\prod _{i=1} ^{d}
\frac {1} {n-\rk T_i}\binom {n-\rk T_i}{m_1^{(i)},m_2^{(i)},\dots,
m_n^{(i)}}\Bigg)
\sum _{j=1} ^{d}\frac{m^{(j)}_1(n-\rk T_j)} {m^{(j)}_0+1},
\end{equation}
where $m^{(j)}_0=n-\rk T_j-\sum _{s=1} ^{n}m^{(j)}_s$.

\smallskip
{\em(iii)} All other decomposition numbers $N_{B_n}(T_1,T_2,\dots,T_d)$
and $N_{B_n}^{\text {comb}}(T_1,T_2,\dots,T_d)$ with
$\rk T_1+\rk T_2+\dots+\rk T_d=n$ are zero.
\end{theorem}

\begin{proof}
Determining the decomposition numbers 
$$N_{B_n}(T_1,T_2,\dots,T_d)=
N_{B_n}(T_d,\dots,T_2,T_1)$$
(recall \eqref{Aa}), respectively
$$N_{B_n}^{\text {comb}}(T_1,T_2,\dots,T_d)=
N_{B_n}^{\text {comb}}(T_d,\dots,T_2,T_1),$$ 
amounts to counting all possible factorisations 
\begin{equation} \label{eq:facB}
[1,2,\dots,n]=\si_d\cdots\si_2\si_1,
\end{equation}
where $\si_i$ has type $T_i$ as a parabolic Coxeter element, 
respectively has combinatorial type $T_i$.
The reader should observe that
the factorisation \eqref{eq:facB} is {\it minimal}, in the
sense that
$$n=\ell_T\big([1,2,\dots,n]\big)=
\ell_T(\si_1)+\ell_T(\si_2)+\dots+\ell_T(\si_d),$$
since $\ell_T(\si_i)=\rk T_i$, and since, by
our assumption, the sum of the ranks of the $T_i$'s equals $n$.
A further observation is that, in a factorisation \eqref{eq:facB},
there must be at least one factor $\si_i$ which contains a type $B$
cycle in its
(type $B$) disjoint cycle decomposition, because the
sign of $[1,2,\dots,n]$ as an element of the group $S_{2n}$ of all
permutations of $\{1,2,\dots,n,\bar 1,\bar 2,\dots,\bar n\}$ is
$-1$, while the sign of any type $A$ cycle is $+1$. 

\medskip
We first prove Claim~(iii). Let us assume, by contradiction, that
there is a minimal decomposition \eqref{eq:facB}
in which, altogether, we find {\it at least two} type $B$ cycles in the
(type $B$) disjoint cycle decompositions of the $\si_i$'s. In that case,
\eqref{eq:facB} has the form
\begin{equation} \label{eq:facu}
[1,2,\dots,n]=u_1\ka_1u_2\ka_2u_3,
\end{equation}
where $\ka_1$ and $\ka_2$ are two type $B$ cycles, and
$u_1,u_2,u_3$ are the factors in between. Moreover, the factorisation
\eqref{eq:facu} is minimal, meaning that 
\begin{equation} \label{eq:minB}
n=\ell_T(u_1)+\ell_T(\ka_1)+\ell_T(u_2)+\ell_T(\ka_2)+\ell_T(u_3).
\end{equation}
We may rewrite \eqref{eq:facu} as
$$
[1,2,\dots,n]=\ka_1\ka_2(\ka_2^{-1}\ka_1^{-1}u_1\ka_1\ka_2)
(\ka_2^{-1}u_2\ka_2)u_3,
$$
or, setting $u_1'=\ka_2^{-1}\ka_1^{-1}u_1\ka_1\ka_2$
and $u_2'=\ka_2^{-1}u_2\ka_2$, as
\begin{equation} \label{eq:facu1}
[1,2,\dots,n]=\ka_1\ka_2u_1'u_2'u_3.
\end{equation}
This factorisation is still minimal since $u_1'$ is conjugate to $u_1$
and $u_2'$ is conjugate to $u_2$. At this point, we observe that
$\ka_1$ must be a cycle of the form \eqref{eq:bcycle}
with $a_1<a_2<\dots<a_k<
\overline{a_1}<\overline{a_2}<\dots<\overline{a_k}$ in the order
$1<2<\dots<n<\bar 1<\bar2<\dots<\bar n$, because
otherwise $\ka_1\not\le_T [1,2,\dots,n]$, which would contradict
\eqref{eq:facu1}. A similar argument applies to $\ka_2$.
Now, if $\ka_1$ and $\ka_2$ are not
disjoint, then it is easy to see that
$\ell_T(\ka_1\ka_2)<\ell_T(\ka_1)+\ell_T(\ka_2)$, whence
\begin{align*}
n&=\ell_T([1,2,\dots,n])\\
&=\ell_T(\ka_1\ka_2u_1'u_2'u_3)\\
&\le\ell_T(\ka_1\ka_2)+\ell_T(u_1')+\ell_T(u_2')+\ell_T(u_3)\\
&\le\ell_T(\ka_1\ka_2)+\ell_T(u_1)+\ell_T(u_2)+\ell_T(u_3)\\
&<\ell_T(\ka_1)+\ell_T(\ka_2)+\ell_T(u_1)+\ell_T(u_2)+\ell_T(u_3),
\end{align*}
a contradiction to \eqref{eq:minB}. If, on the other hand, $\ka_1$
and $\ka_2$ are disjoint, then we can find 
$i,j\in\{1,2,\dots,n,\bar 1,\bar2,\dots,\bar n\}$, such that
$i<j<\ka_1(i)<\ka_2(j)$ (in the above order of
$\{1,2,\dots,n,\bar 1,\bar2,\dots,\bar n\}$). 
In other words, if we represent $\ka_1$ and
$\ka_2$ in the obvious way in a cyclic diagram 
(cf.\ \cite[Sec.~2]{ReivAG}), then they cross each
other. However, in that case we have
\begin{equation*} 
\ka_1\ka_2\not\le_T [1,2,\dots,n],
\end{equation*}
contradicting the fact that \eqref{eq:facu1} is a minimal
factorisation. (This is one of the consequences of Biane's
group-theoretic characterisation \cite[Theorem~1]{BianAA} of
non-crossing partitions.)

\begin{figure}
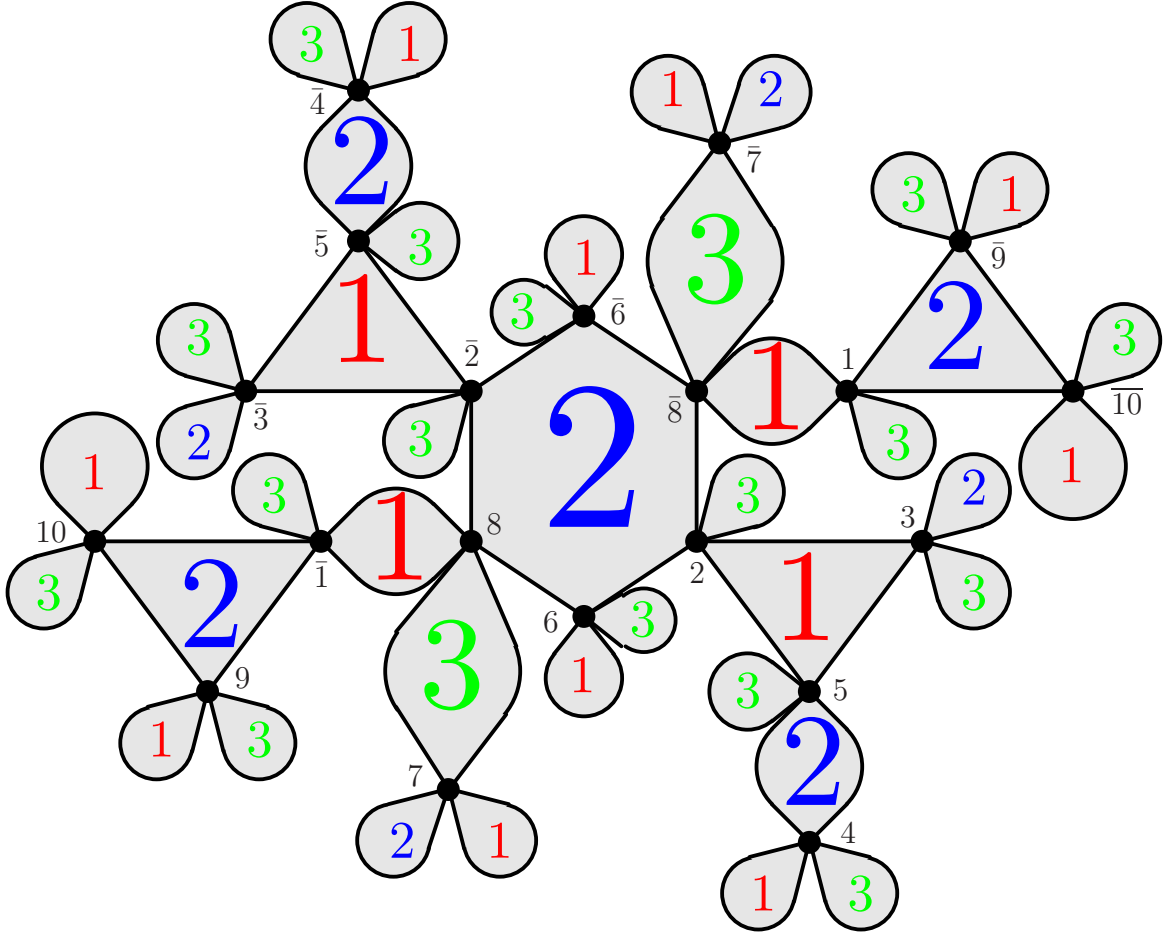

\centertexdraw{
  \drawdim truecm  \linewd.05
  \move(-1.5 -1)
  \rlvec(1.5 -1)
  \rlvec(1.5 1)
  \rlvec(0 2)
  \rlvec(-1.5 1)
  \rlvec(-1.5 -1)
  \rlvec(0 -2)
  \lfill f:.9
  \htext(-.6 -.8){\Blau{\Riesig2}}

  \move(1.5 -1)
  \rlvec(3 0)
  \rlvec(-1.5 -2)
  \rlvec(-1.5 2)
  \lfill f:.9
  \htext(2.5 -2.3){\Rot{\riesig1}}

  \move(-1.5 1)
  \rlvec(-3 0)
  \rlvec(1.5 2)
  \rlvec(1.5 -2)
  \lfill f:.9
  \htext(-3.4 1.4){\Rot{\riesig1}}

  \move(3 -4)
  \fcir f:.9 r:.7
  \larc r:.7 sd:135 ed:225
  \move(3 -4)
  \larc r:.7 sd:-45 ed:45
  \move(2.5 -3.5)
  \rlvec (.5 .5)
  \rlvec (.5 -.5)
  \lfill f:.9
  \move(2.5 -4.5)
  \rlvec (.5 -.5)
  \rlvec (.5 .5)
  \lfill f:.9
  \htext(2.6 -4.5){\riesig\Blau{2}}

  \move(-3 4)
  \fcir f:.9 r:.7
  \larc r:.7 sd:135 ed:225
  \move(-3 4)
  \larc r:.7 sd:-45 ed:45
  \move(-2.5 3.5)
  \rlvec (-.5 -.5)
  \rlvec (-.5 .5)
  \lfill f:.9
  \move(-2.5 4.5)
  \rlvec (-.5 .5)
  \rlvec (-.5 -.5)
  \lfill f:.9
  \htext(-3.4 3.5){\riesig\Blau{2}}

  \move(2.5 1)
  \fcir f:.9 r:.7
  \larc r:.7 sd:45 ed:135
  \move(2.5 1)
  \larc r:.7 sd:225 ed:-45
  \move(2 .5)
  \rlvec (-.5 .5)
  \rlvec (.5 .5)
  \lfill f:.9
  \move(3 .5)
  \rlvec (.5 .5)
  \rlvec (-.5 .5)
  \lfill f:.9
  \htext(2.1 .5){\riesig\Rot{1}}

  \move(-2.5 -1)
  \fcir f:.9 r:.7
  \larc r:.7 sd:45 ed:135
  \move(-2.5 -1)
  \larc r:.7 sd:225 ed:-45
  \move(-2 -.5)
  \rlvec (.5 -.5)
  \rlvec (-.5 -.5)
  \lfill f:.9
  \move(-3 -.5)
  \rlvec (-.5 -.5)
  \rlvec (.5 -.5)
  \lfill f:.9
  \htext(-2.9 -1.5){\riesig\Rot{1}}

  \move(2.18 -.34)
  \fcir f:.9 r:.48
  \larc r:.48 sd:-80 ed:170
  \move(2.3 -.8)
  \rlvec (-.8 -.2)
  \rlvec (.2 .8)
  \lfill f:.9
  \htext(2 -.6){\grosz\Gruen{3}}

  \move(-2.18 .34)
  \fcir f:.9 r:.48
  \larc r:.48 sd:100 ed:-10
  \move(-2.3 .8)
  \rlvec (.8 .2)
  \rlvec (-.2 -.8)
  \lfill f:.9
  \htext(-2.35 0.15){\grosz\Gruen{3}}

  \move(3.5 1)
  \rlvec(3 0)
  \rlvec(-1.5 2)
  \rlvec(-1.5 -2)
  \lfill f:.9
  \htext(4.5 1.3){\Blau{\riesig2}}

  \move(-3.5 -1)
  \rlvec(-3 0)
  \rlvec(1.5 -2)
  \rlvec(1.5 2)
  \lfill f:.9
  \htext(-5.4 -2.4){\Blau{\riesig2}}

  \move(4.18 .32)
  \fcir f:.9 r:.48
  \larc r:.48 sd:190 ed:80
  \move(4.3 .8)
  \rlvec (-.8 .2)
  \rlvec (.2 -.8)
  \lfill f:.9
  \htext(4 .1){\grosz\Gruen{3}}

  \move(-4.18 -.32)
  \fcir f:.9 r:.48
  \larc r:.48 sd:10 ed:260
  \move(-4.3 -.8)
  \rlvec (.8 -.2)
  \rlvec (-.2 .8)
  \lfill f:.9
  \htext(-4.3 -.6){\grosz\Gruen{3}}

  \move(-1.745 -2.72891)
  \fcir f:.9 r:.9
  \larc r:.9 sd:145 ed:215
  \larc r:.9 sd:-40 ed:30
  \move(-2.5 -2.2)
  \rlvec (1 1.2)
  \rlvec (0.6 -1.4)
  \lfill f:.9
  \move(-2.5 -3.2)
  \rlvec (.70 -1.1)
  \rlvec (0.76 1)
  \lfill f:.9
  \htext(-2.2 -3.2){\riesig\Gruen{3}}

  \move(1.745 2.72891)
  \fcir f:.9 r:.9
  \larc r:.9 sd:140 ed:210
  \larc r:.9 sd:-37 ed:30
  \move(2.5 2.2)
  \rlvec (-1 -1.2)
  \rlvec (-0.6 1.4)
  \lfill f:.9
  \move(2.5 3.2)
  \rlvec (-.70 1.1)
  \rlvec (-0.76 -1)
  \lfill f:.9
  \htext(1.3 2.2){\riesig\Gruen{3}}

  \move(2.18 -3)
  \fcir f:.9 r:.51
  \larc r:.51 sd:45 ed:-45
  \move(2.5 -2.6)
  \rlvec (0.5 -0.4)
  \rlvec (-0.5 -0.4)
  \lfill f:.9
  \htext(2.0 -3.2){\grosz\Gruen{3}}

  \move(-2.18 3)
  \fcir f:.9 r:.51
  \larc r:.51 sd:225 ed:135
  \move(-2.5 2.6)
  \rlvec (-0.5 0.4)
  \rlvec (0.5 0.4)
  \lfill f:.9
  \htext(-2.35 2.7){\grosz\Gruen{3}}

  \move(2.32 -5.68)
  \fcir f:.9 r:.48
  \larc r:.48 sd:100 ed:-10
  \move(2.2 -5.2)
  \rlvec (.8 .2)
  \rlvec (-.2 -.8)
  \lfill f:.9
  \htext(2.2 -5.9){\grosz\Rot{1}}

  \move(3.68 -5.68)
  \fcir f:.9 r:.48
  \larc r:.48 sd:190 ed:80
  \move(3.8 -5.2)
  \rlvec (-.8 .2)
  \rlvec (.2 -.8)
  \lfill f:.9
  \htext(3.5 -5.9){\grosz\Gruen{3}}

  \move(-2.32 5.68)
  \fcir f:.9 r:.48
  \larc r:.48 sd:280 ed:170
  \move(-2.2 5.2)
  \rlvec (-.8 -.2)
  \rlvec (.2 .8)
  \lfill f:.9
  \htext(-2.5 5.4){\grosz\Rot{1}}

  \move(-3.68 5.68)
  \fcir f:.9 r:.48
  \larc r:.48 sd:10 ed:260
  \move(-3.8 5.2)
  \rlvec (.8 -.2)
  \rlvec (-.2 .8)
  \lfill f:.9
  \htext(-3.8 5.4){\grosz\Gruen{3}}

  \move(0 -2.82)
  \fcir f:.9 r:.51
  \larc r:.51 sd:135 ed:45
  \move(-0.4 -2.5)
  \rlvec (0.4 0.5)
  \rlvec (0.4 -0.5)
  \lfill f:.9
  \htext(-.2 -3){\grosz\Rot{1}}

  \move(0 2.82)
  \fcir f:.9 r:.51
  \larc r:.51 sd:-45 ed:225
  \move(0.4 2.5)
  \rlvec (-0.4 -0.5)
  \rlvec (-0.4 0.5)
  \lfill f:.9
  \htext(-.15 2.56){\grosz\Rot{1}}

  \move(.8 -2.05)
  \fcir f:.9 r:.42
  \larc r:.42 sd:235 ed:130
  \move(0.5 -2.4)
  \rlvec (-0.5 0.4)
  \rlvec (0.5 0.3)
  \lfill f:.9
  \htext(.6 -2.25){\grosz\Gruen{3}}

  \move(-.8 2.05)
  \fcir f:.9 r:.43
  \larc r:.43 sd:52 ed:310
  \move(-0.5 2.4)
  \rlvec (0.5 -0.4)
  \rlvec (-0.5 -0.3)
  \lfill f:.9
  \htext(-1 1.85){\grosz\Gruen{3}}

  \move(-1.12 -4.98)
  \fcir f:.9 r:.48
  \larc r:.48 sd:190 ed:80
  \move(-1 -4.5)
  \rlvec (-.8 .2)
  \rlvec (.2 -.8)
  \lfill f:.9
  \htext(-1.3 -5.2){\grosz\Rot{1}}

  \move(-2.53 -4.98)
  \fcir f:.9 r:.48
  \larc r:.48 sd:100 ed:-10
  \move(-2.6 -4.5)
  \rlvec (.8 .2)
  \rlvec (-.26 -.8)
  \lfill f:.9
  \htext(-2.6 -5.2){\grosz\Blau{2}}

  \move(1.12 4.98)
  \fcir f:.9 r:.48
  \larc r:.48 sd:10 ed:260
  \move(1 4.5)
  \rlvec (.8 -.2)
  \rlvec (-.2 .8)
  \lfill f:.9
  \htext(1 4.8){\grosz\Rot{1}}

  \move(2.53 4.98)
  \fcir f:.9 r:.48
  \larc r:.48 sd:280 ed:170
  \move(2.6 4.5)
  \rlvec (-.8 -.2)
  \rlvec (.26 .8)
  \lfill f:.9
  \htext(2.3 4.8){\grosz\Blau{2}}

  \move(5.68 3.68)
  \fcir f:.9 r:.48
  \larc r:.48 sd:280 ed:170
  \move(5.8 3.2)
  \rlvec (-.8 -.2)
  \rlvec (.2 .8)
  \lfill f:.9
  \htext(5.5 3.4){\grosz\Rot{1}}

  \move(4.32 3.68)
  \fcir f:.9 r:.48
  \larc r:.48 sd:10 ed:260
  \move(4.2 3.2)
  \rlvec (.8 -.2)
  \rlvec (-.2 .8)
  \lfill f:.9
  \htext(4.2 3.4){\grosz\Gruen{3}}

  \move(-5.68 -3.68)
  \fcir f:.9 r:.48
  \larc r:.48 sd:100 ed:-10
  \move(-5.8 -3.2)
  \rlvec (.8 .2)
  \rlvec (-.2 -.8)
  \lfill f:.9
  \htext(-5.8 -3.85){\grosz\Rot{1}}

  \move(-4.32 -3.68)
  \fcir f:.9 r:.48
  \larc r:.48 sd:190 ed:80
  \move(-4.2 -3.2)
  \rlvec (-.8 .2)
  \rlvec (.2 -.8)
  \lfill f:.9
  \htext(-4.5 -3.85){\grosz\Gruen{3}}

  \move(6.5 0)
  \fcir f:.9 r:.7
  \larc r:.7 sd:135 ed:45
  \move(6 .5)
  \rlvec (0.5 0.5)
  \rlvec (0.5 -0.5)
  \lfill f:.9
  \htext(6.3 -.2){\grosz\Rot{1}}

  \move(7.18 1.68)
  \fcir f:.9 r:.48
  \larc r:.48 sd:-80 ed:170
  \move(7.3 1.2)
  \rlvec (-.8 -.2)
  \rlvec (.2 .8)
  \lfill f:.9
  \htext(7 1.5){\grosz\Gruen{3}}

  \move(-6.5 0)
  \fcir f:.9 r:.7
  \larc r:.7 sd:-45 ed:225
  \move(-6 -.5)
  \rlvec (-0.5 -0.5)
  \rlvec (-0.5 0.5)
  \lfill f:.9
  \htext(-6.7 -.3){\grosz\Rot{1}}

  \move(-7.18 -1.68)
  \fcir f:.9 r:.48
  \larc r:.48 sd:100 ed:-10
  \move(-7.3 -1.2)
  \rlvec (.8 .2)
  \rlvec (-.2 -.8)
  \lfill f:.9
  \htext(-7.3 -1.9){\grosz\Gruen{3}}

  \move(-5.18 .32)
  \fcir f:.9 r:.48
  \larc r:.48 sd:100 ed:-10
  \move(-5.3 .8)
  \rlvec (.8 .2)
  \rlvec (-.2 -.8)
  \lfill f:.9
  \htext(-5.3 .1){\grosz\Blau{2}}

  \move(-5.18 1.68)
  \fcir f:.9 r:.48
  \larc r:.48 sd:10 ed:-100
  \move(-5.3 1.2)
  \rlvec (.8 -.2)
  \rlvec (-.2 .8)
  \lfill f:.9
  \htext(-5.3 1.5){\grosz\Gruen{3}}

  \move(5.18 -.32)
  \fcir f:.9 r:.48
  \larc r:.48 sd:280 ed:170
  \move(5.3 -.8)
  \rlvec (-.8 -.2)
  \rlvec (.2 .8)
  \lfill f:.9
  \htext(5 -.5){\grosz\Blau{2}}

  \move(5.18 -1.68)
  \fcir f:.9 r:.48
  \larc r:.48 sd:190 ed:80
  \move(5.3 -1.2)
  \rlvec (-.8 .2)
  \rlvec (.2 -.8)
  \lfill f:.9
  \htext(5 -1.9){\grosz\Gruen{3}}

  \Ringerl(-1.5 -1)
  \Ringerl(1.5 -1)
  \Ringerl(1.5 1)
  \Ringerl(-1.5 1)
  \Ringerl(0 2)
  \Ringerl(0 -2)

  \Ringerl(4.5 -1)
  \Ringerl(-4.5 1)

  \Ringerl(3 -3)
  \Ringerl(3 -5)
  \Ringerl(-3 3)
  \Ringerl(-3 5)

  \Ringerl(3.5 1)
  \Ringerl(6.5 1)
  \Ringerl(5 3)
  \Ringerl(-3.5 -1)
  \Ringerl(-6.5 -1)
  \Ringerl(-5 -3)

  \Ringerl(1.8 4.3)
  \Ringerl(-1.8 -4.3)

  \htext(3.4 1.3){$1$}
  \htext(-3.6 -1.6){$\bar 1$}

  \htext(1.4 -1.55){$2$}
  \htext(-1.6 1.35){$\bar 2$}

  \htext(4.2 -.8){$3$}
  \htext(-4.4 .5){$\bar 3$}

  \htext(3.4 -5){$4$}
  \htext(-3.65 4.7){$\bar 4$}

  \htext(3.3 -3.1){$5$}
  \htext(-3.6 2.8){$\bar 5$}

  \htext(-.55 -2.2){$6$}
  \htext(.35 1.9){$\bar 6$}

  \htext(-2.35 -4.25){$7$}
  \htext(2.15 3.93){$\bar 7$}

  \htext(-1.3 -.9){$8$}
  \htext(1.1 .6){$\bar 8$}

  \htext(-4.65 -2.95){$9$}
  \htext(5.4 2.65){$\bar 9$}

  \htext(-7.3 -1){$10$}
  \htext(7.0 .7){$\overline{10}$}
}
\caption{The $3$-cactus corresponding to the factorisation
\eqref{eq:fac0B}}
\label{fig:3}
\end{figure}

\medskip
We turn now to Claims~(i) and (ii). In what follows, we shall show
that the formulae \eqref{eq:2} and \eqref{eq:3} follow from results
of B\'ona, Bousquet, Labelle and Leroux \cite{BoBLAA} on the
enumeration of $m$-ary cacti with a rotational symmetry. In order to
explain this, we must first define a bijection between minimal factorisations
\eqref{eq:facB} and certain planar maps. 
By a {\it map}, we mean a connected 
graph embedded in the plane
such that edges do not intersect except in vertices. The maps which
are of relevance here are maps in which faces different from the
outer face intersect only in vertices, and
are coloured with colours from $\{1,2,\dots,d\}$. 
Such maps will be referred to
as {\it $d$-cacti} from now on.\footnote{\label{foot:1}%
We warn the reader that our
terminology deviates from the one in \cite{BoBLAA,GoJaAS}. We follow
loosely the conventions in \cite{IrviAA}. To be precise, our
$d$-cacti in which the rotator around every vertex is 
$(1,2,\dots,d)^{\text{O}}$
are dual to the coloured $d$-cacti in \cite{GoJaAS}, respectively
$d$-ary cacti in \cite{BoBLAA}, in 
the following sense: one is obtained from the other by ``interchanging" the
roles of vertices and faces, that is, given a $d$-cactus in our
sense, one obtains a $d$-cactus in the sense of Goulden and Jackson
by shrinking faces to vertices and blowing up vertices of degree $\de$
to faces with $\de$ vertices, keeping the incidence relations between
faces and vertices. Another minor difference is that colours are
arranged in counter-clockwise order in \cite{BoBLAA,GoJaAS}, 
while we arrange colours in clockwise order.}
Examples of $3$-cacti
can be found in Figures~\ref{fig:3} and \ref{fig:4}. 
In the figures, the faces different from the outer face are the shaded
ones. Their colours are indicated by the numbers $1$, $2$,
respectively $3$, placed in the centre of the faces. Figure~\ref{fig:3}
shows a $3$-cactus in which the vertices are labelled, while
Figure~\ref{fig:4} shows one in which the vertices are not labelled.
(That one of the vertices in Figure~\ref{fig:4} 
is marked by a bold dot should be ignored
for the moment.)

In what follows, we need the concept of the {\it rotator} around a
vertex $v$ in a $d$-cactus, which, by definition, is the cyclic list of
colours of faces encountered in a clockwise journey around $v$.
If, while travelling around $v$, we encounter the colours
$b_1,b_2,\dots,b_k$, in this order, then we will write
$(b_1,b_2,\dots,b_k)^{\text{O}}$ for the rotator, meaning that
$(b_1,b_2,\dots,b_k)^{\text{O}}= (b_2,\dots,b_k,b_1)^{\text{O}}$, etc.
For example, the rotator of {\it all\/} the vertices in the map
in Figure~\ref{fig:3} is $(1,2,3)^{\text{O}}$.

We illustrate the bijection between minimal factorisations
\eqref{eq:facB} and $d$-cacti
with an example. Take $n=10$ and $d=3$,
and consider the factorisation
\begin{equation} \label{eq:fac0B}
[1,2,\dots,10]=\si_3\si_2\si_1,
\end{equation}
where $\si_3=((7,8))$,
$\si_2=[2,6,8]\,((1,\bar 9,\overline{10}))\,((4,5))$,
and $\si_1=((1,\bar 8))\,((2,3,5))$.
For each cycle $(a_1,a_2,\dots,a_k)$ (sic!) of $\si_i$, we create a $k$-gon
coloured $i$, and label its vertices $a_1,a_2,\dots,a_k$ in clockwise
order. (The warning ``sic!" is there to avoid misunderstandings:
for each type $A$ ``cycle" $((b_1,b_2,\dots,b_k))$ we create {\it two}
$k$-gons, the vertices of one being labelled $b_1,b_2,\dots,b_k$, and
the vertices of the other being labelled 
$\overline{b_1},\overline{b_2},\dots,\overline{b_k}$, while 
for each type $B$ ``cycle" $[b_1,b_2,\dots,b_k]$ we create {\it one}
$2k$-gon with vertices labelled $b_1,b_2,\dots,b_k,
\overline{b_1},\overline{b_2},\dots,\overline{b_k}$.) 
We glue these polygons
into a $d$-cactus, the faces of which are these
polygons plus the outer face, by identifying equally labelled vertices such that
the rotator of each vertex is $(1,2,\dots,d)$.
Figure~\ref{fig:3} shows the outcome of this procedure for the
factorisation \eqref{eq:fac0B}.

The fact that the result of the procedure can be realised as 
a $d$-cactus follows from Euler's 
formula. Namely, the number of faces corresponding to the polygons is 
$1+2{\sum _{i=1} ^{d}}\sum _{k=0} ^{n}m^{(i)}_k$ (the $1$
coming from the polygon corresponding to the type $B$ cycle), 
the number of edges
is $2\al+2\sum _{i=1} ^{d}\sum _{k=0} ^{n}m^{(i)}_k(k+1)$, and the number
of vertices is $2n$. Hence, if we include the outer face, the number
of vertices minus the number of edges plus the number of faces is
\begin{align} \notag
2n-
2\al-2\sum _{i=1} ^{d}\sum _{k=0} ^{n}m^{(i)}_k(k+1)+{}
&2\sum _{i=1} ^{d}\sum _{k=0} ^{n}m^{(i)}_k+2\\
\notag
&=
2n+2-2\al-
2\sum _{i=1} ^{d}\sum _{k=0} ^{n}k\cdot m^{(i)}_k\\
\notag
&=2n+2-2\rk T_1-2\rk T_2-\dots-2\rk T_d\\
&=2,
\label{eq:EulerB}
\end{align}
according to our assumption concerning the sum of the ranks of the
types $T_i$.

\begin{figure}
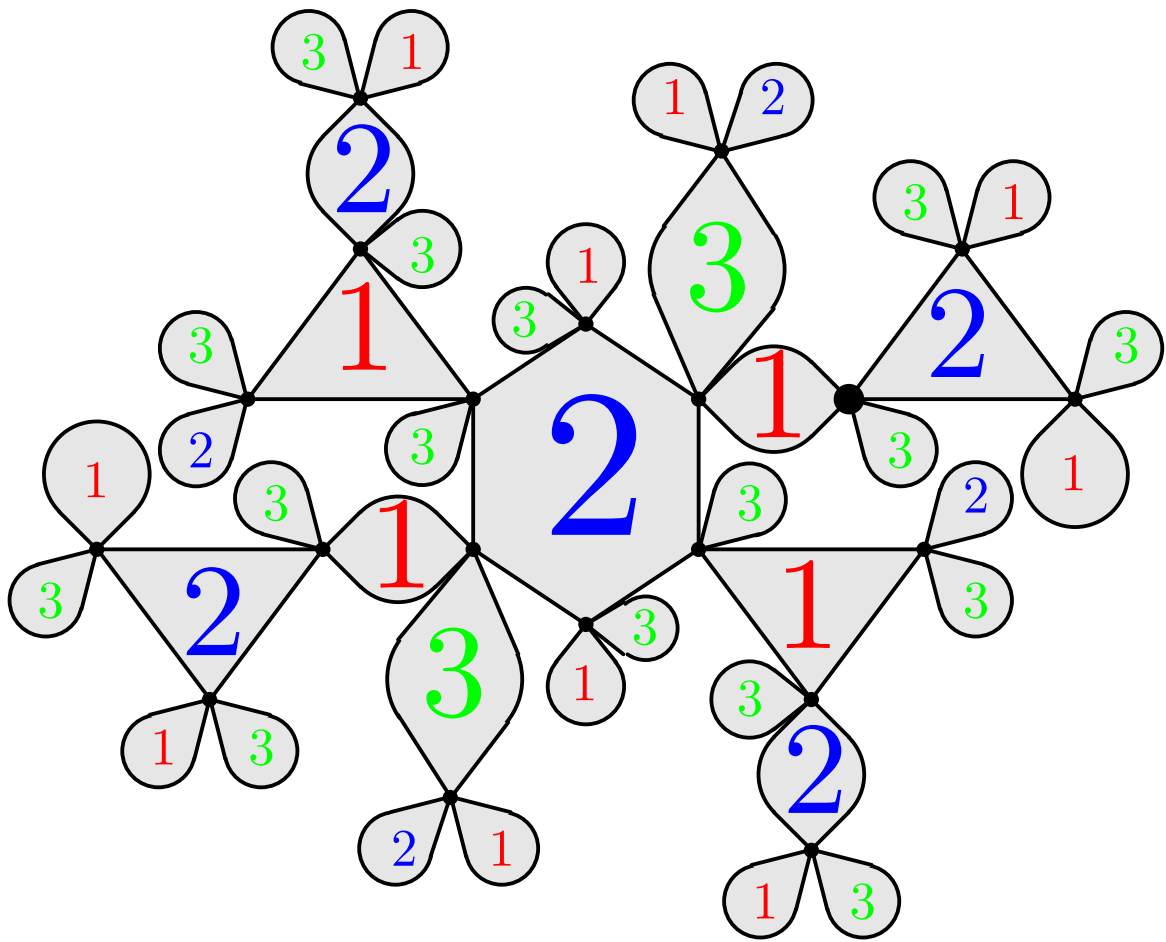

\centertexdraw{
  \drawdim truecm  \linewd.05
  \move(-1.5 -1)
  \rlvec(1.5 -1)
  \rlvec(1.5 1)
  \rlvec(0 2)
  \rlvec(-1.5 1)
  \rlvec(-1.5 -1)
  \rlvec(0 -2)
  \lfill f:.9
  \htext(-.6 -.8){\Blau{\Riesig2}}

  \move(1.5 -1)
  \rlvec(3 0)
  \rlvec(-1.5 -2)
  \rlvec(-1.5 2)
  \lfill f:.9
  \htext(2.5 -2.3){\Rot{\riesig1}}

  \move(-1.5 1)
  \rlvec(-3 0)
  \rlvec(1.5 2)
  \rlvec(1.5 -2)
  \lfill f:.9
  \htext(-3.4 1.4){\Rot{\riesig1}}

  \move(3 -4)
  \fcir f:.9 r:.7
  \larc r:.7 sd:135 ed:225
  \move(3 -4)
  \larc r:.7 sd:-45 ed:45
  \move(2.5 -3.5)
  \rlvec (.5 .5)
  \rlvec (.5 -.5)
  \lfill f:.9
  \move(2.5 -4.5)
  \rlvec (.5 -.5)
  \rlvec (.5 .5)
  \lfill f:.9
  \htext(2.6 -4.5){\riesig\Blau{2}}

  \move(-3 4)
  \fcir f:.9 r:.7
  \larc r:.7 sd:135 ed:225
  \move(-3 4)
  \larc r:.7 sd:-45 ed:45
  \move(-2.5 3.5)
  \rlvec (-.5 -.5)
  \rlvec (-.5 .5)
  \lfill f:.9
  \move(-2.5 4.5)
  \rlvec (-.5 .5)
  \rlvec (-.5 -.5)
  \lfill f:.9
  \htext(-3.4 3.5){\riesig\Blau{2}}

  \move(2.5 1)
  \fcir f:.9 r:.7
  \larc r:.7 sd:45 ed:135
  \move(2.5 1)
  \larc r:.7 sd:225 ed:-45
  \move(2 .5)
  \rlvec (-.5 .5)
  \rlvec (.5 .5)
  \lfill f:.9
  \move(3 .5)
  \rlvec (.5 .5)
  \rlvec (-.5 .5)
  \lfill f:.9
  \htext(2.1 .5){\riesig\Rot{1}}

  \move(-2.5 -1)
  \fcir f:.9 r:.7
  \larc r:.7 sd:45 ed:135
  \move(-2.5 -1)
  \larc r:.7 sd:225 ed:-45
  \move(-2 -.5)
  \rlvec (.5 -.5)
  \rlvec (-.5 -.5)
  \lfill f:.9
  \move(-3 -.5)
  \rlvec (-.5 -.5)
  \rlvec (.5 -.5)
  \lfill f:.9
  \htext(-2.9 -1.5){\riesig\Rot{1}}

  \move(2.18 -.34)
  \fcir f:.9 r:.48
  \larc r:.48 sd:-80 ed:170
  \move(2.3 -.8)
  \rlvec (-.8 -.2)
  \rlvec (.2 .8)
  \lfill f:.9
  \htext(2 -.6){\grosz\Gruen{3}}

  \move(-2.18 .34)
  \fcir f:.9 r:.48
  \larc r:.48 sd:100 ed:-10
  \move(-2.3 .8)
  \rlvec (.8 .2)
  \rlvec (-.2 -.8)
  \lfill f:.9
  \htext(-2.35 0.15){\grosz\Gruen{3}}

  \move(3.5 1)
  \rlvec(3 0)
  \rlvec(-1.5 2)
  \rlvec(-1.5 -2)
  \lfill f:.9
  \htext(4.5 1.3){\Blau{\riesig2}}

  \move(-3.5 -1)
  \rlvec(-3 0)
  \rlvec(1.5 -2)
  \rlvec(1.5 2)
  \lfill f:.9
  \htext(-5.4 -2.4){\Blau{\riesig2}}

  \move(4.18 .32)
  \fcir f:.9 r:.48
  \larc r:.48 sd:190 ed:80
  \move(4.3 .8)
  \rlvec (-.8 .2)
  \rlvec (.2 -.8)
  \lfill f:.9
  \htext(4 .1){\grosz\Gruen{3}}

  \move(-4.18 -.32)
  \fcir f:.9 r:.48
  \larc r:.48 sd:10 ed:260
  \move(-4.3 -.8)
  \rlvec (.8 -.2)
  \rlvec (-.2 .8)
  \lfill f:.9
  \htext(-4.3 -.6){\grosz\Gruen{3}}

  \move(-1.745 -2.72891)
  \fcir f:.9 r:.9
  \larc r:.9 sd:145 ed:215
  \larc r:.9 sd:-40 ed:30
  \move(-2.5 -2.2)
  \rlvec (1 1.2)
  \rlvec (0.6 -1.4)
  \lfill f:.9
  \move(-2.5 -3.2)
  \rlvec (.70 -1.1)
  \rlvec (0.76 1)
  \lfill f:.9
  \htext(-2.2 -3.2){\riesig\Gruen{3}}

  \move(1.745 2.72891)
  \fcir f:.9 r:.9
  \larc r:.9 sd:140 ed:210
  \larc r:.9 sd:-37 ed:30
  \move(2.5 2.2)
  \rlvec (-1 -1.2)
  \rlvec (-0.6 1.4)
  \lfill f:.9
  \move(2.5 3.2)
  \rlvec (-.70 1.1)
  \rlvec (-0.76 -1)
  \lfill f:.9
  \htext(1.3 2.2){\riesig\Gruen{3}}

  \move(2.18 -3)
  \fcir f:.9 r:.51
  \larc r:.51 sd:45 ed:-45
  \move(2.5 -2.6)
  \rlvec (0.5 -0.4)
  \rlvec (-0.5 -0.4)
  \lfill f:.9
  \htext(2.0 -3.2){\grosz\Gruen{3}}

  \move(-2.18 3)
  \fcir f:.9 r:.51
  \larc r:.51 sd:225 ed:135
  \move(-2.5 2.6)
  \rlvec (-0.5 0.4)
  \rlvec (0.5 0.4)
  \lfill f:.9
  \htext(-2.35 2.7){\grosz\Gruen{3}}

  \move(2.32 -5.68)
  \fcir f:.9 r:.48
  \larc r:.48 sd:100 ed:-10
  \move(2.2 -5.2)
  \rlvec (.8 .2)
  \rlvec (-.2 -.8)
  \lfill f:.9
  \htext(2.2 -5.9){\grosz\Rot{1}}

  \move(3.68 -5.68)
  \fcir f:.9 r:.48
  \larc r:.48 sd:190 ed:80
  \move(3.8 -5.2)
  \rlvec (-.8 .2)
  \rlvec (.2 -.8)
  \lfill f:.9
  \htext(3.5 -5.9){\grosz\Gruen{3}}

  \move(-2.32 5.68)
  \fcir f:.9 r:.48
  \larc r:.48 sd:280 ed:170
  \move(-2.2 5.2)
  \rlvec (-.8 -.2)
  \rlvec (.2 .8)
  \lfill f:.9
  \htext(-2.5 5.4){\grosz\Rot{1}}

  \move(-3.68 5.68)
  \fcir f:.9 r:.48
  \larc r:.48 sd:10 ed:260
  \move(-3.8 5.2)
  \rlvec (.8 -.2)
  \rlvec (-.2 .8)
  \lfill f:.9
  \htext(-3.8 5.4){\grosz\Gruen{3}}

  \move(0 -2.82)
  \fcir f:.9 r:.51
  \larc r:.51 sd:135 ed:45
  \move(-0.4 -2.5)
  \rlvec (0.4 0.5)
  \rlvec (0.4 -0.5)
  \lfill f:.9
  \htext(-.2 -3){\grosz\Rot{1}}

  \move(0 2.82)
  \fcir f:.9 r:.51
  \larc r:.51 sd:-45 ed:225
  \move(0.4 2.5)
  \rlvec (-0.4 -0.5)
  \rlvec (-0.4 0.5)
  \lfill f:.9
  \htext(-.15 2.56){\grosz\Rot{1}}

  \move(.8 -2.05)
  \fcir f:.9 r:.42
  \larc r:.42 sd:235 ed:130
  \move(0.5 -2.4)
  \rlvec (-0.5 0.4)
  \rlvec (0.5 0.3)
  \lfill f:.9
  \htext(.6 -2.25){\grosz\Gruen{3}}

  \move(-.8 2.05)
  \fcir f:.9 r:.43
  \larc r:.43 sd:52 ed:310
  \move(-0.5 2.4)
  \rlvec (0.5 -0.4)
  \rlvec (-0.5 -0.3)
  \lfill f:.9
  \htext(-1 1.85){\grosz\Gruen{3}}

  \move(-1.12 -4.98)
  \fcir f:.9 r:.48
  \larc r:.48 sd:190 ed:80
  \move(-1 -4.5)
  \rlvec (-.8 .2)
  \rlvec (.2 -.8)
  \lfill f:.9
  \htext(-1.3 -5.2){\grosz\Rot{1}}

  \move(-2.53 -4.98)
  \fcir f:.9 r:.48
  \larc r:.48 sd:100 ed:-10
  \move(-2.6 -4.5)
  \rlvec (.8 .2)
  \rlvec (-.26 -.8)
  \lfill f:.9
  \htext(-2.6 -5.2){\grosz\Blau{2}}

  \move(1.12 4.98)
  \fcir f:.9 r:.48
  \larc r:.48 sd:10 ed:260
  \move(1 4.5)
  \rlvec (.8 -.2)
  \rlvec (-.2 .8)
  \lfill f:.9
  \htext(1 4.8){\grosz\Rot{1}}

  \move(2.53 4.98)
  \fcir f:.9 r:.48
  \larc r:.48 sd:280 ed:170
  \move(2.6 4.5)
  \rlvec (-.8 -.2)
  \rlvec (.26 .8)
  \lfill f:.9
  \htext(2.3 4.8){\grosz\Blau{2}}

  \move(5.68 3.68)
  \fcir f:.9 r:.48
  \larc r:.48 sd:280 ed:170
  \move(5.8 3.2)
  \rlvec (-.8 -.2)
  \rlvec (.2 .8)
  \lfill f:.9
  \htext(5.5 3.4){\grosz\Rot{1}}

  \move(4.32 3.68)
  \fcir f:.9 r:.48
  \larc r:.48 sd:10 ed:260
  \move(4.2 3.2)
  \rlvec (.8 -.2)
  \rlvec (-.2 .8)
  \lfill f:.9
  \htext(4.2 3.4){\grosz\Gruen{3}}

  \move(-5.68 -3.68)
  \fcir f:.9 r:.48
  \larc r:.48 sd:100 ed:-10
  \move(-5.8 -3.2)
  \rlvec (.8 .2)
  \rlvec (-.2 -.8)
  \lfill f:.9
  \htext(-5.8 -3.85){\grosz\Rot{1}}

  \move(-4.32 -3.68)
  \fcir f:.9 r:.48
  \larc r:.48 sd:190 ed:80
  \move(-4.2 -3.2)
  \rlvec (-.8 .2)
  \rlvec (.2 -.8)
  \lfill f:.9
  \htext(-4.5 -3.85){\grosz\Gruen{3}}

  \move(6.5 0)
  \fcir f:.9 r:.7
  \larc r:.7 sd:135 ed:45
  \move(6 .5)
  \rlvec (0.5 0.5)
  \rlvec (0.5 -0.5)
  \lfill f:.9
  \htext(6.3 -.2){\grosz\Rot{1}}

  \move(7.18 1.68)
  \fcir f:.9 r:.48
  \larc r:.48 sd:-80 ed:170
  \move(7.3 1.2)
  \rlvec (-.8 -.2)
  \rlvec (.2 .8)
  \lfill f:.9
  \htext(7 1.5){\grosz\Gruen{3}}

  \move(-6.5 0)
  \fcir f:.9 r:.7
  \larc r:.7 sd:-45 ed:225
  \move(-6 -.5)
  \rlvec (-0.5 -0.5)
  \rlvec (-0.5 0.5)
  \lfill f:.9
  \htext(-6.7 -.3){\grosz\Rot{1}}

  \move(-7.18 -1.68)
  \fcir f:.9 r:.48
  \larc r:.48 sd:100 ed:-10
  \move(-7.3 -1.2)
  \rlvec (.8 .2)
  \rlvec (-.2 -.8)
  \lfill f:.9
  \htext(-7.3 -1.9){\grosz\Gruen{3}}

  \move(-5.18 .32)
  \fcir f:.9 r:.48
  \larc r:.48 sd:100 ed:-10
  \move(-5.3 .8)
  \rlvec (.8 .2)
  \rlvec (-.2 -.8)
  \lfill f:.9
  \htext(-5.3 .1){\grosz\Blau{2}}

  \move(-5.18 1.68)
  \fcir f:.9 r:.48
  \larc r:.48 sd:10 ed:-100
  \move(-5.3 1.2)
  \rlvec (.8 -.2)
  \rlvec (-.2 .8)
  \lfill f:.9
  \htext(-5.3 1.5){\grosz\Gruen{3}}

  \move(5.18 -.32)
  \fcir f:.9 r:.48
  \larc r:.48 sd:280 ed:170
  \move(5.3 -.8)
  \rlvec (-.8 -.2)
  \rlvec (.2 .8)
  \lfill f:.9
  \htext(5 -.5){\grosz\Blau{2}}

  \move(5.18 -1.68)
  \fcir f:.9 r:.48
  \larc r:.48 sd:190 ed:80
  \move(5.3 -1.2)
  \rlvec (-.8 .2)
  \rlvec (.2 -.8)
  \lfill f:.9
  \htext(5 -1.9){\grosz\Gruen{3}}

  \Mark(3.5 1)

  \ringerl(-1.5 -1)
  \ringerl(1.5 -1)
  \ringerl(1.5 1)
  \ringerl(-1.5 1)
  \ringerl(0 2)
  \ringerl(0 -2)

  \ringerl(4.5 -1)
  \ringerl(-4.5 1)

  \ringerl(3 -3)
  \ringerl(3 -5)
  \ringerl(-3 3)
  \ringerl(-3 5)

  \ringerl(6.5 1)
  \ringerl(5 3)
  \ringerl(-3.5 -1)
  \ringerl(-6.5 -1)
  \ringerl(-5 -3)

  \ringerl(1.8 4.3)
  \ringerl(-1.8 -4.3)

}
\caption{A rotation-symmetric $3$-cactus with a marked vertex}
\label{fig:4}
\end{figure}

We may further simplify this geometric representation of 
a minimal factorisation \eqref{eq:facB} by deleting all vertex labels and
marking the vertex which had label $1$. If this simplification is
applied to the $3$-cactus in Figure~\ref{fig:3}, we obtain the
$3$-cactus in Figure~\ref{fig:4}.
Indeed, the knowledge of which vertex carries label $1$ allows us
to reconstruct all other vertex labels as follows: starting from the
vertex labelled $1$, we travel clockwise along the boundary of the face
coloured $1$ until we reach the next vertex (that is, we traverse
only a single edge); from there, 
we travel clockwise along the boundary of the face
coloured $2$ until we reach the next vertex; etc., until we have
travelled along an edge bounding a face of colour $d$. The vertex
that we have reached must carry label $2$. Etc.
Clearly, if drawn appropriately into the plane,
a $d$-cactus resulting from an application of the above
procedure to a minimal factorisation \eqref{eq:facB} is 
symmetric with respect to a rotation by $180^\circ$, the centre of
the rotation being the centre of the regular $2\al$-gon corresponding
to the unique type $B$ cycle of $\si_j$; cf.\ Figure~\ref{fig:4}.
In what follows, we shall abbreviate this property as {\it
rotation-symmetric}.

In summary, under the assumptions of Claim~(i),
the decomposition number\break 
$N_{B_n}^{\text {comb}}(T_1,T_2,\dots,T_d)$, respectively, if
$\al\ge2$, the decomposition number $N_{B_n}(T_1,T_2,\dots,\break T_d)$ also,
equals the number of all rotation-symmetric 
$d$-cacti on $2n$ vertices in which one
vertex is marked and all vertices have rotator $(1,2,\dots,d)^{\text
{O}}$, with exactly $m^{(i)}_k$ pairs of faces of colour
$i$ having $k+1$ vertices, arranged symmetrically
around a central face of colour $j$ with $2\al$ vertices.

Aside from the marking of one vertex, equivalent objects are counted
in \cite[Theorem~25]{BoBLAA}. In our language,
modulo the ``dualisation" described
in Footnote~\ref{foot:1}, and upon replacing $m$ by $d$, 
the objects which are counted in the cited theorem are $d$-cacti
in which all vertices have rotator $(1,2,\dots,d)^{\text
{O}}$, and which are invariant under a rotation (not necessarily by
$180^\circ$). To be precise, from the proof of \cite[(81)]{BoBLAA}
(not given in full detail in \cite{BoBLAA}) 
it can be extracted that the number of
$d$-cacti on $2n$ vertices, in which 
all vertices have rotator $(1,2,\dots,d)^{\text
{O}}$, which are invariant under a rotation by $(360/s)^\circ$, $s$
being {\it maximal\/} with this property,
and which have exactly $2 m^{(i)}_k$ faces of colour
$i$ having $k+1$ vertices arranged 
around a central face of colour $j$ with $2\al$ vertices, equals
\begin{multline} \label{eq:BoBLAA}
(2n)^{d-2}s
{\sum _{t} ^{}}{}^{\displaystyle\prime}\mu(t/s)
\binom {2(n-\rk T_j)/t}{2m_1^{(j)}/t,2m_2^{(j)}/t,\dots,
2m_n^{(j)}/t}\\
\cdot
\underset{i\ne j}{\prod _{i=1} ^{d}}
\frac {1} {2(n-\rk T_i)}\binom {2(n-\rk
T_i)/t}{2m_1^{(i)}/t,2m_2^{(i)}/t,\dots,
2m_n^{(i)}/t},
\end{multline}
where the sum extends over all $t$ with $s\mid t$, $t\mid 2\al$,
and $t\mid 2m_k^{(i)}$ for all $i=1,2,\dots,d$ and
$k=1,2,\dots,n$. Here,
$\mu(\cdot)$ is the M\"obius function from number
theory.\footnote{Formula~(81) in \cite{BoBLAA} does not distinguish
the colour or the size of the central face (that is,
in the language of \cite{BoBLAA}:
the colour or the degree of the central vertex),
therefore it is in fact a sum over all possible colours and
sizes, represented there by the summations over $i$ and $h$,
respectively.}
In presenting the formula in the above form, we have also used the
observation that, for all $i$ (including $i=j$~!), the number of type
$A$ cycles of $\si_i$ is $n-\rk T_i$.

As we said above, the $d$-cacti that we want to enumerate have one
marked vertex, whereas the $d$-cacti counted by \eqref{eq:BoBLAA}
have no marked vertex. However, given a $d$-cactus counted by
\eqref{eq:BoBLAA}, we have exactly $2n/s$ inequivalent ways of
marking a vertex. Hence, recalling that the $d$-cacti that we want to
count are invariant under a rotation by $180^\circ$, we must multiply
the expression \eqref{eq:BoBLAA} by $2n/s$, and then sum the result
over all {\it even} $s$. Since, by definition of the M\"obius
function, we have
$$\sum _{2\mid s\mid t} ^{}\mu(t/s)=
\sum _{s'\mid \frac {t} {2}} ^{}\mu(t/2s')=\begin{cases} 
1&\text {if $\frac {t} {2}=1$,}\\
0&\text {otherwise,}\end{cases}
$$
the result of this summation is exactly the right-hand side of \eqref{eq:2}.

\medskip
Finally, we prove Claim~(ii). From what we already know,
in a minimal factorisation \eqref{eq:facB} exactly one of the
factors on the right-hand side must contain a type $B$ cycle of length
$1$ in its (type $B$) disjoint cycle decomposition, 
$\si_j$ say. As a parabolic
Coxeter element, a type $B$ cycle of length $1$ has type $A_1$.
Since all considerations in the proof of Claim~(i) are also valid for
$\al=1$, we may use Formula~\eqref{eq:2} with $\al=1$, and with
$m^{(j)}_1$ replaced by $m^{(j)}_1-1$, to count the number of 
these factorisations, to obtain
$$n^{d-1}
\binom {n-\rk T_j}{m_1^{(j)}-1,m_2^{(j)},\dots,
m_n^{(j)}}
\underset{i\ne j}{\prod _{i=1} ^{d}}
\frac {1} {n-\rk T_i}\binom {n-\rk T_i}{m_1^{(i)},m_2^{(i)},\dots,
m_n^{(i)}}.$$
This has to be summed over $j=1,2,\dots,d$. The result is exactly
\eqref{eq:3}.

\medskip
The proof of the theorem is now complete.
\end{proof}

Combining the previous theorem with the summation formula of 
Lemma~\ref{lem:binsum}, we can now derive compact formulae
for {\it all\/} type $B_n$ decomposition numbers.

\begin{theorem} \label{thm:21}
{\em(i)} 
Let the types $T_1,T_2,\dots,T_d$ be given,
where 
$$T_i=A_1^{m_1^{(i)}}*A_2^{m_2^{(i)}}*\dots*A_n^{m_n^{(i)}},\quad 
i=1,2,\dots,j-1,j+1,\dots,d,$$ 
and
$$T_j=B_\al*A_1^{m_1^{(j)}}*A_2^{m_2^{(j)}}*\dots*A_n^{m_n^{(j)}},
$$ 
for some $\al\ge1$.
Then
\begin{multline} \label{eq:21a}
N_{B_n}^{\text {comb}}(T_1,T_2,\dots,T_d)=n^{d-1}
\binom {n}{\rk T_1+\rk T_2+\dots+\rk T_d}\\
\times
\binom {n-\rk T_j}{m_1^{(j)},m_2^{(j)},\dots,
m_n^{(j)}}
\underset{i\ne j}{\prod _{i=1} ^{d}}
\frac {1} {n-\rk T_i}\binom {n-\rk T_i}{m_1^{(i)},m_2^{(i)},\dots,
m_n^{(i)}},
\end{multline}
where the multinomial coefficient is defined as in
Lemma~{\em\ref{lem:binsum}}.
For $\al\ge2,$ the number 
$N_{B_n}(T_1,T_2,\dots, T_d)$ is given by the same formula.

\smallskip
{\em(ii)} 
Let the types $T_1,T_2,\dots,T_d$ be given,
where 
$$T_i=A_1^{m_1^{(i)}}*A_2^{m_2^{(i)}}*\dots*A_n^{m_n^{(i)}},\quad 
i=1,2,\dots,d.$$ 
Then
\begin{multline} \label{eq:21b}
N_{B_n}^{\text {comb}}(T_1,T_2,\dots,T_d)\\
=n^{d}
\binom {n-1}{\rk T_1+\rk T_2+\dots+\rk T_d}
\Bigg(\prod _{i=1} ^{d}
\frac {1} {n-\rk T_i}\binom {n-\rk T_i}{m_1^{(i)},m_2^{(i)},\dots,
m_n^{(i)}}\Bigg),
\end{multline}
whereas
\begin{multline} \label{eq:21c}
N_{B_n}(T_1,T_2,\dots,T_d)\\
=n^{d-1}
\binom {n}{\rk T_1+\rk T_2+\dots+\rk T_d}
\Bigg(\prod _{i=1} ^{d}
\frac {1} {n-\rk T_i}\binom {n-\rk T_i}{m_1^{(i)},m_2^{(i)},\dots,
m_n^{(i)}}\Bigg)\\
\times
\(n-\rk T_1-\rk T_2-\dots-\rk T_d
+\sum _{j=1} ^{d}\frac{m^{(j)}_1(n-\rk T_j)} {m^{(j)}_0+1}\),
\end{multline}
with $m^{(j)}_0=n-\rk T_j-\sum _{s=1} ^{n}m^{(j)}_s$.

\smallskip
{\em(iii)} All other decomposition numbers $N_{B_n}(T_1,T_2,\dots,T_d)$
and $N_{B_n}^{\text {comb}}(T_1,T_2,\dots,T_d)$ are zero.
\end{theorem}

\begin{proof} 
If we write $r$ for $n-\rk T_1-\rk T_2-\dots-\rk T_d$,
then for $\Phi=B_n$ the relation \eqref{Ab} becomes
\begin{equation} \label{eq:relB}
N_{B_n}(T_1,T_2,\dots,T_d)=
{\sum _{T:\rk T=r}}
^{}N_{B_n}(T_1,T_2,\dots,T_d,T),
\end{equation}
with the same relation holding for $N_{B_n}^{\text {comb}}$ in place
of $N_{B_n}$.

In order to prove \eqref{eq:21a}, we
let $T=A_1^{m_1}*A_2^{m_2}*\dots*A_n^{m_n}$ and use
\eqref{eq:2} in \eqref{eq:relB}, to obtain
\begin{multline*}
N_{B_n}^{\text {comb}}(T_1,T_2,\dots,T_d)=\sum _{m_1+2m_2+\dots+nm_n=r} ^{}
n^{d}\frac {1} {n-r}\binom {n-r}{m_1,m_2,\dots,m_n}\\
\cdot
\binom {n-\rk T_j}{m_1^{(j)},m_2^{(j)},\dots,
m_n^{(j)}}
\underset{i\ne j}{\prod _{i=1} ^{d}}
\frac {1} {n-\rk T_i}\binom {n-\rk T_i}{m_1^{(i)},m_2^{(i)},\dots,
m_n^{(i)}}.
\end{multline*}
If we use \eqref{eq:binsum} with $M=n-r$, we arrive at our claim after
little simplification.

In order to prove \eqref{eq:21b}, we
let $T=B_\al*A_1^{m_1}*A_2^{m_2}*\dots*A_n^{m_n}$ in \eqref{eq:relB}.
The important point to be observed here is that, in contrast to the
previous argument, in the present case 
$T$ {\it must\/} have a factor
$B_\al$. Subsequently, use of \eqref{eq:2} in \eqref{eq:relB}
yields
\begin{multline} \label{eq:alsum}
N_{B_n}^{\text {comb}}(T_1,T_2,\dots,T_d)=\sum _{\al=1} ^{n}\sum
_{m_1+2m_2+\dots+nm_n=r-\al} ^{}
n^{d}\binom {n-r}{m_1,m_2,\dots,m_n}\\
\cdot
{\prod _{i=1} ^{d}}
\frac {1} {n-\rk T_i}\binom {n-\rk T_i}{m_1^{(i)},m_2^{(i)},\dots,
m_n^{(i)}}.
\end{multline}
Now we use \eqref{eq:binsum} with $r$ replaced by $r-\al$ and $M=n-r$, 
and subsequently the elementary summation formula
\begin{equation} \label{eq:rksum0}
\sum _{\al=1} ^{n}\binom {n-\al-1}{r-\al}=
\sum _{\al=1} ^{n}\binom {n-\al-1}{n-r-1}=\binom {n-1}{n-r}=\binom
{n-1}{r-1}.
\end{equation}
Then, after little rewriting, we arrive at our claim.

To establish \eqref{eq:21c}, we must recall that the {\it
group-theoretic} type $A_1$ does not distinguish between a type $A$
cycle $((i,j))=(i,j)\,(\bar i,\bar j)$ and a type $B$ cycle
$[i]=(i,\bar i)$. Hence, to obtain $N_{B_n}(T_1,T_2,\dots,T_d)$ in
the case that no $T_i$ contains a $B_\al$ for $\al\ge2$, we
must add the expression \eqref{eq:21b} and the expressions
\eqref{eq:21a} with $m_1^{(j)}$ replaced by $m_1^{(j)}-1$ 
over $j=1,2,\dots,d$. As is not difficult to see, this
sum is indeed equal to \eqref{eq:21c}.
\end{proof}

\section{Decomposition numbers for type $D$}
\label{sec:4} 

In this section we compute the decomposition numbers for the type $D_n$.
Theorem~\ref{thm:3} gives the formulae for the full rank
decomposition numbers, while Theorem~\ref{thm:22} presents the
implied formulae for the decomposition numbers of arbitrary rank.
To our knowledge, these are new results, which did not appear earlier
in the literature on map enumeration or on the connection coefficients
in the symmetric group or other Coxeter groups.
Nevertheless, the proof of Theorem~\ref{thm:3} is entirely in the
spirit of the fundamental paper \cite{GoJaAS}, in that the problem of
counting factorisations is translated into a problem of map
enumeration, which is then solved by a generating function approach
that requires the use of the Lagrange--Good formula for coefficient
extraction.

We begin with the result concerning the full rank decomposition
numbers in type $D_n$.
While reading the statement of the theorem below, the reader should again
recall from Section~\ref{sec:1} the distinction between group-theoretic and
combinatorial decomposition numbers.

\begin{theorem} \label{thm:3}
{\em(i)} If $T_1,T_2,\dots,T_d$ are types
with $\rk T_1+\rk T_2+\dots+\rk T_d=n,$
where 
$$T_i=A_1^{m_1^{(i)}}*A_2^{m_2^{(i)}}*\dots*A_n^{m_n^{(i)}},\quad 
i=1,2,\dots,j-1,j+1,\dots,d,$$ 
and
$$T_j=D_\al*A_1^{m_1^{(j)}}*A_2^{m_2^{(j)}}*\dots*A_n^{m_n^{(j)}},
$$ 
for some $\al\ge2,$
then
\begin{multline} \label{eq:6}
N_{D_n}^{\text {comb}}(T_1,T_2,\dots,T_d)=(n-1)^{d-1}
\binom {n-\rk T_j}{m_1^{(j)},m_2^{(j)},\dots,
m_n^{(j)}}\\
\times
\underset{i\ne j}{\prod _{i=1} ^{d}}
\frac {1} {n-\rk T_i-1}\binom {n-\rk T_i-1}{m_1^{(i)},m_2^{(i)},\dots,
m_n^{(i)}},
\end{multline}
where the multinomial coefficient is defined as in
Lemma~{\em\ref{lem:binsum}}.
For $\al\ge4,$ the number  
$N_{D_n}(T_1,T_2,\dots, T_d)$ is given by the same formula.

\smallskip
{\em(ii)} If $T_1,T_2,\dots,T_d$ are types
with $\rk T_1+\rk T_2+\dots+\rk T_d=n,$
where 
$$T_i=A_1^{m_1^{(i)}}*A_2^{m_2^{(i)}}*\dots*A_n^{m_n^{(i)}},\quad 
i=1,2,\dots,d,$$ 
then
\begin{multline} \label{eq:7comb}
N_{D_n}^{\text {comb}}(T_1,T_2,\dots,T_d)\\=(n-1)^{d-1}
\(2\sum _{j=1} ^{d}
\binom {n-\rk T_j}{m_1^{(j)},m_2^{(j)},\dots,m_n^{(j)}}
\underset{i\ne j} {\prod _{i=1} ^{d}}
\frac {1} {n-\rk T_i-1}\binom {n-\rk T_i-1}{m_1^{(i)},m_2^{(i)},\dots,
m_n^{(i)}}\right.\\
\left.
-2(d-1)(n-1)\underset{\vphantom{f}}{\prod _{i=1} ^{d}}
\frac {1} {n-\rk T_i-1}\binom {n-\rk T_i-1}{m_1^{(i)},m_2^{(i)},\dots,
m_n^{(i)}}\),
\end{multline}
while
\begin{multline} \label{eq:7}
N_{D_n}(T_1,T_2,\dots,T_d)\\=(n-1)^{d-1}
\(\sum _{j=1} ^{d}\Bigg(\underset{i\ne j}
{\prod _{i=1} ^{d}}
\frac {1} {n-\rk T_i-1}\binom {n-\rk T_i-1}{m_1^{(i)},m_2^{(i)},\dots,
m_n^{(i)}}\Bigg)
\Bigg(
2\binom {n-\rk T_j}{m_1^{(j)},m_2^{(j)},\dots,m_n^{(j)}}\right.\\
+\binom {n-\rk T_j}{m_1^{(j)},m_2^{(j)},m_3^{(j)}-1,m_4^{(j)},\dots,m_n^{(j)}}
+\binom {n-\rk T_j}{m_1^{(j)}-2,m_2^{(j)},\dots,m_n^{(j)}}
\Bigg)\\
\left.
-2(d-1)(n-1)\underset{\vphantom{f}}{\prod _{i=1} ^{d}}
\frac {1} {n-\rk T_i-1}\binom {n-\rk T_i-1}{m_1^{(i)},m_2^{(i)},\dots,
m_n^{(i)}}\).
\end{multline}

\smallskip
{\em(iii)} All other decomposition numbers $N_{D_n}(T_1,T_2,\dots,T_d)$
and $N_{D_n}^{\text {comb}}(T_1,T_2,\dots,T_d)$ with
$\rk T_1+\rk T_2+\dots+\rk T_d=n$ are zero.
\end{theorem}

\begin{remark}
These formulae must be correctly interpreted when $T_i$ contains
no $D_\al$ and $\rk T_i=n-1$.
In that case, because of $n-1=\rk
T_i=m_1^{(i)}+2m_2^{(i)}+\dots+nm_n^{(i)}$, there must be an $\ell$,
$1\le \ell\le n-1$, with $m_\ell^{(i)}\ge1$. 
We then interpret the term
$$\frac {1} {n-\rk T_i-1}\binom {n-\rk T_i-1}{m_1^{(i)},m_2^{(i)},\dots,
m_n^{(i)}}$$
as
$$\frac {1} {n-\rk T_i-1}\binom {n-\rk T_i-1}{m_1^{(i)},m_2^{(i)},\dots,
m_n^{(i)}}=\frac {1} {m_\ell^{(i)}}\binom {n-\rk T_i-2}{m_1^{(i)},
\dots,m_\ell^{(i)}-1,\dots,
m_n^{(i)}},$$
where the multinomial coefficient is zero whenever 
$$-1=n-\rk T_i-2<
m_1^{(i)}+\dots+(m_\ell^{(i)}-1)+\dots+m_n^{(i)},$$ 
except when all of
$m_1^{(i)},\dots,m_\ell^{(i)}-1,\dots,m_n^{(i)}$ are zero.
Explicitly, one must read
$$\frac {1} {n-\rk T_i-1}\binom {n-\rk T_i-1}{m_1^{(i)},m_2^{(i)},\dots,
m_n^{(i)}}=0$$
if $\rk T_i=n-1$ but $T_i\ne A_{n-1}$, and
$$\frac {1} {n-\rk T_i-1}\binom {n-\rk T_i-1}{m_1^{(i)},m_2^{(i)},\dots,
m_n^{(i)}}=1$$
if $T_i=A_{n-1}$.
\end{remark}

\begin{proof}[Proof of Theorem~{\em\ref{thm:3}}]
Determining the decomposition number 
$$N_{D_n}(T_1,T_2,\dots,T_d)=
N_{D_n}(T_d,\dots,T_2,T_1)$$ 
(recall \eqref{Aa}), respectively
$$N_{D_n}^{\text {comb}}(T_1,T_2,\dots,T_d)=
N_{D_n}^{\text {comb}}(T_d,\dots,T_2,T_1),$$ 
amounts to counting all possible factorisations 
\begin{equation} \label{eq:facD}
(1,2,\dots,n-1,\bar1,\bar2,\dots,
\overline{n-1})\,(n,\bar n)=\si_d\cdots\si_2\si_1,
\end{equation}
where $\si_i$ has type $T_i$ as a parabolic Coxeter element, 
respectively has combinatorial type $T_i$.
Here also, the factorisation \eqref{eq:facD} is {\it minimal} in the
sense that
$$n=\ell_T\big((1,2,\dots,n-1,\bar1,\bar2,\dots,
\overline{n-1})\,(n,\bar n)\big)=
\ell_T(\si_1)+\ell_T(\si_2)+\dots+\ell_T(\si_d),$$
since $\ell_T(\si_i)=\rk T_i$, and since, by
our assumption, the sum of the ranks of the $T_i$'s equals $n$.

\medskip
We first prove Claim~(iii). Let us assume, for contradiction, that
there is a minimal factorisation \eqref{eq:facD}, 
in which, altogether, we find {\it at least two} type $B$ cycles 
of length $\ge2$ in the
(type $B$) disjoint cycle decompositions of the $\si_i$'s. 
It can then be shown by arguments similar to those in the proof of
Claim~(iii) in Theorem~\ref{thm:2} that this leads to a contradiction.
Hence, ``at worst," we may find a type $B$ cycle of length $1$,
$(a,\bar a)$ say, and
another type $B$ cycle, $\ka$ say.
Both of them must be contained in the
disjoint cycle decomposition of one of the $\si_i$'s since all the
$\si_i$'s are elements of $W(D_n)$. Given that
$\ka$ has length $\al-1$, the product of both, $(a,\bar
a)\,\ka$, is of combinatorial type $D_\al$, $\al\ge2$, whereas, as
a parabolic Coxeter element, it is of type $D_\al$ only if $\al\ge4$.
If $\al=3$, then it is a
parabolic Coxeter element of type $A_3$, and if $\al=2$ it is of type
$A_1^2$. Thus, we are actually
in the cases to which Claims~(i) and (ii) apply.

\medskip
To prove Claim~(i), we continue this line of argument. By a
variation of the conjugation argument
\eqref{eq:facu}--\eqref{eq:facu1}, we may assume that these two type $B$
cycles are contained in $\si_d$, $\si_d=(a,\bar a)\,\ka\,\si_d'$ say,
where, as above, $(a,\bar a)$ is the type $B$
cycle of length $1$ and $\ka$ is the other type $B$ cycle, and where
$\si_d'$ is free of type $B$ cycles. 
In that case,
\eqref{eq:facD} takes the form
\begin{equation} \label{eq:facuD}
c=(1,2,\dots,n-1,\bar1,\bar2,\dots,
\overline{n-1})\,(n,\bar n)=(a,\bar a)\,\ka\,\si_d'\cdots\si_1.
\end{equation}
If $a\ne n$, $\ka\ne(n,\bar n)$, and if $\ka$ does not fix $n$,
then $(a,\bar a)\ka\not\le_Tc$, a contradiction. Likewise, 
if $a\ne n$, $\ka=[b_1,b_2,\dots,b_k]$ with
$n\notin\{b_1,b_2,\dots,b_k\}$, then
$(a,\bar a)\,\ka\not\le_T[1,2,\dots,n-1]$,
again a contradiction. Hence, we may assume that
$a=n$, whence $(a,\bar a)\,\ka=\ka\,(n,\bar n)$ forms a parabolic Coxeter
element of type $D_\al$, given that $\ka$ has length $\al-1$.
We are
then in the position to determine all possible factorisations of the
form \eqref{eq:facuD}, which reduces to
\begin{equation} \label{eq:reduD}
(1,2,\dots,n-1,\bar1,\bar2,\dots,
\overline{n-1})=[1,2,\dots,n-1]=\ka\si_d'\cdots\si_1.
\end{equation}
This is now a minimal type $B$ factorisation of the form \eqref{eq:facB} with
$n$ replaced by $n-1$. We may therefore use Formula~\eqref{eq:2} with
$n$ replaced by $n-1$, and with $\rk T_j$ replaced by $\rk T_j-1$.
These substitutions lead exactly to \eqref{eq:6}.

\medskip
Finally, we turn to Claim~(ii). 
First we discuss two degenerate cases which come from the
identifications $D_3\sim A_3$, respectively $D_2\sim A_1^2$, and
which only occur for $N_{D_n}(T_1,T_2,\dots,T_d)$ (but not for the
combinatorial decomposition numbers 
$N_{D_n}^{\text {comb}}(T_1,T_2,\break\dots,T_d)$).
It may happen that one of the factors in \eqref{eq:facD}, 
let us say, without loss of generality, $\si_d$,
contains a type $B$ cycle of length $1$ and one of length $2$
in its disjoint cycle decomposition;
that is, $\si_d$ may contain
$$(n,\bar n)\,[a,b]=(n,\bar n)\,(a,b,\bar a,\bar b)=
[a,b]\,[b,n]\,[b,\bar n].$$
As a parabolic Coxeter element, this is of type $A_3$.
By the reduction \eqref{eq:facuD}--\eqref{eq:reduD}, we may count the
number of these possibilities by Formula~\eqref{eq:2} with $n$ replaced
by $n-1$, $\rk T_j$ replaced by $\rk T_j-1$, and $m^{(j)}_3$ replaced
by $m^{(j)}_3-1$. This explains the second term in the factor 
in big parentheses on
the right-hand side of \eqref{eq:7}.
On the other hand, it may happen that one of the factors in \eqref{eq:facD}, 
let us say again, without loss of generality, $\si_d$,
contains two type $B$ cycles of length $1$ 
in its disjoint cycle decomposition;
that is, $\si_d$ may contain
$(n,\bar n)\,(a,\bar a).$
As a parabolic Coxeter element, this is of type $A_1^2$.
By the reduction \eqref{eq:facuD}--\eqref{eq:reduD}, we may count the
number of these possibilities by Formula~\eqref{eq:2} with $n$ replaced
by $n-1$, $\rk T_j$ replaced by $\rk T_j-1$, and $m^{(j)}_1$ replaced
by $m^{(j)}_1-2$. This explains the third term in the factor 
in big parentheses on
the right-hand side of \eqref{eq:7}.

\begin{figure}
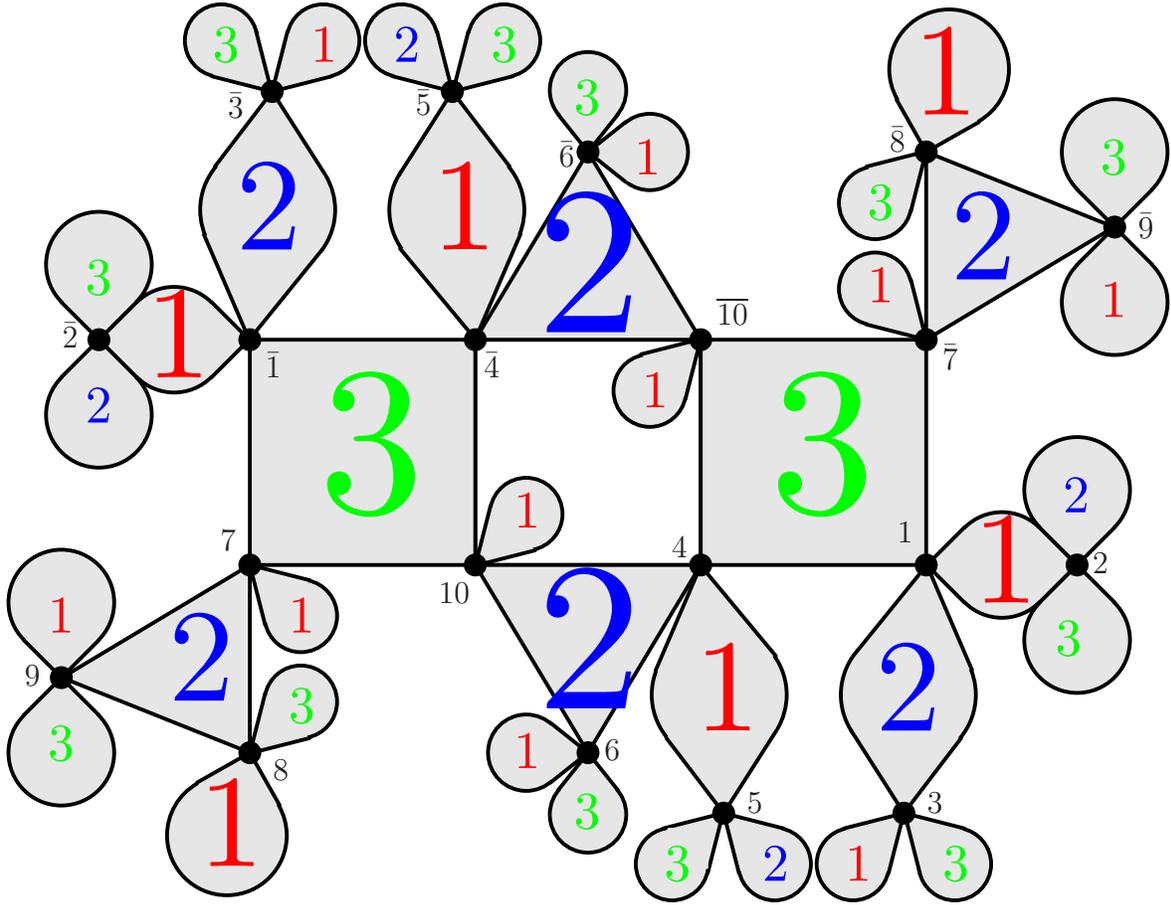

\centertexdraw{
  \drawdim truecm  \linewd.05
  \move(-1.5 -1.5)
  \rlvec(3 0)
  \rlvec(0 3)
  \rlvec(-3 0)
  \rlvec(0 -3)

  \move(1.5 -1.5)
  \rlvec(3 0)
  \rlvec(0 3)
  \rlvec(-3 0)
  \rlvec(0 -3)
  \lfill f:.9
  \htext(2.4 -.8){\Gruen{\Riesig3}}

  \move(-4.5 -1.5)
  \rlvec(3 0)
  \rlvec(0 3)
  \rlvec(-3 0)
  \rlvec(0 -3)
  \lfill f:.9
  \htext(-3.6 -.8){\Gruen{\Riesig3}}

  \move(-1.5 1.5)
  \rlvec(3 0)
  \rlvec(-1.5 2.5)
  \rlvec(-1.5 -2.5)
  \lfill f:.9
  \htext(-.7 1.6){\Blau{\Riesig2}}

  \move(-1.5 -1.5)
  \rlvec(3 0)
  \rlvec(-1.5 -2.5)
  \rlvec(-1.5 2.5)
  \lfill f:.9
  \htext(-.7 -3.4){\Blau{\Riesig2}}

  \move(4.5 1.5)
  \rlvec(2.5 1.5)
  \rlvec(-2.5 1)
  \rlvec(0 -2.5)
  \lfill f:.9
  \htext(4.8 2.3){\Blau{\riesig2}}

  \move(-4.5 -1.5)
  \rlvec(-2.5 -1.5)
  \rlvec(2.5 -1)
  \rlvec(0 2.5)
  \lfill f:.9
  \htext(-5.6 -3.3){\Blau{\riesig2}}

  \move(.82 .82)
  \fcir f:.9 r:.48
  \larc r:.48 sd:100 ed:-10
  \move(.7 1.3)
  \rlvec (.8 .2)
  \rlvec (-.2 -.8)
  \lfill f:.9
  \htext(.7 .6){\grosz\Rot{1}}

  \move(-.82 -.82)
  \fcir f:.9 r:.48
  \larc r:.48 sd:-80 ed:170
  \move(-.7 -1.3)
  \rlvec (-.8 -.2)
  \rlvec (.2 .8)
  \lfill f:.9
  \htext(-1 -1){\grosz\Rot{1}}

  \move(3.82 2.18)
  \fcir f:.9 r:.48
  \larc r:.48 sd:10 ed:-100
  \move(3.7 1.7)
  \rlvec (.8 -.2)
  \rlvec (-.2 .8)
  \lfill f:.9
  \htext(3.7 2){\grosz\Rot{1}}

  \move(-3.82 -2.18)
  \fcir f:.9 r:.48
  \larc r:.48 sd:190 ed:80
  \move(-3.7 -1.7)
  \rlvec (-.8 .2)
  \rlvec (.2 -.8)
  \lfill f:.9
  \htext(-4 -2.4){\grosz\Rot{1}}

  \move(4.8 5.1)
  \fcir f:.9 r:.8
  \larc r:.8 sd:-60 ed:210
  \move(4.1 4.7)
  \rlvec (0.4 -0.7)
  \rlvec (0.7 0.4)
  \lfill f:.9
  \htext(4.3 4.5){\riesig\Rot{1}}

  \move(3.82 3.32)
  \fcir f:.9 r:.48
  \larc r:.48 sd:100 ed:-10
  \move(3.7 3.8)
  \rlvec (.8 .2)
  \rlvec (-.2 -.8)
  \lfill f:.9
  \htext(3.7 3.1){\grosz\Gruen{3}}

  \move(-4.8 -5.1)
  \fcir f:.9 r:.8
  \larc r:.8 sd:120 ed:30
  \move(-4.1 -4.7)
  \rlvec (-0.4 0.7)
  \rlvec (-0.7 -0.4)
  \lfill f:.9
  \htext(-5.2 -5.5){\riesig\Rot{1}}

  \move(-3.82 -3.32)
  \fcir f:.9 r:.48
  \larc r:.48 sd:-80 ed:170
  \move(-3.7 -3.8)
  \rlvec (-.8 -.2)
  \rlvec (.2 .8)
  \lfill f:.9
  \htext(-4 -3.6){\grosz\Gruen{3}}

  \move(7 4)
  \fcir f:.9 r:.7
  \larc r:.7 sd:-45 ed:225
  \move(6.5 3.5)
  \rlvec (0.5 -0.5)
  \rlvec (0.5 0.5)
  \lfill f:.9
  \htext(6.8 3.7){\grosz\Gruen{3}}

  \move(7 2)
  \fcir f:.9 r:.7
  \larc r:.7 sd:135 ed:45
  \move(6.5 2.5)
  \rlvec (0.5 0.5)
  \rlvec (0.5 -0.5)
  \lfill f:.9
  \htext(6.8 1.8){\grosz\Rot{1}}

  \move(-7 -2)
  \fcir f:.9 r:.7
  \larc r:.7 sd:-45 ed:225
  \move(-7.5 -2.5)
  \rlvec (0.5 -0.5)
  \rlvec (0.5 0.5)
  \lfill f:.9
  \htext(-7.2 -2.4){\grosz\Rot{1}}

  \move(-7 -4)
  \fcir f:.9 r:.7
  \larc r:.7 sd:135 ed:45
  \move(-7.5 -3.5)
  \rlvec (0.5 0.5)
  \rlvec (0.5 -0.5)
  \lfill f:.9
  \htext(-7.2 -4.1){\grosz\Gruen{3}}

  \move(-.82 -4)
  \fcir f:.9 r:.51
  \larc r:.51 sd:45 ed:-45
  \move(-0.5 -3.6)
  \rlvec (0.5 -0.4)
  \rlvec (-0.5 -0.4)
  \lfill f:.9
  \htext(-1.0 -4.2){\grosz\Rot{1}}

  \move(0 -4.82)
  \fcir f:.9 r:.51
  \larc r:.51 sd:135 ed:45
  \move(-0.4 -4.5)
  \rlvec (0.4 0.5)
  \rlvec (0.4 -0.5)
  \lfill f:.9
  \htext(-.2 -5){\grosz\Gruen{3}}

  \move(.82 4)
  \fcir f:.9 r:.51
  \larc r:.51 sd:235 ed:135
  \move(0.5 3.6)
  \rlvec (-0.5 0.4)
  \rlvec (0.5 0.4)
  \lfill f:.9
  \htext(.6 3.7){\grosz\Rot{1}}

  \move(0 4.82)
  \fcir f:.9 r:.51
  \larc r:.51 sd:-45 ed:235
  \move(0.4 4.5)
  \rlvec (-0.4 -0.5)
  \rlvec (-0.4 0.5)
  \lfill f:.9
  \htext(-.2 4.5){\grosz\Gruen{3}}

  \move(1.74511 -3.22891)
  \fcir f:.9 r:.9
  \larc r:.9 sd:160 ed:218
  \larc r:.9 sd:-30 ed:35
  \move(.9 -2.9)
  \rlvec (0.6 1.4)
  \rlvec (1 -1.2)
  \lfill f:.9
  \move(1.04 -3.8)
  \rlvec (0.76 -1)
  \rlvec (.70 1.1)
  \lfill f:.9
  \htext(1.4 -3.7){\riesig\Rot{1}}

  \move(-1.74511 3.22891)
  \fcir f:.9 r:.9
  \larc r:.9 sd:145 ed:215
  \larc r:.9 sd:-30 ed:40
  \move(-.9 2.9)
  \rlvec (-0.6 -1.4)
  \rlvec (-1 1.2)
  \lfill f:.9
  \move(-1.04 3.8)
  \rlvec (-0.76 1)
  \rlvec (-.70 -1.1)
  \lfill f:.9
  \htext(-2.1 2.7){\riesig\Rot{1}}

  \move(4.255 -3.22891)
  \fcir f:.9 r:.9
  \larc r:.9 sd:145 ed:215
  \larc r:.9 sd:-40 ed:30
  \move(3.53 -2.7)
  \rlvec (.97 1.2)
  \rlvec (0.6 -1.4)
  \lfill f:.9
  \move(3.5 -3.7)
  \rlvec (.70 -1.1)
  \rlvec (0.76 1)
  \lfill f:.9
  \htext(3.8 -3.7){\riesig\Blau{2}}

  \move(-4.255 3.22891)
  \fcir f:.9 r:.9
  \larc r:.9 sd:140 ed:210
  \larc r:.9 sd:-37 ed:30
  \move(-3.53 2.7)
  \rlvec (-.97 -1.2)
  \rlvec (-0.6 1.4)
  \lfill f:.9
  \move(-3.5 3.7)
  \rlvec (-.70 1.1)
  \rlvec (-0.76 -1)
  \lfill f:.9
  \htext(-4.7 2.7){\riesig\Blau{2}}

  \move(1 -5)
  \rlvec (.8 .2)
  \rlvec (-.2 -.8)
  \lfill f:.9
  \move(1.12 -5.48)
  \fcir f:.9 r:.48
  \larc r:.48 sd:100 ed:-10
  \htext(1 -5.7){\grosz\Gruen{3}}

  \move(2.48 -5.48)
  \fcir f:.9 r:.48
  \larc r:.48 sd:190 ed:80
  \move(2.6 -5)
  \rlvec (-.8 .2)
  \rlvec (.2 -.8)
  \lfill f:.9
  \htext(2.3 -5.7){\grosz\Blau{2}}

  \move(3.52 -5.48)
  \fcir f:.9 r:.48
  \larc r:.48 sd:100 ed:-10
  \move(3.4 -5)
  \rlvec (.8 .2)
  \rlvec (-.2 -.8)
  \lfill f:.9
  \htext(3.4 -5.7){\grosz\Rot{1}}

  \move(4.88 -5.48)
  \fcir f:.9 r:.48
  \larc r:.48 sd:190 ed:80
  \move(5 -5)
  \rlvec (-.8 .2)
  \rlvec (.2 -.8)
  \lfill f:.9
  \htext(4.7 -5.7){\grosz\Gruen{3}}

  \move(-1.12 5.48)
  \fcir f:.9 r:.48
  \larc r:.48 sd:-80 ed:170
  \move(-1 5)
  \rlvec (-.8 -.2)
  \rlvec (.2 .8)
  \lfill f:.9
  \htext(-1.3 5.2){\grosz\Gruen{3}}

  \move(-2.48 5.48)
  \fcir f:.9 r:.48
  \larc r:.48 sd:10 ed:260
  \move(-2.6 5)
  \rlvec (.8 -.2)
  \rlvec (-.2 .8)
  \lfill f:.9
  \htext(-2.6 5.2){\grosz\Blau{2}}

  \move(-3.52 5.48)
  \fcir f:.9 r:.48
  \larc r:.48 sd:-80 ed:170
  \move(-3.4 5)
  \rlvec (-.8 -.2)
  \rlvec (.2 .8)
  \lfill f:.9
  \htext(-3.7 5.2){\grosz\Rot{1}}

  \move(-4.88 5.48)
  \fcir f:.9 r:.48
  \larc r:.48 sd:10 ed:260
  \move(-5 5)
  \rlvec (.8 -.2)
  \rlvec (-.2 .8)
  \lfill f:.9
  \htext(-5 5.2){\grosz\Gruen{3}}

  \move(5.5 -1.5)
  \fcir f:.9 r:.7
  \larc r:.7 sd:45 ed:135
  \larc r:.7 sd:225 ed:-45
  \move(5 -2)
  \rlvec (-.5 .5)
  \rlvec (.5 .5)
  \lfill f:.9
  \move(6 -2)
  \rlvec (.5 .5)
  \rlvec (-.5 .5)
  \lfill f:.9
  \htext(5.1 -2){\riesig\Rot{1}}

  \move(6.5 -.5)
  \fcir f:.9 r:.7
  \larc r:.7 sd:-45 ed:225
  \move(6 -1)
  \rlvec (0.5 -0.5)
  \rlvec (0.5 0.5)
  \lfill f:.9
  \htext(6.3 -.8){\grosz\Blau{2}}

  \move(6.5 -2.5)
  \fcir f:.9 r:.7
  \larc r:.7 sd:135 ed:45
  \move(6 -2)
  \rlvec (0.5 0.5)
  \rlvec (0.5 -0.5)
  \lfill f:.9
  \htext(6.2 -2.7){\grosz\Gruen{3}}

  \move(-5.5 1.5)
  \fcir f:.9 r:.7
  \larc r:.7 sd:45 ed:135
  \larc r:.7 sd:225 ed:-45
  \move(-5 2)
  \rlvec (.5 -.5)
  \rlvec (-.5 -.5)
  \lfill f:.9
  \move(-6 2)
  \rlvec (-.5 -.5)
  \rlvec (.5 -.5)
  \lfill f:.9
  \htext(-5.9 1){\riesig\Rot{1}}

  \move(-6.5 2.5)
  \fcir f:.9 r:.7
  \larc r:.7 sd:-45 ed:225
  \move(-7 2)
  \rlvec (0.5 -0.5)
  \rlvec (0.5 0.5)
  \lfill f:.9
  \htext(-6.7 2.1){\grosz\Gruen{3}}

  \move(-6.5 .5)
  \fcir f:.9 r:.7
  \larc r:.7 sd:135 ed:45
  \move(-7 1)
  \rlvec (0.5 0.5)
  \rlvec (0.5 -0.5)
  \lfill f:.9
  \htext(-6.7 .4){\grosz\Blau{2}}

  \Ringerl(-4.5 -1.5)
  \Ringerl(-1.5 -1.5)
  \Ringerl(1.5 -1.5)
  \Ringerl(4.5 -1.5)
  \Ringerl(-4.5 1.5)
  \Ringerl(-1.5 1.5)
  \Ringerl(1.5 1.5)
  \Ringerl(4.5 1.5)

  \Ringerl(4.5 4)
  \Ringerl(7 3)
  \Ringerl(-4.5 -4)
  \Ringerl(-7 -3)

  \Ringerl(0 4)
  \Ringerl(0 -4)

  \Ringerl(1.8 -4.8)
  \Ringerl(-1.8 4.8)
  \Ringerl(4.2 -4.8)
  \Ringerl(-4.2 4.8)

  \Ringerl(6.5 -1.5)
  \Ringerl(-6.5 1.5)

  \htext(4.1 -1.2){1}
  \htext(-4.3 1){$\bar 1$}

  \htext(6.7 -1.6){2}
  \htext(-7 1.4){$\bar 2$}

  \htext(4.5 -4.8){3}
  \htext(-4.8 4.45){$\bar 3$}

  \htext(1.1 -1.4){4}
  \htext(-1.4 1){$\bar 4$}

  \htext(2.1 -4.8){5}
  \htext(-2.3 4.5){$\bar 5$}

  \htext(.2 -4.1){6}
  \htext(-.4 3.8){$\bar 6$}

  \htext(-4.9 -1.3){7}
  \htext(4.7 1.1){$\bar 7$}

  \htext(-4.2 -4.35){$8$}
  \htext(4.0 4){$\bar 8$}

  \htext(-7.5 -3.1){9}
  \htext(7.3 2.85){$\bar 9$}

  \htext(-2 -2){10}
  \htext(1.7 1.7){$\overline{10}$}
}
\caption{The $3$-atoll corresponding to the factorisation
\eqref{eq:fac0D}}
\label{fig:5}
\end{figure}

From now on we may assume that none of the $\si_i$'s contains a type
$B$ cycle in its (type $B$) disjoint cycle decomposition.
To determine the number of minimal factorisations \eqref{eq:facD} in
this case,
we construct again a bijection between these factorisations and
certain maps. In what follows, we will still use the concept
of a rotator, introduced in the proof of Theorem~\ref{thm:2}.
We apply again the procedure described in that proof.
That is, for each (ordinary) 
cycle $(a_1,a_2,\dots,a_k)$ of $\si_i$, we create a $k$-gon
coloured $i$, label its vertices $a_1,a_2,\dots,a_k$ in clockwise
order, and glue these polygons
into a map by identifying equally labelled vertices such that
the rotator of each vertex is $(1,2,\dots,d)$.
However, this map can be embedded in the plane only if we allow
the creation of an inner face corresponding to the cycle $(n,\bar n)$
on the left-hand side of \eqref{eq:facD} (the outer face 
corresponding to the large cycle $(1,2,\dots,n-1,\bar1,\bar2,\dots,
\overline{n-1})$). Moreover, this inner face must be bounded by $2d$
edges. We call such a map, in which all faces except the
outer face and an inner face intersect only in vertices, 
and are coloured with colours from $\{1,2,\dots,d\}$, 
and in which the inner face is bounded by $2d$ edges, a {\it
$d$-atoll}.
For example, if we take $n=10$ and $d=3$,
and consider the factorisation
\begin{equation} \label{eq:fac0D}
(1,2,\dots,9,\bar1,\bar2,\dots,
\bar{9})\,(10,\overline{10})=\si_3\si_2\si_1,
\end{equation}
where $\si_3=((1,4,\overline{10},\bar 7))$,
$\si_2=((1,3))\,((4,6,10))\,((7,8,9))$,
and $\si_1=((1,2))\,((4,5))$, and apply this procedure, we obtain 
the $3$-atoll in Figure~\ref{fig:5}.
In the figure, the faces corresponding to cycles are shaded.
As in Figures~\ref{fig:3} and \ref{fig:4}, the outer face is not
shaded. Here, there is in addition an inner face which is not shaded,
the face formed by the vertices $4, 10, \bar 4,\overline{10}$.
Again, the colours of the shaded faces are indicated by the numbers $1$, $2$,
respectively $3$, placed in the centre of the faces.

Unsurprisingly, 
the fact that the result of the procedure can be realised as 
a $d$-atoll follows again from Euler's 
formula. More precisely, the number of faces corresponding to the polygons is 
$2{\sum _{i=1} ^{d}}\sum _{k=0} ^{n}m^{(i)}_k$, 
the number of edges
is $2\sum _{i=1} ^{d}\sum _{k=0} ^{n}m^{(i)}_k(k+1)$, and the number
of vertices is $2n$. Hence, if we include the outer face and the inner
face, the number
of vertices minus the number of edges plus the number of faces is
\begin{align} \notag
2n-
2\sum _{i=1} ^{d}\sum _{k=0} ^{n}m^{(i)}_k(k+1)+{}
&2\sum _{i=1} ^{d}\sum _{k=0} ^{n}m^{(i)}_k+2=
2n+2-
2\sum _{i=1} ^{d}\sum _{k=0} ^{n}k\cdot m^{(i)}_k\\
\notag
&=2n+2-2\rk T_1-2\rk T_2-\dots-2\rk T_d\\
&=2,
\label{eq:EulerD}
\end{align}
according to our assumption concerning the sum of the ranks of the
types $T_i$.

\begin{figure}
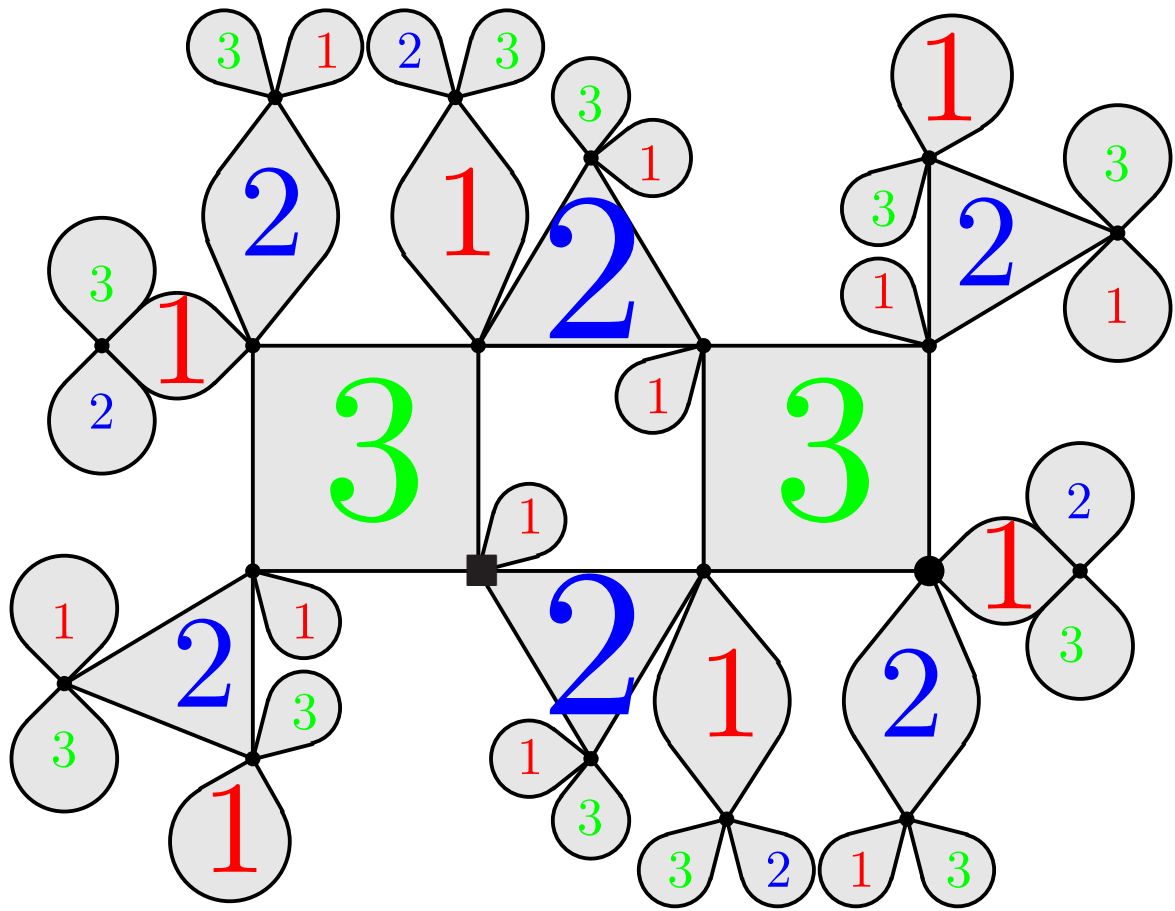

\centertexdraw{
  \drawdim truecm  \linewd.05
  \move(-1.5 -1.5)
  \rlvec(3 0)
  \rlvec(0 3)
  \rlvec(-3 0)
  \rlvec(0 -3)

  \move(1.5 -1.5)
  \rlvec(3 0)
  \rlvec(0 3)
  \rlvec(-3 0)
  \rlvec(0 -3)
  \lfill f:.9
  \htext(2.4 -.8){\Gruen{\Riesig3}}

  \move(-4.5 -1.5)
  \rlvec(3 0)
  \rlvec(0 3)
  \rlvec(-3 0)
  \rlvec(0 -3)
  \lfill f:.9
  \htext(-3.6 -.8){\Gruen{\Riesig3}}

  \move(-1.5 1.5)
  \rlvec(3 0)
  \rlvec(-1.5 2.5)
  \rlvec(-1.5 -2.5)
  \lfill f:.9
  \htext(-.7 1.6){\Blau{\Riesig2}}

  \move(-1.5 -1.5)
  \rlvec(3 0)
  \rlvec(-1.5 -2.5)
  \rlvec(-1.5 2.5)
  \lfill f:.9
  \htext(-.7 -3.4){\Blau{\Riesig2}}

  \move(4.5 1.5)
  \rlvec(2.5 1.5)
  \rlvec(-2.5 1)
  \rlvec(0 -2.5)
  \lfill f:.9
  \htext(4.8 2.3){\Blau{\riesig2}}

  \move(-4.5 -1.5)
  \rlvec(-2.5 -1.5)
  \rlvec(2.5 -1)
  \rlvec(0 2.5)
  \lfill f:.9
  \htext(-5.6 -3.3){\Blau{\riesig2}}

  \move(.82 .82)
  \fcir f:.9 r:.48
  \larc r:.48 sd:100 ed:-10
  \move(.7 1.3)
  \rlvec (.8 .2)
  \rlvec (-.2 -.8)
  \lfill f:.9
  \htext(.7 .6){\grosz\Rot{1}}

  \move(-.82 -.82)
  \fcir f:.9 r:.48
  \larc r:.48 sd:-80 ed:170
  \move(-.7 -1.3)
  \rlvec (-.8 -.2)
  \rlvec (.2 .8)
  \lfill f:.9
  \htext(-1 -1){\grosz\Rot{1}}

  \move(3.82 2.18)
  \fcir f:.9 r:.48
  \larc r:.48 sd:10 ed:-100
  \move(3.7 1.7)
  \rlvec (.8 -.2)
  \rlvec (-.2 .8)
  \lfill f:.9
  \htext(3.7 2){\grosz\Rot{1}}

  \move(-3.82 -2.18)
  \fcir f:.9 r:.48
  \larc r:.48 sd:190 ed:80
  \move(-3.7 -1.7)
  \rlvec (-.8 .2)
  \rlvec (.2 -.8)
  \lfill f:.9
  \htext(-4 -2.4){\grosz\Rot{1}}

  \move(4.8 5.1)
  \fcir f:.9 r:.8
  \larc r:.8 sd:-60 ed:210
  \move(4.1 4.7)
  \rlvec (0.4 -0.7)
  \rlvec (0.7 0.4)
  \lfill f:.9
  \htext(4.3 4.5){\riesig\Rot{1}}

  \move(3.82 3.32)
  \fcir f:.9 r:.48
  \larc r:.48 sd:100 ed:-10
  \move(3.7 3.8)
  \rlvec (.8 .2)
  \rlvec (-.2 -.8)
  \lfill f:.9
  \htext(3.7 3.1){\grosz\Gruen{3}}

  \move(-4.8 -5.1)
  \fcir f:.9 r:.8
  \larc r:.8 sd:120 ed:30
  \move(-4.1 -4.7)
  \rlvec (-0.4 0.7)
  \rlvec (-0.7 -0.4)
  \lfill f:.9
  \htext(-5.2 -5.5){\riesig\Rot{1}}

  \move(-3.82 -3.32)
  \fcir f:.9 r:.48
  \larc r:.48 sd:-80 ed:170
  \move(-3.7 -3.8)
  \rlvec (-.8 -.2)
  \rlvec (.2 .8)
  \lfill f:.9
  \htext(-4 -3.6){\grosz\Gruen{3}}

  \move(7 4)
  \fcir f:.9 r:.7
  \larc r:.7 sd:-45 ed:225
  \move(6.5 3.5)
  \rlvec (0.5 -0.5)
  \rlvec (0.5 0.5)
  \lfill f:.9
  \htext(6.8 3.7){\grosz\Gruen{3}}

  \move(7 2)
  \fcir f:.9 r:.7
  \larc r:.7 sd:135 ed:45
  \move(6.5 2.5)
  \rlvec (0.5 0.5)
  \rlvec (0.5 -0.5)
  \lfill f:.9
  \htext(6.8 1.8){\grosz\Rot{1}}

  \move(-7 -2)
  \fcir f:.9 r:.7
  \larc r:.7 sd:-45 ed:225
  \move(-7.5 -2.5)
  \rlvec (0.5 -0.5)
  \rlvec (0.5 0.5)
  \lfill f:.9
  \htext(-7.2 -2.4){\grosz\Rot{1}}

  \move(-7 -4)
  \fcir f:.9 r:.7
  \larc r:.7 sd:135 ed:45
  \move(-7.5 -3.5)
  \rlvec (0.5 0.5)
  \rlvec (0.5 -0.5)
  \lfill f:.9
  \htext(-7.2 -4.1){\grosz\Gruen{3}}

  \move(-.82 -4)
  \fcir f:.9 r:.51
  \larc r:.51 sd:45 ed:-45
  \move(-0.5 -3.6)
  \rlvec (0.5 -0.4)
  \rlvec (-0.5 -0.4)
  \lfill f:.9
  \htext(-1.0 -4.2){\grosz\Rot{1}}

  \move(0 -4.82)
  \fcir f:.9 r:.51
  \larc r:.51 sd:135 ed:45
  \move(-0.4 -4.5)
  \rlvec (0.4 0.5)
  \rlvec (0.4 -0.5)
  \lfill f:.9
  \htext(-.2 -5){\grosz\Gruen{3}}

  \move(.82 4)
  \fcir f:.9 r:.51
  \larc r:.51 sd:235 ed:135
  \move(0.5 3.6)
  \rlvec (-0.5 0.4)
  \rlvec (0.5 0.4)
  \lfill f:.9
  \htext(.6 3.7){\grosz\Rot{1}}

  \move(0 4.82)
  \fcir f:.9 r:.51
  \larc r:.51 sd:-45 ed:235
  \move(0.4 4.5)
  \rlvec (-0.4 -0.5)
  \rlvec (-0.4 0.5)
  \lfill f:.9
  \htext(-.2 4.5){\grosz\Gruen{3}}

  \move(1.74511 -3.22891)
  \fcir f:.9 r:.9
  \larc r:.9 sd:160 ed:218
  \larc r:.9 sd:-30 ed:35
  \move(.9 -2.9)
  \rlvec (0.6 1.4)
  \rlvec (1 -1.2)
  \lfill f:.9
  \move(1.04 -3.8)
  \rlvec (0.76 -1)
  \rlvec (.70 1.1)
  \lfill f:.9
  \htext(1.4 -3.7){\riesig\Rot{1}}

  \move(-1.74511 3.22891)
  \fcir f:.9 r:.9
  \larc r:.9 sd:145 ed:215
  \larc r:.9 sd:-30 ed:40
  \move(-.9 2.9)
  \rlvec (-0.6 -1.4)
  \rlvec (-1 1.2)
  \lfill f:.9
  \move(-1.04 3.8)
  \rlvec (-0.76 1)
  \rlvec (-.70 -1.1)
  \lfill f:.9
  \htext(-2.1 2.7){\riesig\Rot{1}}

  \move(4.255 -3.22891)
  \fcir f:.9 r:.9
  \larc r:.9 sd:145 ed:215
  \larc r:.9 sd:-40 ed:30
  \move(3.53 -2.7)
  \rlvec (.97 1.2)
  \rlvec (0.6 -1.4)
  \lfill f:.9
  \move(3.5 -3.7)
  \rlvec (.70 -1.1)
  \rlvec (0.76 1)
  \lfill f:.9
  \htext(3.8 -3.7){\riesig\Blau{2}}

  \move(-4.255 3.22891)
  \fcir f:.9 r:.9
  \larc r:.9 sd:140 ed:210
  \larc r:.9 sd:-37 ed:30
  \move(-3.53 2.7)
  \rlvec (-.97 -1.2)
  \rlvec (-0.6 1.4)
  \lfill f:.9
  \move(-3.5 3.7)
  \rlvec (-.70 1.1)
  \rlvec (-0.76 -1)
  \lfill f:.9
  \htext(-4.7 2.7){\riesig\Blau{2}}

  \move(1.12 -5.48)
  \fcir f:.9 r:.48
  \larc r:.48 sd:100 ed:-10
  \move(1 -5)
  \rlvec (.8 .2)
  \rlvec (-.2 -.8)
  \lfill f:.9
  \htext(1 -5.7){\grosz\Gruen{3}}

  \move(2.48 -5.48)
  \fcir f:.9 r:.48
  \larc r:.48 sd:190 ed:80
  \move(2.6 -5)
  \rlvec (-.8 .2)
  \rlvec (.2 -.8)
  \lfill f:.9
  \htext(2.3 -5.7){\grosz\Blau{2}}

  \move(3.52 -5.48)
  \fcir f:.9 r:.48
  \larc r:.48 sd:100 ed:-10
  \move(3.4 -5)
  \rlvec (.8 .2)
  \rlvec (-.2 -.8)
  \lfill f:.9
  \htext(3.4 -5.7){\grosz\Rot{1}}

  \move(4.88 -5.48)
  \fcir f:.9 r:.48
  \larc r:.48 sd:190 ed:80
  \move(5 -5)
  \rlvec (-.8 .2)
  \rlvec (.2 -.8)
  \lfill f:.9
  \htext(4.7 -5.7){\grosz\Gruen{3}}

  \move(-1.12 5.48)
  \fcir f:.9 r:.48
  \larc r:.48 sd:-80 ed:170
  \move(-1 5)
  \rlvec (-.8 -.2)
  \rlvec (.2 .8)
  \lfill f:.9
  \htext(-1.3 5.2){\grosz\Gruen{3}}

  \move(-2.48 5.48)
  \fcir f:.9 r:.48
  \larc r:.48 sd:10 ed:260
  \move(-2.6 5)
  \rlvec (.8 -.2)
  \rlvec (-.2 .8)
  \lfill f:.9
  \htext(-2.6 5.2){\grosz\Blau{2}}

  \move(-3.52 5.48)
  \fcir f:.9 r:.48
  \larc r:.48 sd:-80 ed:170
  \move(-3.4 5)
  \rlvec (-.8 -.2)
  \rlvec (.2 .8)
  \lfill f:.9
  \htext(-3.7 5.2){\grosz\Rot{1}}

  \move(-4.88 5.48)
  \fcir f:.9 r:.48
  \larc r:.48 sd:10 ed:260
  \move(-5 5)
  \rlvec (.8 -.2)
  \rlvec (-.2 .8)
  \lfill f:.9
  \htext(-5 5.2){\grosz\Gruen{3}}

  \move(5.5 -1.5)
  \fcir f:.9 r:.7
  \larc r:.7 sd:45 ed:135
  \move(5.5 -1.5)
  \larc r:.7 sd:225 ed:-45
  \move(5 -2)
  \rlvec (-.5 .5)
  \rlvec (.5 .5)
  \lfill f:.9
  \move(6 -2)
  \rlvec (.5 .5)
  \rlvec (-.5 .5)
  \lfill f:.9
  \htext(5.1 -2){\riesig\Rot{1}}

  \move(6.5 -.5)
  \fcir f:.9 r:.7
  \larc r:.7 sd:-45 ed:225
  \move(6 -1)
  \rlvec (0.5 -0.5)
  \rlvec (0.5 0.5)
  \lfill f:.9
  \htext(6.3 -.8){\grosz\Blau{2}}

  \move(6.5 -2.5)
  \fcir f:.9 r:.7
  \larc r:.7 sd:135 ed:45
  \move(6 -2)
  \rlvec (0.5 0.5)
  \rlvec (0.5 -0.5)
  \lfill f:.9
  \htext(6.2 -2.7){\grosz\Gruen{3}}

  \move(-5.5 1.5)
  \fcir f:.9 r:.7
  \larc r:.7 sd:45 ed:135
  \move(-5.5 1.5)
  \larc r:.7 sd:225 ed:-45
  \move(-5 2)
  \rlvec (.5 -.5)
  \rlvec (-.5 -.5)
  \lfill f:.9
  \move(-6 2)
  \rlvec (-.5 -.5)
  \rlvec (.5 -.5)
  \lfill f:.9
  \htext(-5.9 1){\riesig\Rot{1}}

  \move(-6.5 2.5)
  \fcir f:.9 r:.7
  \larc r:.7 sd:-45 ed:225
  \move(-7 2)
  \rlvec (0.5 -0.5)
  \rlvec (0.5 0.5)
  \lfill f:.9
  \htext(-6.7 2.1){\grosz\Gruen{3}}

  \move(-6.5 .5)
  \fcir f:.9 r:.7
  \larc r:.7 sd:135 ed:45
  \move(-7 1)
  \rlvec (0.5 0.5)
  \rlvec (0.5 -0.5)
  \lfill f:.9
  \htext(-6.7 .4){\grosz\Blau{2}}

  \ringerl(-4.5 -1.5)
  \htext(-1.7 -1.7){\Large $\blacksquare$}
  \ringerl(1.5 -1.5)
  \Mark(4.5 -1.5)
  \ringerl(-4.5 1.5)
  \ringerl(-1.5 1.5)
  \ringerl(1.5 1.5)
  \ringerl(4.5 1.5)

  \ringerl(4.5 4)
  \ringerl(7 3)
  \ringerl(-4.5 -4)
  \ringerl(-7 -3)

  \ringerl(0 4)
  \ringerl(0 -4)

  \ringerl(1.8 -4.8)
  \ringerl(-1.8 4.8)
  \ringerl(4.2 -4.8)
  \ringerl(-4.2 4.8)

  \ringerl(6.5 -1.5)
  \ringerl(-6.5 1.5)
}
\caption{A rotation-symmetric $3$-atoll with two marked vertices}
\label{fig:6}
\end{figure}

Again, we may further simplify this geometric representation of 
a minimal factorisation \eqref{eq:facD} by deleting all vertex
labels,
marking the vertex which had label $1$ with 
{\Large$\raise-1pt\hbox{$\bullet$}$}, and
marking the vertex that had label $n$ with $\blacksquare$. 
If this simplification is
applied to the $3$-atoll in Figure~\ref{fig:5}, we obtain the
$3$-atoll in Figure~\ref{fig:6}.
Clearly, if drawn appropriately into the plane,
a $d$-atoll resulting from an application of the above
procedure to a minimal factorisation \eqref{eq:facD} is 
symmetric with respect to a rotation by $180^\circ$, the centre of
the rotation being the centre of the inner face; cf.\ Figure~\ref{fig:6}.
As earlier, we shall abbreviate this property as {\it
rotation-symmetric}. In fact, there is not much freedom for the
choice of the vertex marked by $\blacksquare$ once a vertex has
been marked
by {\Large$\raise-1pt\hbox{$\bullet$}$}. Clearly, if we run through 
the vertex labelling process described in the proof of
Theorem~\ref{thm:2}, labelling $1$ the vertex which is marked by
{\Large$\raise-1pt\hbox{$\bullet$}$}, we shall reconstruct the labels
$1,2,\dots,n-1,\bar 1,\bar2,\dots,\overline{n-1}$. This leaves only
$2$ vertices incident to the inner face unlabelled, one of which
will have to carry the mark $\blacksquare$.

In summary, under the assumptions of Claim~(ii),
the number of minimal factorisations \eqref{eq:facD}, in which none of
the $\si_i$'s contains a type $B$ cycle in its disjoint cycle decomposition,
equals twice the number of all rotation-symmetric 
$d$-atolls on $2n$ vertices, in which one
vertex is marked by {\Large$\raise-1pt\hbox{$\bullet$}$},
all vertices have rotator $(1,2,\dots,d)^{\text
{O}}$, and with exactly $m^{(i)}_k$ pairs of faces of colour
$i$ having $k+1$ vertices, arranged symmetrically
around the inner face (which is not coloured).
Let us denote the number of these $d$-atolls by $N'_{D_n}(T_1,T_2,\dots,T_d)$.

We must now enumerate these $d$-atolls. 
First of all, introducing a figure of speech, 
we shall refer to coloured faces of a $d$-atoll
which share an edge with the inner face but {\it not\/} with the
outer face as faces ``inside the $d$-atoll," and all others as
faces ``outside the $d$-atoll."
For example, in Figure~\ref{fig:6} we find two faces inside the
$3$-atoll, namely the two loop faces attached to the vertices
labelled $10$, respectively $\overline{10}$, in Figure~\ref{fig:5}.
Since, in a $d$-atoll,
the inner face is bounded by exactly $2d$ edges,
inside the $d$-atoll, we find only coloured faces
containing exactly one vertex. Next, we travel counter-clockwise around the
inner face and record the coloured faces sharing an edge with both
the inner and outer faces. 
Thus we obtain a list of the form
$$F_1,F_2,\dots,F_\ell,F_{\ell+1},\dots,F_{2\ell},$$
where, except possibly for the marking, 
$F_{h+\ell}$ is an identical copy of $F_h$, $h=1,2,\dots,\ell$.
In Figure~\ref{fig:6}, this list contains four faces, 
$\tilde F_1,\tilde F_2,\tilde F_3,\tilde F_4$, where $\tilde F_1$ and 
$\tilde F_3$ are the two quadrangles of colour $3$, and where 
$\tilde F_2$ and $\tilde F_4$ are the two triangles of colour $2$
connecting the two quadrangles.

Continuing the general argument,
let the colour of $F_h$ be $i_h$. Inside the $d$-atoll, 
because of the rotator condition, there must be 
$\{i_{h+1}-i_h-1\}_d$ faces (containing just one vertex) 
incident to the common vertex of $F_h$ and
$F_{h+1}$ coloured $\{i_h+1\}_d,\dots,\{i_{h+1}-1\}_d$, where, by definition,
$$\{x\}_d:=\begin{cases} x,&\text {if $0\le x\le d$}\\
x+d,&\text {if $x<0$}\\
x-d,&\text {if $x>d$},\end{cases}$$
and where $i_{h+\ell}=i_h$, $h=1,2,\dots,\ell$. Here, if 
$\{i_h+1\}_d>\{i_{h+1}-1\}_d$, the sequence of colours
$\{i_h+1\}_d,\dots,\{i_{h+1}-1\}_d$ must be interpreted ``cyclically,"
that is, as
$\{i_h+1\}_d,\{i_h+1\}_d+1,\dots,d,1,2,\dots,\{i_{h+1}-1\}_d$.
As we observed above, the number of edges bounding the inner face is
$2d$. On the other hand, using the notation just introduced, this
number also equals
$$2\sum _{h=1} ^{\ell}\{i_{h+1}-i_h\}_d=
2\sum _{h=1} ^{\ell}\Big((i_{h+1}-i_h)+d\cdot\chi(i_{h+1}<i_h)\Big)=
2d\sum _{h=1} ^{\ell}\chi(i_{h+1}<i_h).
$$
Hence, there is precisely one $h$ for which $i_{h+1}<i_h$. Without
loss of generality, we may assume that $h=\ell$, so that 
$i_1<i_2<\dots<i_\ell$.

The ascending colouring of the faces $F_1,F_2,\dots,F_\ell$ breaks
the (rotation) symmetry of the $d$-atoll. Therefore, we may 
first enumerate $d$-atolls without any marking, and multiply the result
by the number of all possible markings, which is $n-1$.
More precisely, let $N''_{D_n}(T_1,T_2,\dots,T_d)$ denote the
number of all rotation-symmetric 
$d$-atolls on $2n$ vertices, in which
all vertices have rotator $(1,2,\dots,d)^{\text
{O}}$, and with exactly $m^{(i)}_k$ pairs of faces of colour
$i$ having $k+1$ vertices, arranged symmetrically
around the inner face (which is not coloured). Then,
\begin{align} \notag
N_{D_n}(T_1,T_2,\dots,T_d)
&=2N'_{D_n}(T_1,T_2,\dots,T_d)\\
&=2(n-1)N''_{D_n}(T_1,T_2,\dots,T_d).
\label{eq:N''}
\end{align}

We use a generating function approach 
to determine $N''_{D_n}(T_1,T_2,\dots,T_d)$,
which requires a combinatorial decomposition of our objects.
Let $G(\mathbf z)$ be the generating function
\begin{equation} \label{eq:Gj}
G(\mathbf z)=\sum _{A\in \mathcal A} ^{}w(A),
\end{equation}
where $\mathcal A$ is the set of all rotation-symmetric $d$-atolls,
in which all vertices have rotator $(1,2,\dots,d)^{\text
{O}}$, and where
$$w(A)=\prod _{i=1} ^{d}z_i^{\frac {1} {2}\#(\text {faces of $A$ with colour $i$})}
\prod _{i=1} ^{d}\prod _{k=1}
^{\infty}p_{i,k}^{\frac {1} {2}\#(\text {faces of $A$ 
with colour $i$ and $k$ vertices})}.$$
Here, $\mathbf z=(z_1,z_2,\dots,z_d)$, with the $z_i$'s,
$i=1,2,\dots,d$, and
the $p_{i,k}$, $i=1,2,\dots,d$, $k=1,2,\dots$, being indeterminates.
Clearly, in view of the bijection
between minimal factorisations \eqref{eq:facD} and $d$-atolls 
described earlier, and by \eqref{eq:N''}, we have
\begin{equation} \label{eq:coefD1}
N_{D_n}(T_1,T_2,\dots,T_d)=2(n-1)\coef{\mathbf z^{\mathbf c}
\prod _{i=1} ^{d}\prod _{k=0} ^{n}p_{i,k+1}^{m^{(i)}_k}}
G(\mathbf z),
\end{equation}
where $\mathbf c=(c_1,c_2,\dots,c_d)$, with $c_i$ equal to the number
of type $A$ cycles of $\si_i$; that is, $c_i=\sum _{k=0} ^{n}m^{(i)}_k$,
$i=1,2,\dots,d$. Here, and in the sequel, we use the multi-index
notation introduced at the beginning of Section~\ref{sec:1a}.
For later use, we observe that, for all $i$,
$c_i$ is related to $\rk T_i$ via
\begin{equation} \label{eq:crkD}
c_i=n-\rk T_i.
\end{equation}

\begin{figure}
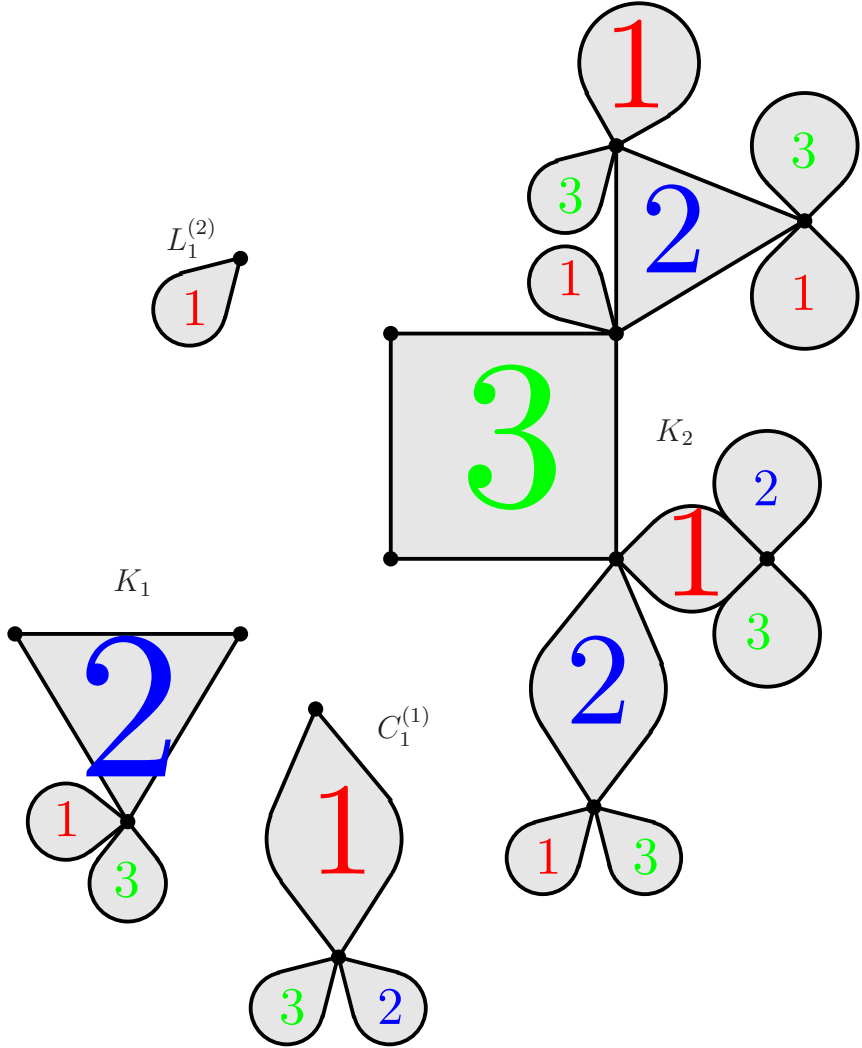

\centertexdraw{
  \drawdim truecm  \linewd.05

\bsegment
  \htext(-.2 -1){$K_1$}

  \move(-1.5 -1.5)
  \rlvec(3 0)
  \rlvec(-1.5 -2.5)
  \rlvec(-1.5 2.5)
  \lfill f:.9
  \htext(-.7 -3.4){\Blau{\Riesig2}}

  \move(-.82 -4)
  \fcir f:.9 r:.51
  \larc r:.51 sd:45 ed:-45
  \move(-0.5 -3.6)
  \rlvec (0.5 -0.4)
  \rlvec (-0.5 -0.4)
  \lfill f:.9
  \htext(-1.0 -4.2){\grosz\Rot{1}}

  \move(0 -4.82)
  \fcir f:.9 r:.51
  \larc r:.51 sd:135 ed:45
  \move(-0.4 -4.5)
  \rlvec (0.4 0.5)
  \rlvec (0.4 -0.5)
  \lfill f:.9
  \htext(-.2 -5){\grosz\Gruen{3}}

  \ringerl(-1.5 -1.5)
  \ringerl(1.5 -1.5)
  \ringerl(0 -4)
\esegment

\move(1 -1)
\bsegment
  \htext(2.3 -2){$C_1^{(1)}$}

  \move(1.74511 -3.22891)
  \fcir f:.9 r:.9
  \larc r:.9 sd:160 ed:218
  \larc r:.9 sd:-30 ed:35
  \move(.9 -2.9)
  \rlvec (0.6 1.4)
  \rlvec (1 -1.2)
  \lfill f:.9
  \move(1.04 -3.8)
  \rlvec (0.76 -1)
  \rlvec (.70 1.1)
  \lfill f:.9
  \htext(1.4 -3.7){\riesig\Rot{1}}

  \move(1.12 -5.48)
  \fcir f:.9 r:.48
  \larc r:.48 sd:100 ed:-10
  \move(1 -5)
  \rlvec (.8 .2)
  \rlvec (-.2 -.8)
  \lfill f:.9
  \htext(1 -5.7){\grosz\Gruen{3}}

  \move(2.48 -5.48)
  \fcir f:.9 r:.48
  \larc r:.48 sd:190 ed:80
  \move(2.6 -5)
  \rlvec (-.8 .2)
  \rlvec (.2 -.8)
  \lfill f:.9
  \htext(2.3 -5.7){\grosz\Blau{2}}

  \ringerl(1.5 -1.5)
  \ringerl(1.8 -4.8)
\esegment

\move(2 1)
\bsegment
  \htext(5 0){$K_2$}

  \move(1.5 -1.5)
  \rlvec(3 0)
  \rlvec(0 3)
  \rlvec(-3 0)
  \rlvec(0 -3)
  \lfill f:.9
  \htext(2.4 -.8){\Gruen{\Riesig3}}

  \move(4.5 1.5)
  \rlvec(2.5 1.5)
  \rlvec(-2.5 1)
  \rlvec(0 -2.5)
  \lfill f:.9
  \htext(4.8 2.3){\Blau{\riesig2}}

  \move(3.82 2.18)
  \fcir f:.9 r:.48
  \larc r:.48 sd:10 ed:-100
  \move(3.7 1.7)
  \rlvec (.8 -.2)
  \rlvec (-.2 .8)
  \lfill f:.9
  \htext(3.7 2){\grosz\Rot{1}}

  \move(4.8 5.1)
  \fcir f:.9 r:.8
  \larc r:.8 sd:-60 ed:210
  \move(4.1 4.7)
  \rlvec (0.4 -0.7)
  \rlvec (0.7 0.4)
  \lfill f:.9
  \htext(4.3 4.5){\riesig\Rot{1}}

  \move(3.82 3.32)
  \fcir f:.9 r:.48
  \larc r:.48 sd:100 ed:-10
  \move(3.7 3.8)
  \rlvec (.8 .2)
  \rlvec (-.2 -.8)
  \lfill f:.9
  \htext(3.7 3.1){\grosz\Gruen{3}}

  \move(7 4)
  \fcir f:.9 r:.7
  \larc r:.7 sd:-45 ed:225
  \move(6.5 3.5)
  \rlvec (0.5 -0.5)
  \rlvec (0.5 0.5)
  \lfill f:.9
  \htext(6.8 3.7){\grosz\Gruen{3}}

  \move(7 2)
  \fcir f:.9 r:.7
  \larc r:.7 sd:135 ed:45
  \move(6.5 2.5)
  \rlvec (0.5 0.5)
  \rlvec (0.5 -0.5)
  \lfill f:.9
  \htext(6.8 1.8){\grosz\Rot{1}}

  \move(4.255 -3.22891)
  \fcir f:.9 r:.9
  \larc r:.9 sd:145 ed:215
  \larc r:.9 sd:-40 ed:30
  \move(3.53 -2.7)
  \rlvec (.97 1.2)
  \rlvec (0.6 -1.4)
  \lfill f:.9
  \move(3.5 -3.7)
  \rlvec (.70 -1.1)
  \rlvec (0.76 1)
  \lfill f:.9
  \htext(3.8 -3.7){\riesig\Blau{2}}

  \move(3.52 -5.48)
  \fcir f:.9 r:.48
  \larc r:.48 sd:100 ed:-10
  \move(3.4 -5)
  \rlvec (.8 .2)
  \rlvec (-.2 -.8)
  \lfill f:.9
  \htext(3.4 -5.7){\grosz\Rot{1}}

  \move(4.88 -5.48)
  \fcir f:.9 r:.48
  \larc r:.48 sd:190 ed:80
  \move(5 -5)
  \rlvec (-.8 .2)
  \rlvec (.2 -.8)
  \lfill f:.9
  \htext(4.7 -5.7){\grosz\Gruen{3}}

  \move(5.5 -1.5)
  \fcir f:.9 r:.7
  \larc r:.7 sd:45 ed:135
  \move(5.5 -1.5)
  \larc r:.7 sd:225 ed:-45
  \move(5 -2)
  \rlvec (-.5 .5)
  \rlvec (.5 .5)
  \lfill f:.9
  \move(6 -2)
  \rlvec (.5 .5)
  \rlvec (-.5 .5)
  \lfill f:.9
  \htext(5.1 -2){\riesig\Rot{1}}

  \move(6.5 -.5)
  \fcir f:.9 r:.7
  \larc r:.7 sd:-45 ed:225
  \move(6 -1)
  \rlvec (0.5 -0.5)
  \rlvec (0.5 0.5)
  \lfill f:.9
  \htext(6.3 -.8){\grosz\Blau{2}}

  \move(6.5 -2.5)
  \fcir f:.9 r:.7
  \larc r:.7 sd:135 ed:45
  \move(6 -2)
  \rlvec (0.5 0.5)
  \rlvec (0.5 -0.5)
  \lfill f:.9
  \htext(6.2 -2.7){\grosz\Gruen{3}}

  \ringerl(4.5 -1.5)
  \ringerl(4.5 1.5)
  \ringerl(1.5 1.5)
  \ringerl(1.5 -1.5)
  \ringerl(4.5 4)
  \ringerl(7 3)
  \ringerl(4.2 -4.8)
  \ringerl(6.5 -1.5)
\esegment

\move(0 2)
\bsegment
  \htext(.5 1.5){$L_1^{(2)}$}

  \move(.82 .82)
  \fcir f:.9 r:.48
  \larc r:.48 sd:100 ed:-10
  \move(.7 1.3)
  \rlvec (.8 .2)
  \rlvec (-.2 -.8)
  \lfill f:.9
  \htext(.7 .6){\grosz\Rot{1}}

  \ringerl(1.5 1.5)
\esegment

}
\caption{The decomposition of the $3$-atoll in Figure~\ref{fig:6}}
\label{fig:6a}
\end{figure}

Now, let $A$ be a $d$-atoll in $\mathcal A$ such that the faces which
share an edge with both the inner and outer faces are 
$$F_1,F_2,\dots,F_\ell,F_{\ell+1},\dots,F_{2\ell},$$
where $F_{h+\ell}$ is an identical copy of $F_h$,
where the colour of $F_h$ is $i_h$, $h=1,2,\dots,\ell$, and with
$i_1<i_2<\dots<i_\ell$. We decompose $A$ by 
separating from each other the polygons which touch in vertices
of the inner face. The decomposition in the case of our example in
Figure~\ref{fig:6} is shown in Figure~\ref{fig:6a}. 
Ignoring identical copies which are there due
to the rotation symmetry, we obtain a list
\begin{multline} \label{eq:decomp}
K_1,L^{(1)}_{i_1+1},\dots,L^{(1)}_{i_2-1},
C^{(1)}_{i_2+1},\dots,C^{(1)}_d,C^{(1)}_1,\dots,C^{(1)}_{i_1-1},\\
K_2,L^{(2)}_{i_2+1},\dots,L^{(2)}_{i_3-1},
C^{(2)}_{i_3+1},\dots,C^{(2)}_d,C^{(2)}_1,\dots,C^{(2)}_{i_2-1},\dots\\
K_\ell,L^{(\ell)}_{i_\ell+1},\dots,L^{(\ell)}_d,
L^{(\ell)}_1,\dots,L^{(\ell)}_{i_1-1},
C^{(\ell)}_{i_1+1},\dots,C^{(\ell)}_{i_\ell-1},
\end{multline}
where $K_h$ is the $d$-cactus containing the face $F_h$, and, hence,
a $d$-cactus in which all but two neighbouring vertices
have rotator $(1,2,\dots,d)^{\text {O}}$, the latter two
vertices being incident to just one face, which is of colour $i_h$, 
where $L^{(h)}_j$ is a face of colour $j$ with just one vertex,
and where $C^{(h)}_j$ is a $d$-cactus in which all but
one vertex have rotator $(1,2,\dots,d)^{\text {O}}$, the 
distinguished
vertex being incident to just one face, which is of colour $j$,
$h=1,2,\dots,\ell$ and $j=1,2,\dots,d$. 
With this notation,
our example in Figure~\ref{fig:6a} is one in which $\ell=2$, $i_1=2$,
$i_2=3$.  

\begin{figure}
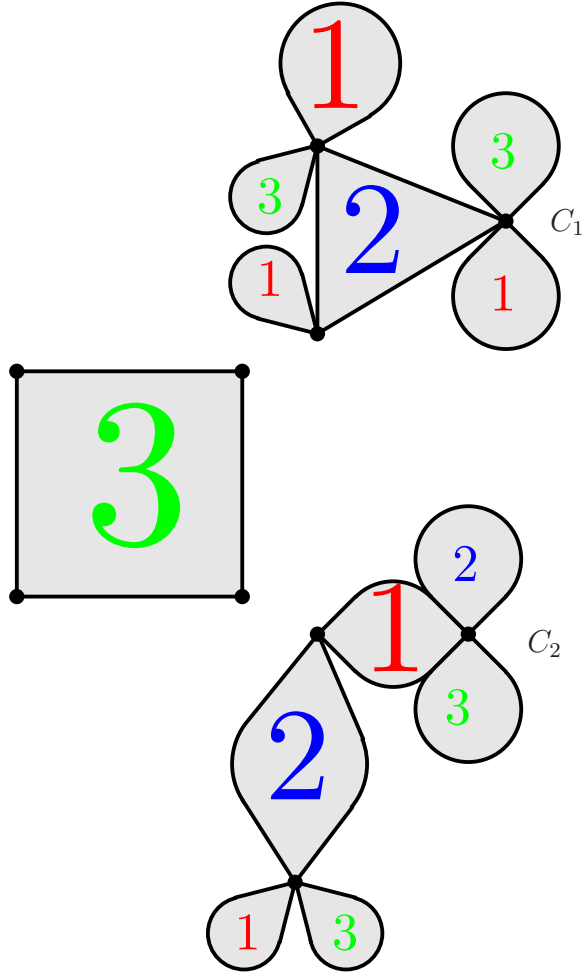

\centertexdraw{
  \drawdim truecm  \linewd.05

\bsegment
  \move(1.5 -1.5)
  \rlvec(3 0)
  \rlvec(0 3)
  \rlvec(-3 0)
  \rlvec(0 -3)
  \lfill f:.9
  \htext(2.4 -.8){\Gruen{\Riesig3}}

  \ringerl(1.5 1.5)
  \ringerl(4.5 1.5)
  \ringerl(1.5 -1.5)
  \ringerl(4.5 -1.5)
\esegment

\move(1 .5)
\bsegment
  \htext(7.6 2.8){$C_1$}

  \move(4.5 1.5)
  \rlvec(2.5 1.5)
  \rlvec(-2.5 1)
  \rlvec(0 -2.5)
  \lfill f:.9
  \htext(4.8 2.3){\Blau{\riesig2}}

  \move(3.82 2.18)
  \fcir f:.9 r:.48
  \larc r:.48 sd:10 ed:-100
  \move(3.7 1.7)
  \rlvec (.8 -.2)
  \rlvec (-.2 .8)
  \lfill f:.9
  \htext(3.7 2){\grosz\Rot{1}}

  \move(4.8 5.1)
  \fcir f:.9 r:.8
  \larc r:.8 sd:-60 ed:210
  \move(4.1 4.7)
  \rlvec (0.4 -0.7)
  \rlvec (0.7 0.4)
  \lfill f:.9
  \htext(4.3 4.5){\riesig\Rot{1}}

  \move(3.82 3.32)
  \fcir f:.9 r:.48
  \larc r:.48 sd:100 ed:-10
  \move(3.7 3.8)
  \rlvec (.8 .2)
  \rlvec (-.2 -.8)
  \lfill f:.9
  \htext(3.7 3.1){\grosz\Gruen{3}}

  \move(7 4)
  \fcir f:.9 r:.7
  \larc r:.7 sd:-45 ed:225
  \move(6.5 3.5)
  \rlvec (0.5 -0.5)
  \rlvec (0.5 0.5)
  \lfill f:.9
  \htext(6.8 3.7){\grosz\Gruen{3}}

  \move(7 2)
  \fcir f:.9 r:.7
  \larc r:.7 sd:135 ed:45
  \move(6.5 2.5)
  \rlvec (0.5 0.5)
  \rlvec (0.5 -0.5)
  \lfill f:.9
  \htext(6.8 1.8){\grosz\Rot{1}}

  \ringerl(4.5 1.5)
  \ringerl(4.5 4)
  \ringerl(7 3)
\esegment

\move(1 -.5)
\bsegment
  \htext(7.3 -1.8){$C_2$}

  \move(4.255 -3.22891)
  \fcir f:.9 r:.9
  \larc r:.9 sd:145 ed:215
  \larc r:.9 sd:-40 ed:30
  \move(3.53 -2.7)
  \rlvec (.97 1.2)
  \rlvec (0.6 -1.4)
  \lfill f:.9
  \move(3.5 -3.7)
  \rlvec (.70 -1.1)
  \rlvec (0.76 1)
  \lfill f:.9
  \htext(3.8 -3.7){\riesig\Blau{2}}

  \move(3.52 -5.48)
  \fcir f:.9 r:.48
  \larc r:.48 sd:100 ed:-10
  \move(3.4 -5)
  \rlvec (.8 .2)
  \rlvec (-.2 -.8)
  \lfill f:.9
  \htext(3.4 -5.7){\grosz\Rot{1}}

  \move(4.88 -5.48)
  \fcir f:.9 r:.48
  \larc r:.48 sd:190 ed:80
  \move(5 -5)
  \rlvec (-.8 .2)
  \rlvec (.2 -.8)
  \lfill f:.9
  \htext(4.7 -5.7){\grosz\Gruen{3}}

  \move(5.5 -1.5)
  \fcir f:.9 r:.7
  \larc r:.7 sd:45 ed:135
  \move(5.5 -1.5)
  \larc r:.7 sd:225 ed:-45
  \move(5 -2)
  \rlvec (-.5 .5)
  \rlvec (.5 .5)
  \lfill f:.9
  \move(6 -2)
  \rlvec (.5 .5)
  \rlvec (-.5 .5)
  \lfill f:.9
  \htext(5.1 -2){\riesig\Rot{1}}

  \move(6.5 -.5)
  \fcir f:.9 r:.7
  \larc r:.7 sd:-45 ed:225
  \move(6 -1)
  \rlvec (0.5 -0.5)
  \rlvec (0.5 0.5)
  \lfill f:.9
  \htext(6.3 -.8){\grosz\Blau{2}}

  \move(6.5 -2.5)
  \fcir f:.9 r:.7
  \larc r:.7 sd:135 ed:45
  \move(6 -2)
  \rlvec (0.5 0.5)
  \rlvec (0.5 -0.5)
  \lfill f:.9
  \htext(6.2 -2.7){\grosz\Gruen{3}}

  \ringerl(4.5 -1.5)
  \ringerl(4.2 -4.8)
  \ringerl(6.5 -1.5)
\esegment
}
\caption{The decomposition of $K_2$ in Figure~\ref{fig:6a}}
\label{fig:6b}
\end{figure}

The $d$-cacti $K_h$ can be further decomposed. Namely, assuming that
the face $F_h$ is a $k$-gon (of colour $i_h$), let
$C_1,C_2,\dots,C_{k-2}$
be the $d$-cacti incident to this $k$-gon, read in clockwise order,
starting with the $d$-cactus to the left of the two distinguished
vertices. Figure~\ref{fig:6b} illustrates this further decomposition
of the $d$-cactus $K_2$ from Figure~\ref{fig:6a}.
After removal of $F_h$, we are left with the
ordered collection
$C_1,C_2,\dots,C_{k-2}$ of $d$-cacti, each of which having 
the property that the
rotator of all but one vertex is $(1,2,\dots,d)^{\text {O}}$, the
exceptional vertex having rotator $(1,\dots,i_h-1,i_h+1,\dots,d)^{\text
{O}}$. By separating from each other 
the polygons of colours $1,\dots,i_h-1,i_h+1,\dots,d$
which touch in the exceptional vertex,
each $d$-cactus $C_i$ in turn can be decomposed into
$d$-cacti $C_{i,1},\dots,C_{i,{i_h}-1},C_{i,{i_h}+1},\dots,C_{i,d}$ with
$C_{i,\II}\in\mathcal C_\II$ for all $k$,
where $\mathcal C_\II$ denotes the set of all $d$-cacti in which all but
one vertex have rotator
$(1,2,\dots,d)^{\text {O}}$, the distinguished
vertex being incident to just one face, which is of colour $\II$.

Let $\om_\II(\mathbf z)$ denote the generating function for the
$d$-cacti in $\mathcal C_\II$, that is,
\begin{equation} \label{eq:omdef}
\om_\II(\mathbf z)=\sum _{C\in\mathcal C_\II} ^{}w(C).
\end{equation}
Furthermore,
for $i=1,2,\dots,d$, define the formal power series $P_i(u)$ in one variable
$u$ via
$$P_i(u)=\sum _{k=1} ^{\infty}p_{i,k}u^{k-1}.$$
Then, by the decomposition \eqref{eq:decomp} and the further
decomposition of the $K_h$'s that we just described, 
the contribution of the above 
$d$-atolls to the generating function \eqref{eq:Gj} is
\begin{align*}
\Bigg(\prod _{j=1} ^{\ell}&\frac {z_{i_j}\om_{i_j}(\mathbf z)} 
{\om_{1}(\mathbf z)\cdots \om_{d}(\mathbf z)}
\(P_{i_j}\(\frac 
{\om_{1}(\mathbf z)\cdots \om_{d}(\mathbf z)}
{\om_{i_j}(\mathbf z)}\)-p_{i_j,1}\)\Bigg)\\
&\kern4cm
\times \(\frac{\prod _{j=1} ^{d}z_jp_{j,1}} 
{\prod _{j=1} ^{\ell}z_{i_j}p_{i_j,1}}\)
\frac {(\om_{1}(\mathbf z)\cdots \om_{d}(\mathbf z))^{\ell-1}} 
{\prod _{j=1} ^{\ell}\om_{i_j}(\mathbf z)}\\
&=\Bigg(\prod _{j=1} ^{d}\frac {z_jp_{j,1}} 
{\om_{j}(\mathbf z)}\Bigg)
\prod _{j=1} ^{\ell}\(\frac {P_{i_j}\(\frac  
{\om_{1}(\mathbf z)\cdots \om_{d}(\mathbf z)}
{\om_{i_j}(\mathbf z)}\)} {p_{i_j,1}}-1\),
\end{align*}
the term in the first line corresponding to the contribution of
the $K_j$'s, the first term in the second line 
corresponding to the contribution of
the $L^{(j)}_k$'s, and the second term in the second line 
corresponding to the contribution of
the $C^{(j)}_k$'s. These expressions must be summed over 
$\ell=2,3,\dots,d$ and all possible choices of $1\le
i_1<i_2<\dots<i_\ell\le d$ to obtain the desired generating function
$G(\mathbf z)$, that is,
\begin{align} \notag
G(\mathbf z)&=
\Bigg(\prod _{j=1} ^{d}\frac {z_jp_{j,1}} 
{\om_{j}(\mathbf z)}\Bigg)
\sum _{\ell=2} ^{d}\sum _{1\le
i_1<i_2<\dots<i_\ell\le d} ^{}
\prod _{j=1} ^{\ell}\(\frac {P_{i_j}\(\frac  
{\om_{1}(\mathbf z)\cdots \om_{d}(\mathbf z)}
{\om_{i_j}(\mathbf z)}\)} {p_{i_j},1}-1\)\\
&=\Bigg(\prod _{j=1} ^{d}\frac {z_jp_{j,1}} 
{\om_{j}(\mathbf z)}\Bigg)
\(\prod _{j=1} ^{d}\frac {P_{j}\(\frac  
{\om_{1}(\mathbf z)\cdots \om_{d}(\mathbf z)}
{\om_{j}(\mathbf z)}\)} {p_{j,1}}-
\sum _{j=1} ^{d}\(\frac {P_{j}\(\frac  
{\om_{1}(\mathbf z)\cdots \om_{d}(\mathbf z)}
{\om_{j}(\mathbf z)}\)} {p_{j,1}}-1\)-1\).
\label{eq:G(z)}
\end{align}
Here we have used the elementary identity
$$\sum _{\ell=0} ^{d}\sum _{1\le
i_1<i_2<\dots<i_\ell\le d} ^{}X_{i_1}X_{i_2}\cdots X_{i_\ell}=
(1+X_1)(1+X_2)\cdots (1+X_d).
$$

Before we are able to proceed, we must
find functional equations for the generating functions $\om_j(\mathbf
z)$, $j=1,2,\dots,d$.
Given a $d$-cactus $C$ in $\mathcal C_\II$ such that the
distinguished vertex
is incident to a $k$-gon (of colour $\II$), we decompose it in a
manner analogous to the decomposition of $K_h$ above. To be more
precise, let $C_1,C_2,\dots,C_{k-1}$
be the $d$-cacti incident to this $k$-gon, read in clockwise order,
starting with the $d$-cactus to the left of the distinguished vertex.
After removal of the $k$-gon, we are left with the
ordered collection
$C_1,C_2,\dots,C_{k-1}$ of $d$-cacti, each of which having 
the property that the
rotator of all but one vertex is $(1,2,\dots,d)^{\text {O}}$, the
exceptional vertex having rotator $(1,\dots,\II-1,\II+1,\dots,d)^{\text
{O}}$. By separating from each other 
the polygons of colours $1,\dots,\II-1,\II+1,\dots,d$
which touch in the exceptional vertex,
each $d$-cactus $C_\JJ$ in turn can be decomposed into
$d$-cacti $C_{\JJ,1},\dots,C_{\JJ,\II-1},C_{\JJ,\II+1},\dots,C_{\JJ,d}$ with
$C_{\JJ,k}\in\mathcal C_k$ for all $k$. The upshot of these
combinatorial considerations is that
$$
\om_\II(\mathbf z)=z_\II P_\II(\om_1(\mathbf z)\cdots\om_d(\mathbf z)/
\om_\II(\mathbf z)),\quad \quad 
\II=1,2,\dots,d,
$$
or, equivalently, 
\begin{equation*} 
z_\II=\frac {\om_\II(\mathbf z)} {P_\II(\om_1(\mathbf z)\cdots\om_d(\mathbf z)/
\om_\II(\mathbf z))},\quad \quad 
\II=1,2,\dots,d.
\end{equation*}

Using this relation, the expression \eqref{eq:G(z)} for $G(\mathbf z)$ may 
now be further simplified, and we obtain
\begin{equation*}
G(\mathbf z)=
1-
\prod _{j=1} ^{d}\frac {p_{j,1}} 
{P_{j}\(\frac  
{\om_{1}(\mathbf z)\cdots \om_{d}(\mathbf z)}
{\om_{j}(\mathbf z)}\)}
\sum _{j=1} ^{d}\frac {P_{j}\(\frac  
{\om_{1}(\mathbf z)\cdots \om_{d}(\mathbf z)}
{\om_{j}(\mathbf z)}\)} {p_{j,1}}+
(d-1)\prod _{j=1} ^{d}\frac {p_{j,1}} 
{P_{j}\(\frac  
{\om_{1}(\mathbf z)\cdots \om_{d}(\mathbf z)}
{\om_{j}(\mathbf z)}\)}.
\end{equation*}
This is substituted in \eqref{eq:coefD1}, to obtain
\begin{multline} \label{eq:coefD}
N_{D_n}(T_1,T_2,\dots,T_d)\\=-2(n-1)\coef{\mathbf z^{\mathbf c}
\prod _{i=1} ^{d}\prod _{k=0} ^{n}p_{i,k+1}^{m^{(i)}_k}}
\(\prod _{j=1} ^{d}\frac {p_{j,1}} 
{P_{j}\(\frac  
{\om_{1}(\mathbf z)\cdots \om_{d}(\mathbf z)}
{\om_{j}(\mathbf z)}\)}\)
\sum _{j=1} ^{d}\frac {P_{j}\(\frac  
{\om_{1}(\mathbf z)\cdots \om_{d}(\mathbf z)}
{\om_{j}(\mathbf z)}\)} {p_{j,1}}\\+
2(n-1)(d-1)\coef{\mathbf z^{\mathbf c}
\prod _{i=1} ^{d}\prod _{k=0} ^{n}p_{i,k+1}^{m^{(i)}_k}}
\prod _{j=1} ^{d}\frac {p_{j,1}} 
{P_{j}\(\frac  
{\om_{1}(\mathbf z)\cdots \om_{d}(\mathbf z)}
{\om_{j}(\mathbf z)}\)}.
\end{multline}
Now the problem is set up for application of the Lagrange--Good
inversion formula. Let $f_i(\mathbf z)=z_i/P_i(z_1\cdots z_d/z_i)$,
$i=1,2,\dots,d$. If we substitute $f_i(\mathbf z)$ in place of $z_i$,
$i=1,2,\dots,d$, in \eqref{eq:coefD}, and apply Theorem~\ref{thm:LG}
with 
$$g(\mathbf z)=\(\prod _{j=1} ^{d}\frac {p_{j,1}} 
{P_{j}\(\frac  
{z_{1}\cdots z_{d}} {z_{j}}\)}\)
\sum _{j=1} ^{d}\frac {P_{j}\(\frac  
{z_{1}\cdots z_{d}}
{z_{j}}\)} {p_{j,1}},$$
respectively
$$g(\mathbf z)=\prod _{j=1} ^{d}\frac {p_{j,1}} 
{P_{j}\(\frac  
{z_{1}\cdots z_{d}} {z_{j}}\)},$$ 
we obtain that
\begin{multline} \label{eq:NDn}
N_{D_n}(T_1,T_2,\dots,T_d)\\=
-2(n-1)\coef{\mathbf z^{\boldsymbol 0}
\prod _{i=1} ^{d}\prod _{k=0} ^{n}p_{i,k+1}^{m^{(i)}_k}}
\(\prod _{i=1} ^{d}{p_{i,1}} \)
\(\sum _{j=1} ^{d}\frac {z_j} {f_j(\mathbf z)p_{j,1}}\)
\mathbf f^{-\mathbf c}(\mathbf z)
\det_{1\le i,k\le d}
\(\frac {\partial f_i} {\partial z_k}(\mathbf z)\)\\
+2(n-1)(d-1)\coef{\mathbf z^{\boldsymbol 0}
\prod _{i=1} ^{d}\prod _{k=0} ^{n}p_{i,k+1}^{m^{(i)}_k}}
\(\prod _{i=1} ^{d} {p_{i,1}} \)
\mathbf f^{-\mathbf c}(\mathbf z)
\det_{1\le i,k\le d}
\(\frac {\partial f_i} {\partial z_k}(\mathbf z)\),
\end{multline}
where $\boldsymbol 0$ stands for the vector $(0,0,\dots,0)$.
We treat the two terms on the right-hand side of \eqref{eq:NDn}
separately. We begin with the second term:
\begin{align*} 
&\coef{\mathbf z^{\boldsymbol 0}
\prod _{i=1} ^{d}\prod _{k=0} ^{n}p_{i,k+1}^{m^{(i)}_k}}
\(\prod _{i=1} ^{d} {p_{i,1}} \)
\mathbf f^{-\mathbf c}(\mathbf z)
\det_{1\le i,k\le d}
\(\frac {\partial f_i} {\partial z_k}(\mathbf z)\)\\
&=\coef{\mathbf z^{\mathbf c}
\prod _{i=1} ^{d}\prod _{k=0} ^{n}p_{i,k+1}^{m^{(i)}_k}}
\(\prod _{i=1} ^{d} {p_{i,1}} \)
\det_{1\le i,k\le d}
\(\left\{ \begin{matrix} 
P_i^{c_i-1}\(\dfrac {z_1\cdots z_d} {z_i}\),\hfill&i=k\\
-P_i^{c_i-2}\(\dfrac {z_1\cdots z_d} {z_i}\)\hfill\\
\quad \times P_i'\(\dfrac {z_1\cdots z_d} {z_i}\)
\dfrac {z_1\cdots z_d} {z_k},&i\ne k
\end{matrix}\right\}\)\\
&=\coef{\mathbf z^{\mathbf c}
\prod _{i=1} ^{d}\prod _{k=0} ^{n}p_{i,k+1}^{m^{(i)}_k}}
\(\prod _{i=1} ^{d} {p_{i,1}} \)\\
&\kern4cm
\times
\det_{1\le i,k\le d}
\(\left\{ \begin{matrix} 
P_i^{c_i-1}\(\dfrac {z_1\cdots z_d} {z_i}\),\hfill&i=k\\
-\dfrac {1} {c_i-1}\(u\dfrac {d} {du}P_i^{c_i-1}\(u\)\)\Bigg\vert_
{u={z_1\cdots z_d}/ {z_i}},&i\ne k
\end{matrix}\right\}\).
\end{align*}
Reading coefficients, we obtain
\begin{multline*}
\prod _{i=1} ^{d}\binom {c_i-1} 
{{m_1^{(i)},m_2^{(i)},\dots,m_n^{(i)}}}
\det_{1\le i,k\le d}\(\left\{ \begin{matrix} 
1,&i=k\\
-\dfrac {\rk T_i} {c_i-1},&i\ne k
\end{matrix}\right\}\)\\
=\prod _{i=1} ^{d}\binom {c_i-1} 
{{m_1^{(i)},m_2^{(i)},\dots,m_n^{(i)}}}
\det_{1\le i,k\le d}\(
1-\chi(i\ne k)\frac {n-1} {c_i-1}
\),
\end{multline*}
the second line being due to \eqref{eq:crkD}.
Now we can apply Lemma~\ref{lem:1} with $X_i=c_i-1$ and $Y_i=n-1$,
$i=1,2,\dots,d$. The term 
\begin{align*}
\sum _{i=1} ^{d}X_i-\sum _{i=2} ^{d}Y_i
&=\sum _{i=1} ^{d}(c_i-1)-(d-1)(n-1)\\
&=\sum _{i=1} ^{d}(n-\rk T_i-1)-(d-1)(n-1)
\end{align*}
on the right-hand side of \eqref{eq:4} simplifies to $-1$ due to our
assumption concerning the sum of the ranks of the
types $T_i$. Hence, 
if we use the relation \eqref{eq:crkD} once more,
the second term on the right-hand side of \eqref{eq:NDn}
is seen to equal
\begin{equation*} 
-2(d-1)(n-1)^d\prod _{i=1} ^{d}\frac {1} {n-\rk T_i-1}\binom {n-\rk T_i-1} 
{{m_1^{(i)},m_2^{(i)},\dots,m_n^{(i)}}}.
\end{equation*}
This explains the second term in the factor in big parentheses in
\eqref{eq:7comb} and
the fourth term in the factor in big parentheses on
the right-hand side of \eqref{eq:7}.

Finally, we come to the first term on the right-hand side of
\eqref{eq:NDn}. We have
\begin{align*} 
&\coef{\mathbf z^{\boldsymbol 0}
\prod _{i=1} ^{d}\prod _{k=0} ^{n}p_{i,k+1}^{m^{(i)}_k}}
\(\prod _{i=1} ^{d} {p_{i,1}} \)
\frac {z_j} {f_j(\mathbf z)p_{j,1}}
\mathbf f^{-\mathbf c}(\mathbf z)
\det_{1\le i,k\le d}
\(\frac {\partial f_i} {\partial z_k}(\mathbf z)\)\\
&=\coef{\mathbf z^{\mathbf c}
\prod _{i=1} ^{d}\prod _{k=0} ^{n}p_{i,k+1}^{m^{(i)}_k}}
\(\underset{i\ne j}{\prod _{i=1} ^{d}} {p_{i,1}} \)
\det_{1\le i,k\le d}
\(\left\{ \begin{matrix} 
P_i^{c_i-1+\chi(i=j)}\(\dfrac {z_1\cdots z_d} {z_i}\),\hfill&i=k\\
-P_i^{c_i-2+\chi(i=j)}\(\dfrac {z_1\cdots z_d} {z_i}\)\hfill\\
\quad \times P_i'\(\dfrac {z_1\cdots z_d} {z_i}\)
\dfrac {z_1\cdots z_d} {z_k},&i\ne k
\end{matrix}\right\}\)\\
&=\coef{\mathbf z^{\mathbf c}
\prod _{i=1} ^{d}\prod _{k=0} ^{n}p_{i,k+1}^{m^{(i)}_k}}
\(\underset{i\ne j}{\prod _{i=1} ^{d}} {p_{i,1}} \)\\
&\kern1cm
\times
\det_{1\le i,k\le d}
\(\left\{ \begin{matrix} 
P_i^{c_i-1+\chi(i=j)}\(\dfrac {z_1\cdots z_d} {z_i}\),\hfill&i=k\\
-\dfrac {1} {c_i-1+\chi(i=j)}\(u\dfrac {d} {du}P_i^{c_i-1+\chi(i=j)}\(u\)\)\Bigg\vert_
{u={z_1\cdots z_d}/ {z_i}},&i\ne k
\end{matrix}\right\}\).
\end{align*}
Reading coefficients, we obtain
\begin{multline*}
\prod _{i=1} ^{d}\binom {c_i-1+\chi(i=j)} 
{{m_1^{(i)},m_2^{(i)},\dots,m_n^{(i)}}}
\det_{1\le i,k\le d}\(\left\{ \begin{matrix} 
1,&i=k\\
-\dfrac {\rk T_i} {c_i-1+\chi(i=j)},&i\ne k
\end{matrix}\right\}\)\\
=\prod _{i=1} ^{d}\binom {c_i-1+\chi(i=j)} 
{{m_1^{(i)},m_2^{(i)},\dots,m_n^{(i)}}}
\det_{1\le i,k\le d}\(\left\{\begin{matrix}
1-\chi(j\ne k)\frac {n} {c_j},&i=j\\
1-\chi(i\ne k)\frac {n-1} {c_i-1},&i\ne j
\end{matrix}\right\}\),
\end{multline*}
the second line being due to \eqref{eq:crkD}.
Now we can apply Corollary~\ref{cor:1} with $r=j$, $X_i=c_i-1$, 
$i=1,\dots,j-1,j+1,\dots,d$, $X_j=c_j$, $Y=n-1$, and $Z=n$. The term 
\begin{align*}
Z\sum _{i=1} ^{d}X_i+(Y-Z)X_j-(d-1)YZ
&=n\(1+\sum _{i=1} ^{d}(c_i-1)\)- c_j-(d-1)(n-1)n\\
&=n\sum _{i=1} ^{d}(n-\rk T_i-1)+n- c_j-(d-1)(n-1)n
\end{align*}
on the right-hand side of \eqref{eq:4} simplifies to $-c_j$ due to our
assumption concerning the sum of the ranks of the
types $T_i$. Hence, 
if we use the relation \eqref{eq:crkD} once more,
the second term on the right-hand side of \eqref{eq:NDn}
is seen to equal the sum over $j=1,2,\dots,d$ of
\begin{equation*} 
2(n-1)^{d-1}
\binom {n-\rk T_j} 
{{m_1^{(j)},m_2^{(j)},\dots,m_n^{(j)}}}
\underset{i\ne j}{\prod _{i=1} ^{d}}
\frac {1} {n-\rk T_i-1}\binom {n-\rk T_i-1} 
{{m_1^{(i)},m_2^{(i)},\dots,m_n^{(i)}}}.
\end{equation*}
This explains the first terms in the factors in big parentheses on
the right-hand sides of \eqref{eq:7comb} and \eqref{eq:7}.

\medskip
The proof of the theorem is complete.
\end{proof}

Combining the previous theorem with the summation formula of 
Lemma~\ref{lem:binsum}, we can now derive compact formulae
for {\it all\/} type $D_n$ decomposition numbers.

\begin{theorem} \label{thm:22}
{\em(i)} 
Let the types $T_1,T_2,\dots,T_d$ be given,
where 
$$T_i=A_1^{m_1^{(i)}}*A_2^{m_2^{(i)}}*\dots*A_n^{m_n^{(i)}},\quad 
i=1,2,\dots,j-1,j+1,\dots,d,$$ 
and
$$T_j=D_\al*A_1^{m_1^{(j)}}*A_2^{m_2^{(j)}}*\dots*A_n^{m_n^{(j)}},
$$ 
for some $\al\ge2$.
Then
\begin{multline} \label{eq:22a}
N_{D_n}^{\text {comb}}(T_1,T_2,\dots,T_d)=(n-1)^{d-1}
\binom {n-1}{\rk T_1+\rk T_2+\dots+\rk T_d-1}\\
\times
\binom {n-\rk T_j}{m_1^{(j)},m_2^{(j)},\dots,
m_n^{(j)}}
\underset{i\ne j}{\prod _{i=1} ^{d}}
\frac {1} {n-\rk T_i-1}\binom {n-\rk T_i-1}{m_1^{(i)},m_2^{(i)},\dots,
m_n^{(i)}},
\end{multline}
where the multinomial coefficient is defined as in
Lemma~{\em\ref{lem:binsum}}.
For $\al\ge4,$ the number  
$N_{D_n}(T_1,T_2,\dots, T_d)$ is given by the same formula.

\smallskip
{\em(ii)} 
Let the types $T_1,T_2,\dots,T_d$ be given,
where 
$$T_i=A_1^{m_1^{(i)}}*A_2^{m_2^{(i)}}*\dots*A_n^{m_n^{(i)}},\quad 
i=1,2,\dots,d.$$ 
Then
\begin{multline} \label{eq:22b}
N_{D_n}^{\text {comb}}(T_1,T_2,\dots,T_d)=(n-1)^{d-1}
\binom {n-1}{\rk T_1+\rk T_2+\dots+\rk T_d-1}\\
\times
\(2\sum _{j=1} ^{d}
\binom {n-\rk T_j}{m_1^{(j)},m_2^{(j)},\dots,m_n^{(j)}}
\Bigg(\underset{i\ne j}
{\prod _{i=1} ^{d}}
\frac {1} {n-\rk T_i-1}\binom {n-\rk T_i-1}{m_1^{(i)},m_2^{(i)},\dots,
m_n^{(i)}}\Bigg)\right.
\\
+\(\frac {\(n-\sum _{\ell=1} ^{d}\rk T_\ell\)\(n-1-\sum _{\ell=1} ^{d}\rk
T_\ell\)} 
{\sum _{\ell=1} ^{d}\rk T_\ell}-2(d-2)(n-1)\)\\
\left.
\cdot\underset{\vphantom{f}}{\prod _{i=1} ^{d}}
\frac {1} {n-\rk T_i-1}\binom {n-\rk T_i-1}{m_1^{(i)},m_2^{(i)},\dots,
m_n^{(i)}}\right),
\end{multline}
whereas
\begin{multline} \label{eq:22c}
N_{D_n}(T_1,T_2,\dots,T_d)=(n-1)^{d-1}
\binom {n-1}{\rk T_1+\rk T_2+\dots+\rk T_d-1}\\
\times
\(\sum _{j=1} ^{d}\Bigg(\underset{i\ne j}
{\prod _{i=1} ^{d}}
\frac {1} {n-\rk T_i-1}\binom {n-\rk T_i-1}{m_1^{(i)},m_2^{(i)},\dots,
m_n^{(i)}}\Bigg)
\Bigg(
2\binom {n-\rk T_j}{m_1^{(j)},m_2^{(j)},\dots,m_n^{(j)}}\right.\\
+\binom {n-\rk T_j}{m_1^{(j)},m_2^{(j)},m_3^{(j)}-1,m_4^{(j)},\dots,m_n^{(j)}}
+\binom {n-\rk T_j}{m_1^{(j)}-2,m_2^{(j)},\dots,m_n^{(j)}}
\Bigg)\\
+\(\frac {\(n-\sum _{\ell=1} ^{d}\rk T_\ell\)\(n-1-\sum _{\ell=1} ^{d}\rk
T_\ell\)} 
{\sum _{\ell=1} ^{d}\rk T_\ell}-2(d-2)(n-1)\)\\
\left.
\cdot\underset{\vphantom{f}}{\prod _{i=1} ^{d}}
\frac {1} {n-\rk T_i-1}\binom {n-\rk T_i-1}{m_1^{(i)},m_2^{(i)},\dots,
m_n^{(i)}}\right).
\end{multline}

\smallskip
{\em(iii)} All other decomposition numbers $N_{D_n}(T_1,T_2,\dots,T_d)$
and $N_{D_n}^{\text {comb}}(T_1,T_2,\dots,T_d)$ are zero.
\end{theorem}

\begin{remark} 
The caveats on interpretations of the formulae in Theorem~\ref{thm:3}
for critical choices of the parameters (cf.\ 
the Remark after the statement of that theorem) apply also to the
formulae of Theorem~\ref{thm:22}.
\end{remark}

\begin{proof} 
We proceed in a manner similar to the proof of Theorem~\ref{thm:21}.
If we write $r$ for $n-\rk T_1-\rk T_2-\dots-\rk T_d$ and set
$\Phi=D_n$, then the relation \eqref{Ab} becomes
\begin{equation} \label{eq:relD}
N_{D_n}(T_1,T_2,\dots,T_d)=
{\sum _{T:\rk T=r}}
^{}N_{D_n}(T_1,T_2,\dots,T_d,T),
\end{equation}
with the same relation holding for $N_{D_n}^{\text {comb}}$ in place
of $N_{D_n}$.

In order to prove \eqref{eq:22a}, we
let $T=A_1^{m_1}*A_2^{m_2}*\dots*A_n^{m_n}$ and use
\eqref{eq:6} in \eqref{eq:relD}, to obtain
\begin{multline*}
N_{D_n}^{\text {comb}}(T_1,T_2,\dots,T_d)=\sum _{m_1+2m_2+\dots+nm_n=r} ^{}
(n-1)^{d}\frac {1} {n-r-1}\binom {n-r-1}{m_1,m_2,\dots,m_n}\\
\cdot
\binom {n-\rk T_j}{m_1^{(j)},m_2^{(j)},\dots,
m_n^{(j)}}
\underset{i\ne j}{\prod _{i=1} ^{d}}
\frac {1} {n-\rk T_i-1}\binom {n-\rk T_i-1}{m_1^{(i)},m_2^{(i)},\dots,
m_n^{(i)}}.
\end{multline*}
If we use \eqref{eq:binsum} with $M=n-r-1$, we arrive at our claim after
little simplification.

Next we prove \eqref{eq:22b}. In contrast to the previous argument,
here the summation on the
right-hand side of \eqref{eq:relD} must be taken over all types $T$
of the form
$T=D_\al*A_1^{m_1}*A_2^{m_2}*\dots*A_n^{m_n}$, $\al\ge2$, 
as well as of the form
$T=A_1^{m_1}*A_2^{m_2}*\dots*A_n^{m_n}$.
For the sum over the former types, we have to
substitute \eqref{eq:6} in \eqref{eq:relD}, to get
\begin{multline} \label{eq:sumD1}
\sum _{\al=2} ^{n}\sum _{m_1+2m_2+\dots+nm_n=r-\al} ^{}
(n-1)^{d}\binom {n-r}{m_1,m_2,\dots,m_n}\\
\cdot
{\prod _{i=1} ^{d}}
\frac {1} {n-\rk T_i-1}\binom {n-\rk T_i-1}{m_1^{(i)},m_2^{(i)},\dots,
m_n^{(i)}}.
\end{multline}
On the other hand, for the sum over the latter types, we have to
substitute \eqref{eq:7comb} in \eqref{eq:relD}, to get
\begin{multline} \label{eq:sumD2}
2\sum _{m_1+2m_2+\dots+nm_n=r} ^{}
(n-1)^{d}\binom {n-r}{m_1,m_2,\dots,m_n}
{\prod _{i=1} ^{d}}
\frac {1} {n-\rk T_i-1}\binom {n-\rk T_i-1}{m_1^{(i)},m_2^{(i)},\dots,
m_n^{(i)}}
\\
+\sum _{m_1+2m_2+\dots+nm_n=r} ^{}
(n-1)^{d}\frac {1} {n-r-1}\binom {n-r-1}{m_1,m_2,\dots,m_n}
\kern4cm\\
\cdot
\(2\sum _{j=1} ^{d}
\binom {n-\rk T_j}{m_1^{(j)},m_2^{(j)},\dots,m_n^{(j)}}
\underset{i\ne j} {\prod _{i=1} ^{d}}
\frac {1} {n-\rk T_i-1}\binom {n-\rk T_i-1}{m_1^{(i)},m_2^{(i)},\dots,
m_n^{(i)}}\right.\\
\left.
-2(d-1)(n-1)\underset{\vphantom{f}}{\prod _{i=1} ^{d}}
\frac {1} {n-\rk T_i-1}\binom {n-\rk T_i-1}{m_1^{(i)},m_2^{(i)},\dots,
m_n^{(i)}}\).
\end{multline}
We simplify \eqref{eq:sumD1} by
using \eqref{eq:binsum} with $r$ replaced by $r-\al$ and $M=n-r$, 
and by subsequently applying the elementary summation formula
\begin{equation} \label{eq:rksum1}
\sum _{\al=2} ^{n}\binom {n-\al-1}{r-\al}=
\sum _{\al=2} ^{n}\binom {n-\al-1}{n-r-1}=\binom {n-2}{n-r}=\binom
{n-2}{r-2}.
\end{equation}
The expression which we obtain in this way 
explains the fraction in the third line of
\eqref{eq:22b} multiplied by the expression in the last line.
On the other hand, we simplify the sums in \eqref{eq:sumD2} by
using \eqref{eq:binsum} with $M=n-r$, respectively $M=n-r-1$.
Thus, the expression \eqref{eq:sumD2} becomes
\begin{multline*}
2
(n-1)^{d}\binom {n-1}r
{\prod _{i=1} ^{d}}
\frac {1} {n-\rk T_i-1}\binom {n-\rk T_i-1}{m_1^{(i)},m_2^{(i)},\dots,
m_n^{(i)}}
\\
+
(n-1)^{d-1}\binom {n-1}r
\Bigg(2\sum _{j=1} ^{d}
\binom {n-\rk T_j}{m_1^{(j)},m_2^{(j)},\dots,m_n^{(j)}}
\underset{i\ne j} {\prod _{i=1} ^{d}}
\frac {1} {n-\rk T_i-1}\binom {n-\rk T_i-1}{m_1^{(i)},m_2^{(i)},\dots,
m_n^{(i)}}\\
-2(d-1)(n-1)\underset{\vphantom{f}}{\prod _{i=1} ^{d}}
\frac {1} {n-\rk T_i-1}\binom {n-\rk T_i-1}{m_1^{(i)},m_2^{(i)},\dots,
m_n^{(i)}}\Bigg),
\end{multline*}
which explains the expression in the second line of \eqref{eq:22b}
and the second expression in the third line of
\eqref{eq:22b} multiplied by the expression in the last line.

The proof of \eqref{eq:22c} is analogous, using \eqref{eq:7} instead
of \eqref{eq:7comb}. We leave the details to the reader. 
\end{proof}

\section{Generalised non-crossing partitions}
\label{sec:6} 

In this section we recall the definition of Armstrong's \cite{ArmDAA}
generalised non-crossing partitions poset, and its combinatorial
realisation from \cite{ArmDAA}
and \cite{KratCG} for the types $A_n$, $B_n$, and $D_n$.

Let again $\Phi$ be a finite root system of rank $n$, and let $W=W(\Phi)$
be the corresponding reflection group.
We define first the {\it non-crossing partition lattice
$NC(\Phi)$} (cf.\ \cite{BesDAA,BRWaAA}). 
Let $c$ be a {\it Coxeter element\/} in $W$.
Then $NC(\Phi)$ is defined to be
the restriction of the partial order $\le_T$ from Section~\ref{sec:1} 
to the set
of all elements which are less than or equal to $c$ in this partial order. 
This definition makes sense since any two Coxeter elements in $W$
are conjugate to each other; the induced inner
automorphism then restricts to an isomorphism of the posets
corresponding to the two Coxeter elements.
It can be shown that $NC(\Phi)$ is in fact a lattice (see
\cite{BRWaAB} for a uniform proof), and moreover self-dual
(this is obvious from the definition). Clearly, the minimal element
in $NC(\Phi)$ is the identity element in $W$, which we denote by
$\ep$, and the maximal element
in $NC(\Phi)$ is the chosen Coxeter element $c$.
The term ``non-crossing partition lattice" is used because
$NC(A_n)$ is isomorphic to the lattice of non-crossing partitions
of $\{1,2,\dots,n+1\}$,
originally introduced by Kreweras \cite{KrewAC}
(see also \cite{FoReAB} and below), and since also
$NC(B_n)$ and $NC(D_n)$ can be realised as lattices of
non-crossing partitions (see \cite{AtReAA,ReivAG} and below).

In addition to a fixed root system, the definition of
Armstrong's {\it generalised
non-crossing partitions} requires a fixed positive integer $m$.
The poset of {\it $m$-divisible non-crossing partitions associated to
the root system $\Phi$} has as 
ground-set the following subset of $(NC(\Phi))^{m+1}$:
\begin{multline} \label{eq:01}
NC^m(\Phi)=\big\{(w_0;w_1,\dots,w_m):w_0w_1\cdots w_m=c\text{ and }\\
\ell_T(w_0)+\ell_T(w_1)+\dots+\ell_T(w_m)=\ell_T(c)\big\}.
\end{multline}
The order relation is defined by
$$(u_0;u_1,\dots,u_m)\le(w_0;w_1,\dots,w_m)\quad \text{if and only
if}\quad u_i\ge_T w_i,\ 1\le i\le m.$$
(According to this definition, $u_0$ and $w_0$ need
not be related in any way. However, it follows from 
\cite[Lemma~3.4.7]{ArmDAA} that, in fact, $u_0\le_T w_0$.) 
The poset $NC^m(\Phi)$ is graded by 
the rank function
\begin{equation} \label{eq:rkm}
\rk\big((w_0;w_1,\dots,w_m)\big)=\ell_T(w_0).
\end{equation}
Thus, there is a unique maximal element, namely $(c;\ep,\dots,\ep)$,
where $\ep$ stands for the identity element in $W$, but, for $m>1$, there 
are many
different minimal elements. In particular, $NC^m(\Phi)$ has no least element
if $m>1$; hence, $NC^m(\Phi)$ is not a lattice
for $m>1$. (It is, however, a graded join-semilattice, see
\cite[Theorem~3.4.4]{ArmDAA}.) 

In what follows, we shall use the notions ``generalised
non-crossing partitions" and ``$m$-divisible non-crossing partitions"
interchangeably, where the latter notion will be employed particularly
in contexts in which we want to underline the presence of the
parameter $m$.

\begin{figure}
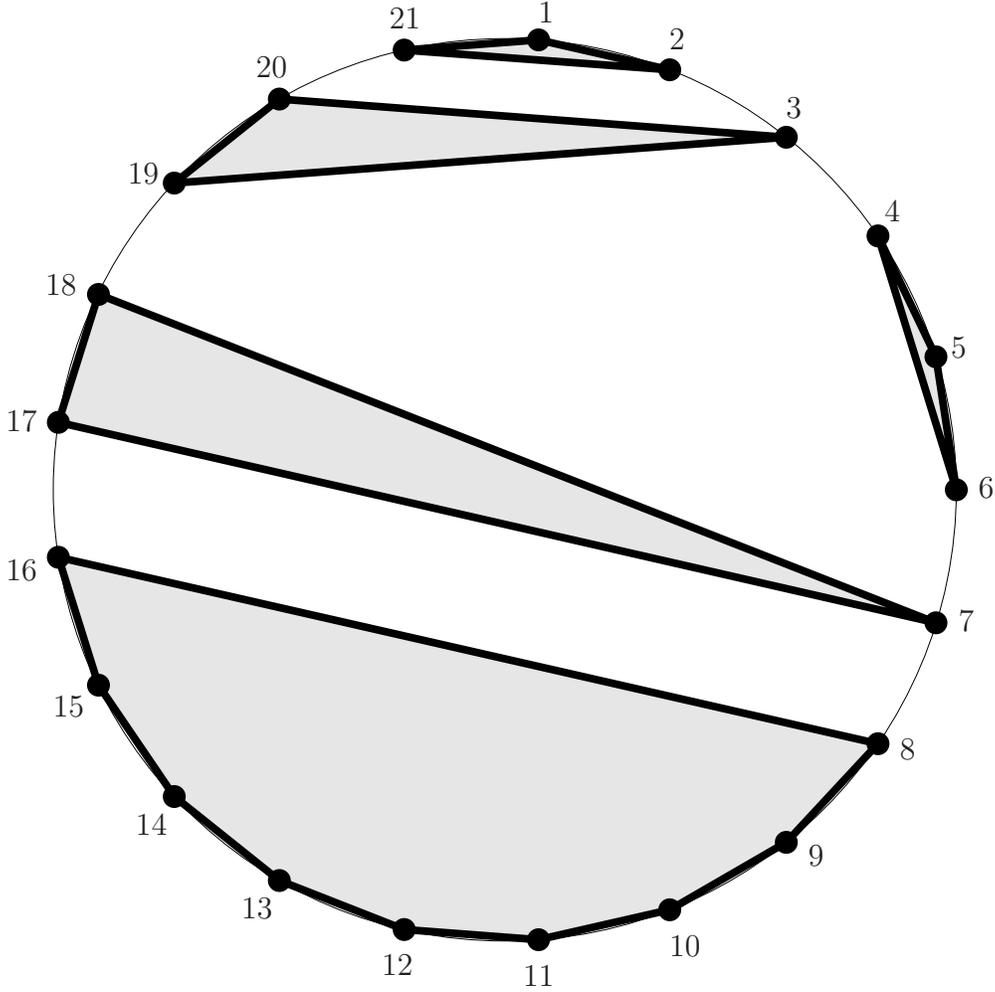

\centertexdraw{
  \drawdim truecm  \linewd.01
  \lcir r:6

  \linewd.10

  \htext(0.448381 6.22322){1} 
  \htext(2.19205 5.86524){2} 
  \htext(3.74094 4.96099){3}
  \htext(5.05743 3.57992){4} 
  \htext(5.93344 1.76853){5} 
  \htext(6.3 -0.1) {6}
  \htext(6.03344 -1.86853){7}
  \htext(5.25743 -3.57992){8} 
  \htext(4.04094 -4.99099){9}
  \htext(2.19205 -6.18524){10} 
  \htext(0.248381 -6.58322){11} 
  \htext(-1.63513 -6.44957){12}
  \htext(-3.5 -5.69615) {13}
  \htext(-4.89831 -4.58104){14} 
  \htext(-6.00581 -3.0033){15}
  \htext(-6.63298 -1.194254){16} 
  \htext(-6.63298 0.794254){17} 
  \htext(-6.10581 2.6033){18}
  \htext(-5.00831 4.08104){19} 
  \htext(-3.3 5.49615) {20}
  \htext(-1.53513 6.14957){21}

  \move(0.448381 5.98322)
  \lvec(2.19205 5.58524)
  \lvec(-1.33513 5.84957)
  \lvec(0.448381 5.98322)
  \lfill f:.9

  \move(3.74094 4.69099)
  \lvec(-4.39831 4.08104)
  \lvec(-3. 5.19615)   
  \lvec(3.74094 4.69099)
  \lfill f:.9

  \move(4.95743 3.37992)
  \lvec(5.73344 1.76853)
  \lvec(6. 0)   
  \lvec(4.95743 3.37992)
  \lfill f:.9

  \move(5.73344 -1.76853)
  \lvec(-5.93298 0.894254)
  \lvec(-5.40581 2.6033)
  \lvec(5.73344 -1.76853)
  \lfill f:.9

  \move(4.95743 -3.37992) 
  \lvec(3.74094 -4.69099) 
  \lvec(2.19205 -5.58524) 
  \lvec(0.448381 -5.98322) 
  \lvec(-1.33513 -5.84957) 
  \lvec(-3. -5.19615)      
  \lvec(-4.39831 -4.08104) 
  \lvec(-5.40581 -2.6033)  
  \lvec(-5.93298 -0.894254) 
  \lvec(4.95743 -3.37992)
  \lfill f:.9

  \Ringerl(0.448381 5.98322) 
  \Ringerl(2.19205 5.58524) 
  \Ringerl(3.74094 4.69099) 
  \Ringerl(4.95743 3.37992) 
  \Ringerl(5.73344 1.76853) 
  \Ringerl(6. 0)            
  \Ringerl(5.73344 -1.76853) 
  \Ringerl(4.95743 -3.37992) 
  \Ringerl(3.74094 -4.69099) 
  \Ringerl(2.19205 -5.58524) 
  \Ringerl(0.448381 -5.98322) 
  \Ringerl(-1.33513 -5.84957) 
  \Ringerl(-3. -5.19615)      
  \Ringerl(-4.39831 -4.08104) 
  \Ringerl(-5.40581 -2.6033)  
  \Ringerl(-5.93298 -0.894254) 
  \Ringerl(-5.93298 0.894254) 
  \Ringerl(-5.40581 2.6033)   
  \Ringerl(-4.39831 4.08104)  
  \Ringerl(-3. 5.19615)       
  \Ringerl(-1.33513 5.84957)  
}
\caption{Combinatorial realisation of a 
$3$-divisible non-crossing partition of type $A_{6}$}
\label{fig:7}
\end{figure}

\medskip
In the remainder of this section, we explain combinatorial
realisations of the $m$-divisible non-crossing partitions of types
$A_{n-1}$, $B_n$, and $D_n$. In order to be able to do so, we need to
recall the definition of Kreweras' non-crossing partitions of
$\{1,2,\dots,N\}$, his
{\it``partitions non crois\'ees d'un cycle"} of \cite{KrewAC}. 
We place $N$ vertices around
a cycle, and label them $1,2,\dots,N$ in clockwise order. The
circular representation of a partition of the set $\{1,2,\dots,N\}$
is the geometric object which arises by representing each block
$\{i_1,i_2,\dots,i_k\}$ of the partition, where $i_1<i_2<\dots<i_k$, by the
polygon consisting of the vertices labelled $i_1,i_2,\dots,i_k$ and
edges which connect these vertices in clockwise order. A partition of
$\{1,2,\dots,N\}$ is called {\it non-crossing} if any two edges 
in its circular representation are disjoint.
Figure~\ref{fig:7} shows the non-crossing partition 
$$\{\{1,2,21\},
\{3,19,20\},
\{4,5,6\},
\{7,17,18\},
\{8,9,10,11,12,13,14,15,16\}\}$$
of $\{1,2,\dots,21\}$. There is a natural partial order on Kreweras'
non-crossing partitions defined by refinement: a partition $\pi_1$ is
less than or equal to the partition $\pi_2$ if every block of
$\pi_1$ is contained in some block of $\pi_2$.

If $\Phi=A_{n-1}$, the $m$-divisible non-crossing partitions are in
bijection with Kreweras-type non-crossing partitions of the set
$\{1,2,\dots,mn\}$, in which all the block sizes are divisible by
$m$. We denote the latter set of non-crossing partitions by
$\widetilde {NC}{}^m(A_{n-1})$. It has been first considered by
Edelman in \cite{EdelAA}.
In fact, Figure~\ref{fig:7} shows an example of a $3$-divisible
non-crossing partition of type $A_{20}$. 

Given an element $(w_0;w_1,\dots,w_m)\in NC^m(A_{n-1})$,
the bijection, $\Na^m_{A_{n-1}}$ say, 
from \cite[Theorem~4.3.8]{ArmDAA} 
works by ``blowing up"
$w_1,w_2,\dots,w_m$, thereby ``interleaving" them,
and then ``gluing" them together by an operation which is called {\it
Kreweras complement\/} in \cite{ArmDAA}. More precisely, for
$i=1,2,\dots,m$, let $\ta_{m,i}$ be
the transformation which maps a permutation $w\in S_{n}$ to a
permutation $\ta_{m,i}(w)\in S_{mn}$ by letting
$$(\ta_{m,i}(w))(mk+i-m)=mw(k)+i-m,\quad k=1,2,\dots,n,$$ 
and 
$(\ta_{m,i}(w))(l)=l$ for all $l\not\equiv i$~(mod~$m$).
At this point, the reader should recall from Section~\ref{sec:1} that
$W(A_{n-1})$ is the symmetric group $S_n$, and that the standard
choice of a Coxeter element in $W(A_{n-1})=S_n$ is $c=(1,2,\dots,n)$. 
With this choice of Coxeter element, the announced bijection maps 
$(w_0;w_1,\dots,w_m)\in NC^m(A_{n-1})$ to
$$\Na^m_{A_{n-1}}(w_0;w_1,\dots,w_m)=
(1,2,\dots,mn)\,(\ta_{m,1}(w_1))^{-1}\,(\ta_{m,2}(w_2))^{-1}\,\cdots\,
(\ta_{m,m}(w_m))^{-1}.$$
We refer the reader to \cite[Sec.~4.3.2]{ArmDAA} for the details.
For example, let $n=7$, $m=3$, $w_0=(4,5,6)$,
$w_1=(3,6)$,
$w_2=(1,7)$, and
$w_3=(1,2,6)$. Then $(w_0;w_1,w_2,w_3)$ is mapped to
\begin{align} \notag 
\Na^3_{A_{6}}(w_0;w_1,w_2,w_3)&=(1,2,\dots,21)\,(7,16)\,
(2,20)\,(18,6,3)\\
&=(1,2,21)\,(3,19,20)\,(4,5,6)\,(7,17,18)\,(8,9,\dots,16).
\label{eq:NCA}
\end{align} 
Figure~\ref{fig:7} shows the graphical representation of
\eqref{eq:NCA} on the circle, in which we represent a cycle
$(i_1,i_2,\dots,i_k)$ as a 
polygon consisting of the vertices labelled $i_1,i_2,\dots,i_k$ and
edges which connect these vertices in clockwise order.

It is shown in \cite[Theorem~4.3.8]{ArmDAA} that $\Na^m_{A_{n-1}}$ 
is in fact an isomorphism between the {\it posets} $NC^m(A_{n-1})$ and 
$\widetilde{NC}{}^m(A_{n-1})$. Furthermore,
it is proved in \cite[Theorem~4.3.13]{ArmDAA} that
\begin{equation} \label{eq:BlockA}
c_i(w_0)=b_{i}(\Na^m_{A_{n-1}}(w_0;w_1,\dots,w_m)),\quad
i=1,2,\dots,n,
\end{equation}
where $c_i(w_0)$ denotes the number of cycles of length $i$ of $w_0$ 
and $b_{i}(\pi)$ denotes the number of blocks of size $mi$ in 
the non-crossing partition $\pi$. 

\begin{figure}
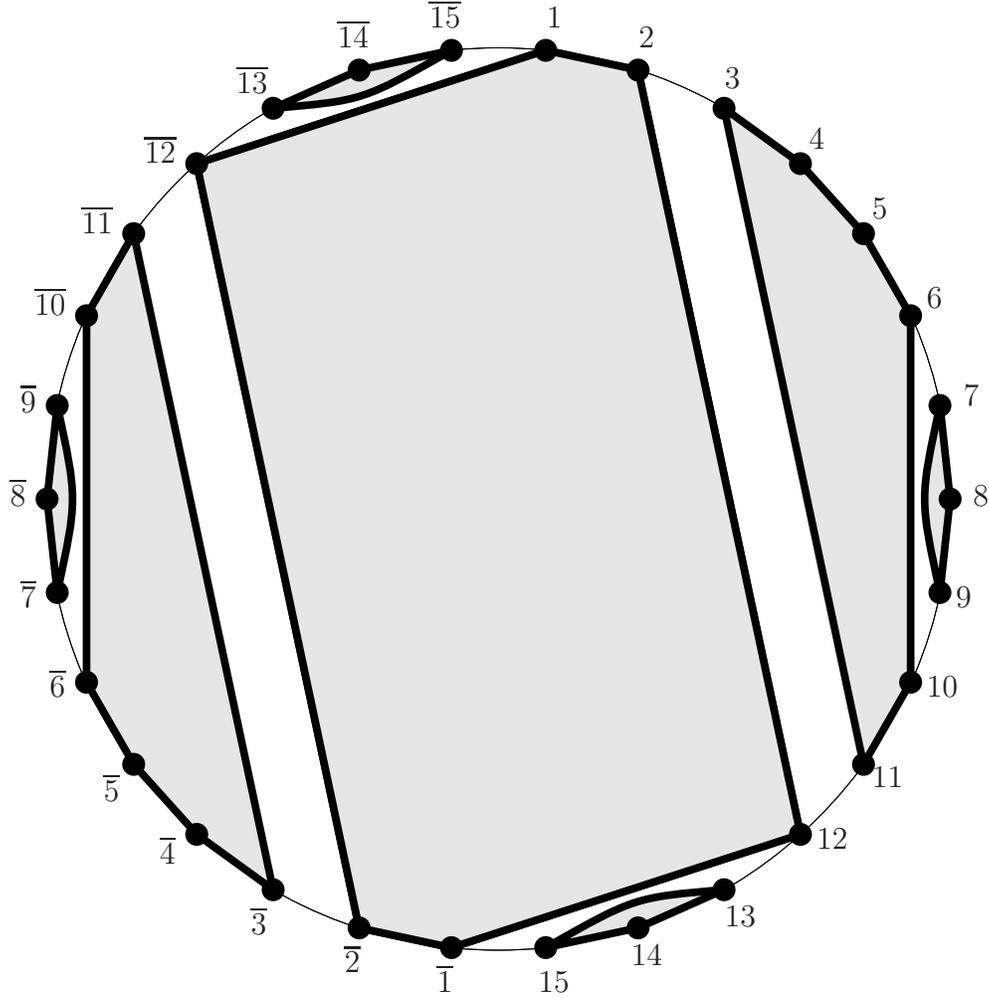

\centertexdraw{
  \drawdim truecm  \linewd.01
  \lcir r:6

  \linewd.10

  \htext(0.627171 6.26713){$1$} 
  \htext(1.8541 6.00634){$2$}   
  \htext(3. 5.45615){$3$}       
  \htext(4.11478 4.65887){$4$} 
  \htext(4.9541 3.72671){$5$}  
  \htext(5.68127 2.54042){$6$} 
  \htext(6.16889 1.24747){$7$} 
  \htext(6.3 -.1){$8$}            
  \htext(6.06889 -1.44747){$9$} 
  \htext(5.68127 -2.64042){$10$} 
  \htext(4.9541 -3.83671){$11$}  
  \htext(4.21478 -4.65887){$12$} 
  \htext(3. -5.69615){$13$}      
  \htext(1.7541 -6.20634){$14$}  
  \htext(0.517171 -6.56713){$15$} 
  \htext(-0.827171 -6.56713){$\overline{1}$} 
  \htext(-2.0541 -6.30634){$\overline{2}$}  
  \htext(-3.3 -5.79615){$\overline{3}$}      
  \htext(-4.51478 -4.85887){$\overline{4}$} 
  \htext(-5.2541 -4.02671){$\overline{5}$}  
  \htext(-5.98127 -2.64042){$\overline{6}$} 
  \htext(-6.36889 -1.44747){$\overline{7}$} 
  \htext(-6.5 -.1){$\overline{8}$}             
  \htext(-6.36889 1.14747){$\overline{9}$}  
  \htext(-6.18127 2.44042){$\overline{10}$}  
  \htext(-5.5541 3.52671){$\overline{11}$}   
  \htext(-4.71478 4.45887){$\overline{12}$}  
  \htext(-3.5 5.39615){$\overline{13}$}       
  \htext(-2.1541 6.00634){$\overline{14}$}   
  \htext(-0.927171 6.26713){$\overline{15}$} 

  \move(0.627171 5.96713)
  \lvec(1.8541 5.70634)
  \lvec(4.01478 -4.45887)
  \lvec(-0.627171 -5.96713)
  \lvec(-1.8541 -5.70634)
  \lvec(-4.01478 4.45887)
  \lvec(0.627171 5.96713)
  \lfill f:.9

  \move(3. 5.19615) 
  \lvec(4.01478 4.45887)
  \lvec(4.8541 3.52671)
  \lvec(5.48127 2.44042)
  \lvec(5.48127 -2.44042)
  \lvec(4.8541 -3.52671)
  \lvec(3. 5.19615) 
  \lfill f:.9

  \move(5.86889 1.24747)
  \lvec(6. 0)  
  \lvec(5.86889 -1.24747)
  \clvec(5.6 0)(5.6 0)(5.86889 1.24747)
  \lfill f:.9

  \move(3. -5.19615)
  \lvec(1.8541 -5.70634)
  \lvec(0.627171 -5.96713)
  \clvec(1.8 -5.3)(1.8 -5.3)(3. -5.19615)
  \lfill f:.9

  \move(-0.627171 -5.96713)
  \lvec(-1.8541 -5.70634)
  \lvec(-4.01478 4.45887)
  \lvec(0.627171 5.96713)
  \lvec(1.8541 5.70634)
  \lvec(4.01478 -4.45887)
  \lvec(-0.627171 -5.96713)
  \lfill f:.9

  \move(-3. -5.19615) 
  \lvec(-4.01478 -4.45887)
  \lvec(-4.8541 -3.52671)
  \lvec(-5.48127 -2.44042)
  \lvec(-5.48127 2.44042)
  \lvec(-4.8541 3.52671)
  \lvec(-3. -5.19615) 
  \lfill f:.9

  \move(-5.86889 -1.24747)
  \lvec(-6. 0)  
  \lvec(-5.86889 1.24747)
  \clvec(-5.6 0)(-5.6 0)(-5.86889 -1.24747)
  \lfill f:.9

  \move(-3. 5.19615)
  \lvec(-1.8541 5.70634)
  \lvec(-0.627171 5.96713)
  \clvec(-1.8 5.3)(-1.8 5.3)(-3. 5.19615)
  \lfill f:.9

  \move(0 0)
  \linewd.01
  \lcir r:6

  \Ringerl(0.627171 5.96713) 
  \Ringerl(1.8541 5.70634)   
  \Ringerl(3. 5.19615)       
  \Ringerl(4.01478 4.45887) 
  \Ringerl(4.8541 3.52671)  
  \Ringerl(5.48127 2.44042) 
  \Ringerl(5.86889 1.24747) 
  \Ringerl(6. 0)            
  \Ringerl(5.86889 -1.24747) 
  \Ringerl(5.48127 -2.44042) 
  \Ringerl(4.8541 -3.52671)  
  \Ringerl(4.01478 -4.45887) 
  \Ringerl(3. -5.19615)      
  \Ringerl(1.8541 -5.70634)  
  \Ringerl(0.627171 -5.96713) 
  \Ringerl(-0.627171 -5.96713) 
  \Ringerl(-1.8541 -5.70634)  
  \Ringerl(-3. -5.19615)      
  \Ringerl(-4.01478 -4.45887) 
  \Ringerl(-4.8541 -3.52671)  
  \Ringerl(-5.48127 -2.44042) 
  \Ringerl(-5.86889 -1.24747) 
  \Ringerl(-6. 0)             
  \Ringerl(-5.86889 1.24747)  
  \Ringerl(-5.48127 2.44042)  
  \Ringerl(-4.8541 3.52671)   
  \Ringerl(-4.01478 4.45887)  
  \Ringerl(-3. 5.19615)       
  \Ringerl(-1.8541 5.70634)   
  \Ringerl(-0.627171 5.96713) 

}
\caption{A $3$-divisible non-crossing partition of type $B_{5}$ with
zero block}
\label{fig:8a}
\end{figure}

\medskip
If $\Phi=B_n$, the $m$-divisible non-crossing partitions are in
bijection with Kreweras-type
non-crossing partitions $\pi$ of the set
$\{1,2,\dots,mn,\bar 1,\bar 2,\dots,\overline{mn}\}$, 
in which all the block sizes are divisible by
$m$, and which have the property that if $B$ is a block of $\pi$ then
also $\overline B:=\{\bar x:x\in B\}$ is a block of $\pi$. (Here, as
earlier, we adopt the convention that $\bar{\bar x}=x$ for all $x$.)
We denote the latter set of non-crossing partitions by
$\widetilde {NC}{}^m(B_{n})$.
A block $B$ with $\overline B=B$ is called a {\it zero block\/}.
A non-crossing partition in $\widetilde {NC}{}^m(B_{n})$ can only have
{\it at most\/} one zero block.
Figures~\ref{fig:8} and \ref{fig:8a} give examples of $3$-divisible
non-crossing partitions of type $B_{5}$. Figure~\ref{fig:8} shows
one without a zero block, while Figure~\ref{fig:8a} shows one with a
zero block. Clearly, the condition that $B$ is a block of the
partition if and only if
$\overline B$ is a block translates into the condition that the
geometric realisation of the partition is invariant under rotation by
$180^\circ$.

\begin{figure}
\centertexdraw{
  \drawdim truecm  \linewd.01
  \lcir r:6

  \linewd.10

  \htext(0.627171 6.26713){$1$} 
  \htext(1.8541 6.00634){$2$}   
  \htext(3. 5.45615){$3$}       
  \htext(4.11478 4.65887){$4$} 
  \htext(4.9541 3.72671){$5$}  
  \htext(5.68127 2.54042){$6$} 
  \htext(6.16889 1.24747){$7$} 
  \htext(6.3 -.1){$8$}            
  \htext(6.06889 -1.44747){$9$} 
  \htext(5.68127 -2.64042){$10$} 
  \htext(4.9541 -3.83671){$11$}  
  \htext(4.21478 -4.65887){$12$} 
  \htext(3. -5.69615){$13$}      
  \htext(1.7541 -6.20634){$14$}  
  \htext(0.517171 -6.56713){$15$} 
  \htext(-0.827171 -6.56713){$\overline{1}$} 
  \htext(-2.0541 -6.30634){$\overline{2}$}  
  \htext(-3.3 -5.79615){$\overline{3}$}      
  \htext(-4.51478 -4.85887){$\overline{4}$} 
  \htext(-5.2541 -4.02671){$\overline{5}$}  
  \htext(-5.98127 -2.64042){$\overline{6}$} 
  \htext(-6.36889 -1.44747){$\overline{7}$} 
  \htext(-6.5 -.1){$\overline{8}$}             
  \htext(-6.36889 1.14747){$\overline{9}$}  
  \htext(-6.18127 2.44042){$\overline{10}$}  
  \htext(-5.5541 3.52671){$\overline{11}$}   
  \htext(-4.71478 4.45887){$\overline{12}$}  
  \htext(-3.5 5.39615){$\overline{13}$}       
  \htext(-2.1541 6.00634){$\overline{14}$}   
  \htext(-0.927171 6.26713){$\overline{15}$} 

  \move(0.627171 5.96713)
  \lvec(-1.8541 -5.70634)
  \lvec(-4.01478 4.45887)
  \lvec(0.627171 5.96713)
  \lfill f:.9

  \move(1.8541 5.70634)
  \lvec(4.01478 -4.45887)
  \lvec(-0.627171 -5.96713)
  \lvec(1.8541 5.70634)
  \lfill f:.9

  \move(3. 5.19615) 
  \lvec(4.01478 4.45887)
  \lvec(4.8541 3.52671)
  \lvec(5.48127 2.44042)
  \lvec(5.48127 -2.44042)
  \lvec(4.8541 -3.52671)
  \lvec(3. 5.19615) 
  \lfill f:.9

  \move(5.86889 1.24747)
  \lvec(6. 0)  
  \lvec(5.86889 -1.24747)
  \clvec(5.6 0)(5.6 0)(5.86889 1.24747)
  \lfill f:.9

  \move(3. -5.19615)
  \lvec(1.8541 -5.70634)
  \lvec(0.627171 -5.96713)
  \clvec(1.8 -5.3)(1.8 -5.3)(3. -5.19615)
  \lfill f:.9

  \move(-3. -5.19615) 
  \lvec(-4.01478 -4.45887)
  \lvec(-4.8541 -3.52671)
  \lvec(-5.48127 -2.44042)
  \lvec(-5.48127 2.44042)
  \lvec(-4.8541 3.52671)
  \lvec(-3. -5.19615) 
  \lfill f:.9

  \move(-5.86889 -1.24747)
  \lvec(-6. 0)  
  \lvec(-5.86889 1.24747)
  \clvec(-5.6 0)(-5.6 0)(-5.86889 -1.24747)
  \lfill f:.9

  \move(-3. 5.19615)
  \lvec(-1.8541 5.70634)
  \lvec(-0.627171 5.96713)
  \clvec(-1.8 5.3)(-1.8 5.3)(-3. 5.19615)
  \lfill f:.9

  \move(0 0)
  \linewd.01
  \lcir r:6

  \Ringerl(0.627171 5.96713) 
  \Ringerl(1.8541 5.70634)   
  \Ringerl(3. 5.19615)       
  \Ringerl(4.01478 4.45887) 
  \Ringerl(4.8541 3.52671)  
  \Ringerl(5.48127 2.44042) 
  \Ringerl(5.86889 1.24747) 
  \Ringerl(6. 0)            
  \Ringerl(5.86889 -1.24747) 
  \Ringerl(5.48127 -2.44042) 
  \Ringerl(4.8541 -3.52671)  
  \Ringerl(4.01478 -4.45887) 
  \Ringerl(3. -5.19615)      
  \Ringerl(1.8541 -5.70634)  
  \Ringerl(0.627171 -5.96713) 
  \Ringerl(-0.627171 -5.96713) 
  \Ringerl(-1.8541 -5.70634)  
  \Ringerl(-3. -5.19615)      
  \Ringerl(-4.01478 -4.45887) 
  \Ringerl(-4.8541 -3.52671)  
  \Ringerl(-5.48127 -2.44042) 
  \Ringerl(-5.86889 -1.24747) 
  \Ringerl(-6. 0)             
  \Ringerl(-5.86889 1.24747)  
  \Ringerl(-5.48127 2.44042)  
  \Ringerl(-4.8541 3.52671)   
  \Ringerl(-4.01478 4.45887)  
  \Ringerl(-3. 5.19615)       
  \Ringerl(-1.8541 5.70634)   
  \Ringerl(-0.627171 5.96713) 

}
\caption{Combinatorial realisation of
a $3$-divisible non-crossing partition of type $B_{5}$}
\label{fig:8}
\end{figure}

Given an element $(w_0;w_1,\dots,w_m)\in NC^m(B_n)$,
the bijection, $\Na^m_{B_n}$ say, 
from \cite[Theorem~4.5.6]{ArmDAA} works in the same way as for 
$NC^m(A_{n-1})$. 
That is, recalling from Section~\ref{sec:1} that $W(B_n)$ can be
combinatorially realised as a subgroup of the group of permutations
of $\{1,2,\dots,n,\bar 1,\bar 2,\dots,\bar n\}$, and that, in this
realisation, the standard choice of a Coxeter element is 
$c=[1,2,\dots,n]=(1,2,\dots,n,\bar 1,\bar 2,\dots,\bar n)$,
the announced bijection maps 
$(w_0;w_1,\dots,w_m)\in NC^m(B_n)$ to
$$\Na^m_{B_n}(w_0;w_1,\dots,w_m)=
[1,2,\dots,mn]\,(\bar\ta_{m,1}(w_1))^{-1}\,(\bar\ta_{m,2}(w_2))^{-1}\,\cdots\,
(\bar\ta_{m,m}(w_m))^{-1},$$
where $\bar\ta_{m,i}$ is the obvious extension of the above
transformations $\ta_{m,i}$: namely we let
$$(\bar\ta_{m,i}(w))(mk+i-m)=mw(k)+i-m,\quad 
k=1,2,\dots,n,\bar1,\bar2,\dots,\bar n,$$ 
and 
$(\bar\ta_{m,i}(w))(l)=l$ and
$(\bar\ta_{m,i}(w))(\bar l)=\bar l$ for all $l\not\equiv i$~(mod~$m$),
where $m\bar k+i-m$ is identified with $\overline{mk+i-m}$ for all $k$
and $i$. We refer the reader to \cite[Sec.~4.5]{ArmDAA} for the details.
For example, let $n=5$, $m=3$, $w_0=((2,4))$,
$w_1=[1]=(1,\bar 1)$,
$w_2=((1,4))$, and
$w_3=((2,3))\,((4,5))$. Then $(w_0;w_1,w_2,w_3)$ is mapped to
\begin{align} \notag
\Na^3_{B_{5}}(w_0;w_1,w_2,w_3)&=[1,2,\dots,15]\,[1]\,((2,11))\,
((6,9))\,((12,15))\\
&=((1,\bar2,\overline{12}))\,((3,4,5,6,10,11))\,((7,8,9))\,((13,14,15)).
\label{eq:NCB}
\end{align} 
Figure~\ref{fig:8} shows the graphical representation of \eqref{eq:NCB}.

It is shown in \cite[Theorem~4.5.6]{ArmDAA} that $\Na^m_{B_n}$ is in fact an
isomorphism between the {\it posets} $NC^m(B_n)$ and 
$\widetilde{NC}{}^m(B_n)$. Furthermore,
it is proved in \cite[proof of Theorem~4.3.13]{ArmDAA} that
\begin{equation} \label{eq:BlockB}
c_i(w_0)=b_{i}(\Na^m_{B_n}(w_0;w_1,\dots,w_m)),\quad
i=1,2,\dots,n,
\end{equation}
where $c_i(w_0)$ denotes the number of type $A$
cycles (recall the corresponding terminology from Section~\ref{sec:2}) 
of length $i$ of $w_0$ 
and $b_{i}(\pi)$ denotes one half of
the number of non-zero blocks of size $mi$ in the non-crossing
partition $\pi$. (Recall that non-zero blocks come in
``symmetric" pairs.)
Consequently, under the bijection
$\Na^m_{B_n}$, the element $w_0$ contains a type $B$ cycle of length
$\ell$ if and only if $\Na^m_{B_n}(w_0;w_1,\dots,w_m)$ contains a zero
block of size $m\ell$.

\begin{figure}
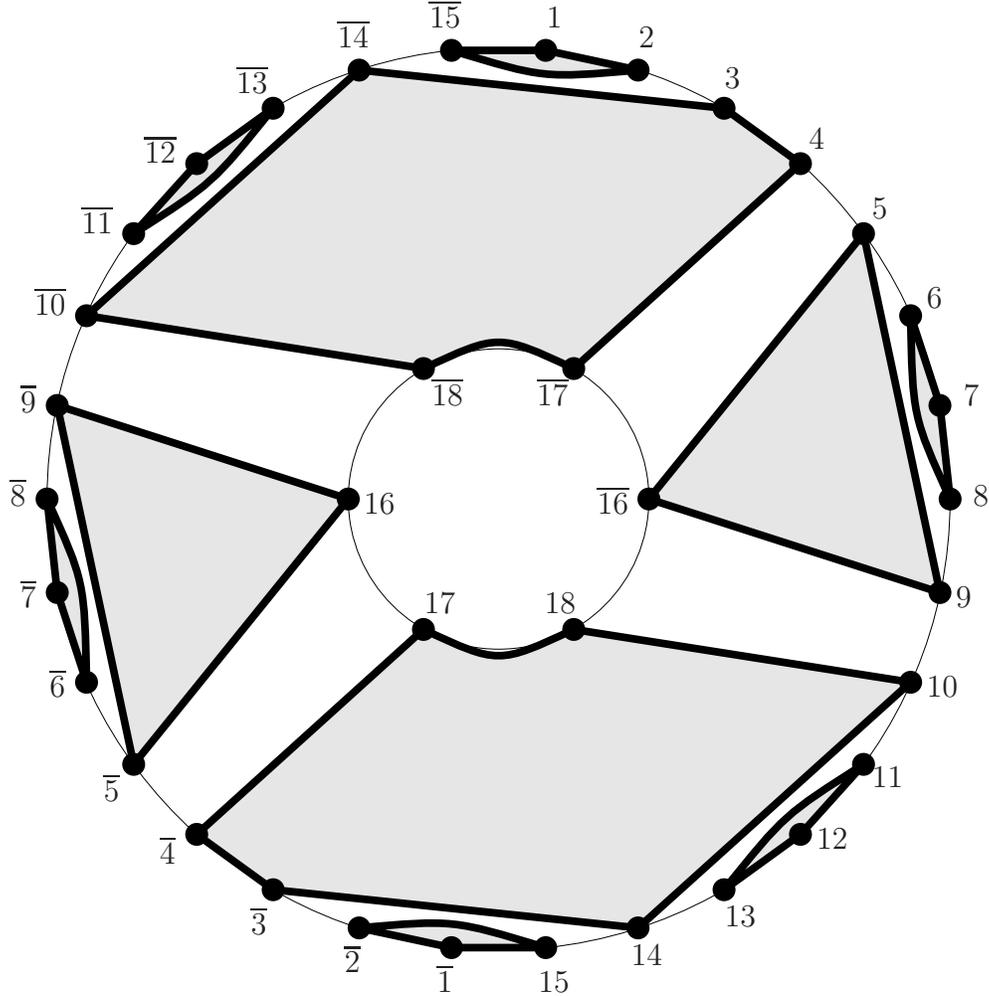

\centertexdraw{
  \drawdim truecm  \linewd.01
  \lcir r:6
  \lcir r:2

  \linewd.10

  \htext(0.627171 6.26713){$1$} 
  \htext(1.8541 6.00634){$2$}   
  \htext(3. 5.45615){$3$}       
  \htext(4.11478 4.65887){$4$} 
  \htext(4.9541 3.72671){$5$}  
  \htext(5.68127 2.54042){$6$} 
  \htext(6.16889 1.24747){$7$} 
  \htext(6.3 -.1){$8$}            
  \htext(6.06889 -1.44747){$9$} 
  \htext(5.68127 -2.64042){$10$} 
  \htext(4.9541 -3.83671){$11$}  
  \htext(4.21478 -4.65887){$12$} 
  \htext(3. -5.69615){$13$}      
  \htext(1.7541 -6.20634){$14$}  
  \htext(0.517171 -6.56713){$15$} 
  \htext(-0.827171 -6.56713){$\overline{1}$} 
  \htext(-2.0541 -6.30634){$\overline{2}$}  
  \htext(-3.3 -5.79615){$\overline{3}$}      
  \htext(-4.51478 -4.85887){$\overline{4}$} 
  \htext(-5.2541 -4.02671){$\overline{5}$}  
  \htext(-5.98127 -2.64042){$\overline{6}$} 
  \htext(-6.36889 -1.44747){$\overline{7}$} 
  \htext(-6.5 -.1){$\overline{8}$}             
  \htext(-6.36889 1.14747){$\overline{9}$}  
  \htext(-6.18127 2.44042){$\overline{10}$}  
  \htext(-5.5541 3.52671){$\overline{11}$}   
  \htext(-4.71478 4.45887){$\overline{12}$}  
  \htext(-3.5 5.39615){$\overline{13}$}       
  \htext(-2.1541 6.00634){$\overline{14}$}   
  \htext(-0.927171 6.26713){$\overline{15}$} 

  \htext(1.3 -.2){$\overline{16}$} 
  \htext(.5 1.23205) {$\overline{17}$} 
  \htext(-.9 1.23205){$\overline{18}$}  
  \htext(-1.8 -.2) {16} 
  \htext(-1. -1.53205){${17}$} 
  \htext(.6 -1.53205){${18}$} 

  \move(0.627171 5.96713)
  \lvec(1.8541 5.70634)
  \clvec(0.6 5.6)(0.6 5.6)(-0.627171 5.96713)
  \lvec(0.627171 5.96713)
  \lfill f:.9

  \move(3. 5.19615)
  \lvec(4.01478 4.45887)
  \lvec(1. 1.73205)
  \clvec(0 2.2)(0 2.2)(-1. 1.73205)
  \lvec(-5.48127 2.44042)
  \lvec(-1.8541 5.70634)
  \lvec(3. 5.19615)
  \lfill f:.9

  \move(5.48127 2.44042)
  \lvec(5.86889 1.24747)
  \lvec(6. 0)     
  \clvec(5.5 1.1)(5.5 1.1)(5.48127 2.44042)
  \lfill f:.9

  \move(4.8541 -3.52671)
  \lvec(4.01478 -4.45887)
  \lvec(3. -5.19615) 
  \clvec(3.8 -4.2)(3.8 -4.2)(4.8541 -3.52671)
  \lfill f:.9

  \move(4.8541 3.52671)
  \lvec(5.86889 -1.24747)
  \lvec(2. 0) 
  \lvec(4.8541 3.52671)
  \lfill f:.9

  \move(-0.627171 -5.96713)
  \lvec(-1.8541 -5.70634)
  \clvec(-0.6 -5.6)(-0.6 -5.6)(0.627171 -5.96713)
  \lvec(-0.627171 -5.96713)
  \lfill f:.9

  \move(-3. -5.19615)
  \lvec(-4.01478 -4.45887)
  \lvec(-1. -1.73205)
  \clvec(-0 -2.2)(-0 -2.2)(1. -1.73205)
  \lvec(5.48127 -2.44042)
  \lvec(1.8541 -5.70634)
  \lvec(-3. -5.19615)
  \lfill f:.9

  \move(-5.48127 -2.44042)
  \lvec(-5.86889 -1.24747)
  \lvec(-6. 0)     
  \clvec(-5.5 -1.1)(-5.5 -1.1)(-5.48127 -2.44042)
  \lfill f:.9

  \move(-4.8541 3.52671)
  \lvec(-4.01478 4.45887)
  \lvec(-3. 5.19615) 
  \clvec(-3.8 4.2)(-3.8 4.2)(-4.8541 3.52671)
  \lfill f:.9

  \move(-4.8541 -3.52671)
  \lvec(-5.86889 1.24747)
  \lvec(-2. 0) 
  \lvec(-4.8541 -3.52671)
  \lfill f:.9

  \Ringerl(0.627171 5.96713) 
  \Ringerl(1.8541 5.70634)   
  \Ringerl(3. 5.19615)       
  \Ringerl(4.01478 4.45887) 
  \Ringerl(4.8541 3.52671)  
  \Ringerl(5.48127 2.44042) 
  \Ringerl(5.86889 1.24747) 
  \Ringerl(6. 0)            
  \Ringerl(5.86889 -1.24747) 
  \Ringerl(5.48127 -2.44042) 
  \Ringerl(4.8541 -3.52671)  
  \Ringerl(4.01478 -4.45887) 
  \Ringerl(3. -5.19615)      
  \Ringerl(1.8541 -5.70634)  
  \Ringerl(0.627171 -5.96713) 
  \Ringerl(-0.627171 -5.96713) 
  \Ringerl(-1.8541 -5.70634)  
  \Ringerl(-3. -5.19615)      
  \Ringerl(-4.01478 -4.45887) 
  \Ringerl(-4.8541 -3.52671)  
  \Ringerl(-5.48127 -2.44042) 
  \Ringerl(-5.86889 -1.24747) 
  \Ringerl(-6. 0)             
  \Ringerl(-5.86889 1.24747)  
  \Ringerl(-5.48127 2.44042)  
  \Ringerl(-4.8541 3.52671)   
  \Ringerl(-4.01478 4.45887)  
  \Ringerl(-3. 5.19615)       
  \Ringerl(-1.8541 5.70634)   
  \Ringerl(-0.627171 5.96713) 

  \Ringerl(2. 0)       
  \Ringerl(1. 1.73205) 
  \Ringerl(-1. 1.73205) 
  \Ringerl(-2. 0)       
  \Ringerl(-1. -1.73205) 
  \Ringerl(1. -1.73205) 
}
\caption{Combinatorial realisation of 
a $3$-divisible non-crossing partition of type $D_{6}$}
\label{fig:9}
\end{figure}

\medskip
The $m$-divisible non-crossing partitions of type $D_n$ cannot be
realised as certain ``partitions non crois\'ees d'un cycle," but as
non-crossing partitions {\it on an annulus} with $2m(n-1)$ vertices
on the outer cycle and $2m$ vertices on the inner cycle, the vertices
on the outer cycle being labelled by $1,2,\dots,mn-m,\bar
1,\bar 2,\dots,\overline{mn-m}$ in {\it clockwise} order, and the 
vertices of the inner cycle being labelled by
$mn-m+1,\dots,mn-1,\break mn,\overline{mn-m+1},\dots,\overline{mn-1},
\overline{mn}$ in {\it counter-clockwise} order. Given a partition
$\pi$ of 
$\{1,2,\dots,mn,\bar 1,\bar 2,\dots,\overline{mn}\}$, we represent it
on this annulus in a manner analogous to Kreweras' graphical
representation of his partitions; namely, we represent each block of
$\pi$ by connecting the vertices labelled by the elements of
the block by curves in clockwise order, the important additional
requirement being here that the curves must be drawn in the {\it interior}
of the annulus. If it is possible to draw the curves in such a way
that no two curves intersect, then the partition is called a {\it
non-crossing partition on the $(2m(n-1),2m)$-annulus}.
Figure~\ref{fig:9} shows a non-crossing partition on the
$(15,6)$-annulus.

\begin{figure}
\centertexdraw{
  \drawdim truecm  

  \linewd.10

  \htext(0.627171 6.26713){$1$} 
  \htext(1.8541 6.00634){$2$}   
  \htext(3. 5.45615){$3$}       
  \htext(4.11478 4.65887){$4$} 
  \htext(4.9541 3.72671){$5$}  
  \htext(5.68127 2.54042){$6$} 
  \htext(6.16889 1.24747){$7$} 
  \htext(6.3 -.1){$8$}            
  \htext(6.06889 -1.44747){$9$} 
  \htext(5.68127 -2.64042){$10$} 
  \htext(4.9541 -3.83671){$11$}  
  \htext(4.21478 -4.65887){$12$} 
  \htext(3. -5.69615){$13$}      
  \htext(1.7541 -6.20634){$14$}  
  \htext(0.517171 -6.56713){$15$} 
  \htext(-0.827171 -6.56713){$\overline{1}$} 
  \htext(-2.0541 -6.30634){$\overline{2}$}  
  \htext(-3.3 -5.79615){$\overline{3}$}      
  \htext(-4.51478 -4.85887){$\overline{4}$} 
  \htext(-5.2541 -4.02671){$\overline{5}$}  
  \htext(-5.98127 -2.64042){$\overline{6}$} 
  \htext(-6.36889 -1.44747){$\overline{7}$} 
  \htext(-6.5 -.1){$\overline{8}$}             
  \htext(-6.36889 1.14747){$\overline{9}$}  
  \htext(-6.18127 2.44042){$\overline{10}$}  
  \htext(-5.5541 3.52671){$\overline{11}$}   
  \htext(-4.71478 4.45887){$\overline{12}$}  
  \htext(-3.5 5.39615){$\overline{13}$}       
  \htext(-2.1541 6.00634){$\overline{14}$}   
  \htext(-0.927171 6.26713){$\overline{15}$} 

  \htext(1.3 -.2){$\overline{16}$} 
  \htext(.5 1.23205) {$\overline{17}$} 
  \htext(-.9 1.23205){$\overline{18}$}  
  \htext(-1.8 -.2) {16} 
  \htext(-1. -1.53205){${17}$} 
  \htext(.6 -1.53205){${18}$} 

  \move(0.627171 5.96713)
  \lvec(1.8541 5.70634)
  \lvec(4.01478 -4.45887)
  \lvec(-0.627171 -5.96713)
  \lvec(-1.8541 -5.70634)
  \lvec(-4.01478 4.45887)
  \lvec(0.627171 5.96713)
  \lfill f:.9

  \move(3. 5.19615) 
  \lvec(4.01478 4.45887)
  \lvec(4.8541 3.52671)
  \lvec(5.48127 2.44042)
  \lvec(5.48127 -2.44042)
  \lvec(4.8541 -3.52671)
  \lvec(3. 5.19615) 
  \lfill f:.9

  \move(5.86889 1.24747)
  \lvec(6. 0)  
  \lvec(5.86889 -1.24747)
  \clvec(5.6 0)(5.6 0)(5.86889 1.24747)
  \lfill f:.9

  \move(3. -5.19615)
  \lvec(1.8541 -5.70634)
  \lvec(0.627171 -5.96713)
  \clvec(1.8 -5.3)(1.8 -5.3)(3. -5.19615)
  \lfill f:.9

  \move(-0.627171 -5.96713)
  \lvec(-1.8541 -5.70634)
  \lvec(-4.01478 4.45887)
  \lvec(0.627171 5.96713)
  \lvec(1.8541 5.70634)
  \lvec(4.01478 -4.45887)
  \lvec(-0.627171 -5.96713)
  \lfill f:.9

  \move(-3. -5.19615) 
  \lvec(-4.01478 -4.45887)
  \lvec(-4.8541 -3.52671)
  \lvec(-5.48127 -2.44042)
  \lvec(-5.48127 2.44042)
  \lvec(-4.8541 3.52671)
  \lvec(-3. -5.19615) 
  \lfill f:.9

  \move(-5.86889 -1.24747)
  \lvec(-6. 0)  
  \lvec(-5.86889 1.24747)
  \clvec(-5.6 0)(-5.6 0)(-5.86889 -1.24747)
  \lfill f:.9

  \move(-3. 5.19615)
  \lvec(-1.8541 5.70634)
  \lvec(-0.627171 5.96713)
  \clvec(-1.8 5.3)(-1.8 5.3)(-3. 5.19615)
  \lfill f:.9

  \move(0 0)
  \linewd.01
  \lcir r:6
  \lcir r:2

  \Ringerl(0.627171 5.96713) 
  \Ringerl(1.8541 5.70634)   
  \Ringerl(3. 5.19615)       
  \Ringerl(4.01478 4.45887) 
  \Ringerl(4.8541 3.52671)  
  \Ringerl(5.48127 2.44042) 
  \Ringerl(5.86889 1.24747) 
  \Ringerl(6. 0)            
  \Ringerl(5.86889 -1.24747) 
  \Ringerl(5.48127 -2.44042) 
  \Ringerl(4.8541 -3.52671)  
  \Ringerl(4.01478 -4.45887) 
  \Ringerl(3. -5.19615)      
  \Ringerl(1.8541 -5.70634)  
  \Ringerl(0.627171 -5.96713) 
  \Ringerl(-0.627171 -5.96713) 
  \Ringerl(-1.8541 -5.70634)  
  \Ringerl(-3. -5.19615)      
  \Ringerl(-4.01478 -4.45887) 
  \Ringerl(-4.8541 -3.52671)  
  \Ringerl(-5.48127 -2.44042) 
  \Ringerl(-5.86889 -1.24747) 
  \Ringerl(-6. 0)             
  \Ringerl(-5.86889 1.24747)  
  \Ringerl(-5.48127 2.44042)  
  \Ringerl(-4.8541 3.52671)   
  \Ringerl(-4.01478 4.45887)  
  \Ringerl(-3. 5.19615)       
  \Ringerl(-1.8541 5.70634)   
  \Ringerl(-0.627171 5.96713) 

  \Ringerl(2. 0)       
  \Ringerl(1. 1.73205) 
  \Ringerl(-1. 1.73205) 
  \Ringerl(-2. 0)       
  \Ringerl(-1. -1.73205) 
  \Ringerl(1. -1.73205) 
}
\caption{A $3$-divisible non-crossing partition of type $D_{6}$ with
zero block}
\label{fig:9a}
\end{figure}

With this definition,
the $m$-divisible non-crossing partitions of type $D_n$ are in
bijection with non-crossing partitions $\pi$ on the $(2m(n-1),2m)$-annulus,
in which successive elements of a block (successive in the circular
order in the graphical representation of the block) are in successive
congruence classes modulo $m$, 
which have the property that if $B$ is a block of $\pi$ then
also $\overline B:=\{\bar x:x\in B\}$ is a block of $\pi$, and which
satisfy an additional restriction concerning their {\it zero block}.
Here again, a zero block is a block $B$ with $\overline B=B$.
The announced additional restriction says that
a zero block can only occur if it contains {\it all\/} the vertices of
the inner cycle, that is, 
$mn-m+1,\dots,mn-1,mn,\overline{mn-m+1},\dots,\overline{mn-1},
\overline{mn}$,
and at least two further elements from the outer cycle.
We denote this set of non-crossing partitions on the
$(2m(n-1),2m)$-annulus by $\widetilde {NC}{}^m(D_{n})$.
A non-crossing partition in $\widetilde {NC}{}^m(D_{n})$ can only have
{\it at most\/} one zero block.
Figures~\ref{fig:9} and \ref{fig:9a} give examples of $3$-divisible
non-crossing partitions of type $D_{6}$, Figure~\ref{fig:9} 
one without a zero block, while Figure~\ref{fig:9a} one with a
zero block. Again, it is clear that 
the condition that $B$ is a block of the partition if and only if
$\overline B$ is a block translates into the condition that the
geometric realisation of the partition is invariant under rotation by
$180^\circ$.

In order to clearly sort out the differences to the earlier combinatorial
realisations of $m$-divisible non-crossing partitions of types
$A_{n-1}$ and $B_n$, we stress that for type $D_n$ there are {\it
three} major features which are not present for the former types: 
(1) here we consider non-crossing partitions on 
an {\it annulus}; (2) it is not sufficient to impose the condition
that the size of every block is divisible by $m$: the condition on
successive elements of a block is stronger; 
(3) there is the above additional restriction
on the zero block (which is not present in type $B_n$).

Given an element $(w_0;w_1,\dots,w_m)\in NC^m(D_n)$,
the bijection, $\Na^m_{D_n}$ say, 
from \cite{KratCG} works as follows.
Recalling from Section~\ref{sec:1} that $W(D_n)$ can be
combinatorially realised as a subgroup of the group of permutations
of $\{1,2,\dots,n,\bar 1,\bar 2,\dots,\bar n\}$, and that, in this
realisation, the standard choice of a Coxeter element is 
$c=[1,2,\dots,n-1]\,[n]=(1,2,\dots,n-1,\bar 1,\bar
2,\dots,\overline{n-1})\,(n,\bar n)$,
the announced bijection maps 
$(w_0;w_1,\dots,w_m)\in NC^m(D_n)$ to
\begin{multline*}
\Na^m_{D_n}(w_0;w_1,\dots,w_m)=
[1,2,\dots,m(n-1)]\,[mn-m+1,\dots,mn-1,mn]\\
\circ(\bar\ta_{m,1}(w_1))^{-1}\,(\bar\ta_{m,2}(w_2))^{-1}\,\cdots\,
(\bar\ta_{m,m}(w_m))^{-1},
\end{multline*}
where $\bar\ta_{m,i}$ is defined as above.
We refer the reader to \cite{KratCG} for the details.
For example, let $n=6$, $m=3$, $w_0=((2,\bar 4))$,
$w_1=((2,\bar 6))\,((4,5))$,
$w_2=((1,\bar 5))\,((2,3))$, and
$w_3=((3,6))$. Then $(w_0;w_1,w_2,w_3)$ is mapped to
\begin{align} \notag
\Na^3_{D_{6}}(w_0;&w_1,w_2,w_3)\\
\notag
&=[1,2,\dots,15]\,[16,17,18]
((4,\overline{16}))\,((10,13))\,((2,\overline{14}))\,
((5,8))\,((9,18))\\
&=((1,2,\overline{15}))\,((3,4,\overline{17},\overline{18},\overline{10},
\overline{14}))\,((5,9,\overline{16}))\,((6,7,8))\,((11,12,13)).
\label{eq:NCD}
\end{align} 
Figure~\ref{fig:9} shows the graphical representation of \eqref{eq:NCD}.

It is shown in \cite{KratCG} that $\Na^m_{D_n}$ is in fact an
isomorphism between the {\it posets} $NC^m(D_n)$ and 
$\widetilde{NC}{}^m(D_n)$. Furthermore,
it is proved in \cite{KratCG} that
\begin{equation} \label{eq:BlockD}
c_i(w_0)=b_{i}(\Na^m_{D_n}(w_0;w_1,\dots,w_m)),\quad
i=1,2,\dots,n,
\end{equation}
where $c_i(w_0)$ denotes the number of type $A$
cycles of length $i$ of $w_0$ 
and $b_{i}(\pi)$ denotes one half of the number of non-zero blocks of size $mi$
in the non-crossing partition $\pi$. (Recall that non-zero blocks come in
``symmetric" pairs.) Consequently, under the bijection
$\Na^m_{D_n}$, the element $w_0$ contains a type $D$ cycle of length
$\ell$ if and only if $\Na^m_{D_n}(w_0;w_1,\dots,w_m)$ contains a zero
block of size $m\ell$.

\section{Decomposition numbers with free factors, and 
enumeration in the poset of generalised non-crossing partitions}
\label{sec:7} 

This section is devoted to applying our formulae from
Sections~\ref{sec:2}--\ref{sec:4} for the
decomposition numbers of the types $A_n$, $B_n$, and $D_n$ to the
enumerative theory of generalised non-crossing partitions for these
types. Theorems~\ref{thm:4}--\ref{thm:6} present formulae for the
number of minimal factorisations of Coxeter elements in types $A_n$,
$B_n$, and $D_n$, respectively, where we do not prescribe the types
of {\it all\/} the factors as for the decomposition numbers, 
but just for some of them, while we
impose rank sum conditions on other factors. Immediate corollaries
are formulae for the number of multi-chains 
$\pi_1\le \pi_2\le \dots\le \pi_{l-1}$, $l$ being given, 
in the posets $\widetilde{NC}{}^m(A_{n-1})$,
$\widetilde{NC}{}^m(B_n)$, and $\widetilde{NC}{}^m(D_n)$,
where the poset rank of $\pi_i$ equals
$r_i$, and where the block structure of $\pi_1$ is prescribed, see
Corollaries~\ref{cor:2}, \ref{cor:3}, and \ref{cor:4}. 
These results in turn imply 
all known enumerative results on
ordinary and generalised non-crossing partitions via appropriate
summations, see the remarks accompanying the corollaries. They also
imply two further new results on chain enumeration in
$\widetilde{NC}{}^m(D_n)$, see Corollaries~\ref{cor:5} and
\ref{cor:6}. We want to stress that, since $\widetilde{NC}{}^m(\Phi)$
and $NC^m(\Phi)$ are isomorphic as posets for $\Phi=A_{n-1},B_n,D_n$,
Corollaries~\ref{cor:2}, \ref{cor:3}, \ref{cor:4}, \ref{cor:5A}, \ref{cor:5}
imply obvious results for $NC^m(\Phi)$ in place of
$\widetilde{NC}{}^m(\Phi)$, $\Phi=A_{n-1},B_n,D_n$, via \eqref{eq:BlockA},
\eqref{eq:BlockB}, respectively \eqref{eq:BlockD}.

\medskip
We begin with our results for type $A_n$. The next theorem
generalises Theorem~\ref{thm:20}, which can be obtained from the
former as the special case in which $l=1$ and $m_1=1$.

\begin{theorem} \label{thm:4}
For a positive integer $d,$ let the types $T_1,T_2,\dots,T_d$ be given,
where 
$$T_i=A_1^{m_1^{(i)}}*A_2^{m_2^{(i)}}*\dots*A_n^{m_n^{(i)}},\quad 
i=1,2,\dots,d,$$ 
and let $l,m_1,m_2,\dots,m_l,s_1,s_2,\dots,s_l$ be given non-negative integers
with
$$\rk T_1+\rk T_2+\dots+\rk T_d+s_1+s_2+\dots+s_l=n.$$
Then the number of factorisations
\begin{equation} \label{eq:FacA}
c=\si_1\si_2\cdots \si_d 
\si_1^{(1)}\si_2^{(1)}\cdots \si_{m_1}^{(1)}
\si_1^{(2)}\si_2^{(2)}\cdots \si_{m_2}^{(2)}
\cdots
\si_1^{(l)}\si_2^{(l)}\cdots \si_{m_l}^{(l)},
\end{equation}
where $c$ is a Coxeter element in $W(A_n),$
such that the type of $\si_i$ is $T_i,$ $i=1,2,\dots,d,$ and such
that
\begin{equation} \label{eq:rkA}
\ell_T( \si_1^{(i)})+\ell_T(\si_2^{(i)})+\cdots
+\ell_T(\si_{m_i}^{(i)})=s_i,\quad i=1,2,\dots,l,
\end{equation}
is given by
\begin{multline} \label{eq:30}
(n+1)^{d-1}
\Bigg(\prod _{i=1} ^{d}
\frac {1} {n-\rk T_i+1}\binom {n-\rk T_i+1}{m_1^{(i)},m_2^{(i)},\dots,
m_n^{(i)}}\Bigg)\\
\times
\binom {m_1(n+1)} {s_1}
\binom {m_2(n+1)} {s_2}\cdots
\binom {m_l(n+1)} {s_l},
\end{multline}
where the multinomial coefficient is defined as in
Lemma~{\em\ref{lem:binsum}}.
\end{theorem}

\begin{proof}
In the factorisation \eqref{eq:FacA},
we first fix also the types of the $\si_i^{(j)}$'s. For
$i=1,2,\dots,m_j$ and $j=1,2,\dots,l$, let the type of $\si_i^{(j)}$ be
$$T_i^{(j)}=A_1^{m_1^{(i,j)}}*A_2^{m_2^{(i,j)}}*\dots*A_n^{m_n^{(i,j)}}.$$
We know that the number of these factorisations is given by
\eqref{eq:1} with $d$ replaced by $d+m_1+m_2+\dots+m_l$ and the
appropriate interpretations of the $m_i^{(j)}$'s. Next we fix
non-negative integers $r_i^{(j)}$ and sum the
expression \eqref{eq:1} over all possible types $T_i^{(j)}$ of
rank $r_i^{(j)}$, $i=1,2,\dots,m_j$, $j=1,2,\dots,l$.
The corresponding summations are completely analogous to the
summation in the proof of Theorem~\ref{thm:20}. As a result, we obtain
\begin{multline*} 
(n+1)^{d-1}
\Bigg(\prod _{i=1} ^{d}
\frac {1} {n-\rk T_i+1}\binom {n-\rk T_i+1}{m_1^{(i)},m_2^{(i)},\dots,
m_n^{(i)}}\Bigg)\\
\times
\binom {n+1} {r_1^{(1)}}
\binom {n+1} {r_2^{(1)}}\cdots
\binom {n+1} {r_{m_1}^{(1)}}
\times
\binom {n+1} {r_1^{(2)}}
\binom {n+1} {r_2^{(2)}}\cdots
\binom {n+1} {r_{m_2}^{(2)}}\\
\times\cdots\times
\binom {n+1} {r_1^{(l)}}
\binom {n+1} {r_2^{(l)}}\cdots
\binom {n+1} {r_{m_l}^{(l)}}
\end{multline*}
for the number of factorisations under consideration.
In view of \eqref{eq:rkA} and \eqref{eq:ellrk},
to obtain the final result, we must sum these expressions
over all non-negative integers $r_1^{(1)},\dots,r_{m_l}^{(l)}$ 
satisfying the equations
\begin{equation} \label{eq:rs}
r_1^{(j)}+r_2^{(j)}+\dots+r_{m_j}^{(j)}=s_j,\quad j=1,2,\dots,l.
\end{equation}
This is easily done by means of the multivariate version of the Chu--Vandermonde
summation. The formula in \eqref{eq:30} follows.
\end{proof}

In view of the combinatorial realisation of 
$m$-divisible non-crossing partitions of type $A_{n-1}$ which we
described in Section~\ref{sec:6},
the special case $d=1$ of the above theorem has the following
enumerative consequence.

\begin{corollary} \label{cor:2}
Let $l$ be a positive integer, and let $s_1,s_2,\dots,s_l$ be
non-negative integers with $s_1+s_2+\dots+s_l=n-1$. 
The number of multi-chains 
$\pi_1\le \pi_2\le \dots\le \pi_{l-1}$ in the poset 
$\widetilde {NC}{}^m(A_{n-1}),$ 
with the property that $\rk(\pi_i)=s_1+s_2+\cdots +s_i,$ 
$i=1,2,\dots,l-1,$ and that the
number of blocks of size $mi$ of $\pi_1$ is $b_i,$ $i=1,2,\dots,n,$
is given by
\begin{equation} \label{eq:31}
\frac {1} {b_1+b_2+\dots+b_n}\binom {b_1+b_2+\dots+b_n}{b_1,b_2,\dots,b_n}
\binom {mn} {s_2}\cdots
\binom {mn} {s_l},
\end{equation}
provided that $b_1+2b_2+\dots+nb_n\le n,$ and is $0$ otherwise.
\end{corollary}

\begin{remark}
The conditions in the statement of the corollary imply that
\begin{equation} \label{eq:sbA}
s_1+b_1+b_2+\dots+b_n=n.
\end{equation}
\end{remark}

\begin{proof} 
Let 
\begin{equation} \label{eq:KetteA}
\pi_1\le \pi_2\le \dots\le \pi_{l-1}
\end{equation} 
be a multi-chain in 
$\widetilde {NC}{}^m(A_{n-1})$. Suppose that, under the bijection
$\Na^m_{A_{n-1}}$, the element $\pi_j$ corresponds to the tuple
$(w_0^{(j)};w_1^{(j)},\dots,w_m^{(j)})$, $j=1,2,\dots,l-1$.
The inequalities in \eqref{eq:KetteA} imply that
$w_1^{(1)},w_2^{(1)},\dots,w_m^{(1)}$ can be factored in the form
$$w_i^{(1)}=u_i^{(2)}u_i^{(3)}\cdots u_i^{(l)},\quad i=1,2,\dots,m,$$
where $u_i^{(l)}=w_i^{(l-1)}$ and, more generally,
\begin{equation} \label{eq:32}
w_i^{(j)}=u_i^{(j+1)}u_i^{(j+2)}\cdots u_i^{(l)},\quad 
i=1,2,\dots,m,\ j=1,2,\dots,l-1.
\end{equation}
For later use, we record that
\begin{align} \notag
c&=w_0^{(j)}w_1^{(j)}\cdots w_m^{(j)}\\
&=w_0^{(j)}
\big(u_1^{(j+1)}u_1^{(j+2)}\cdots u_1^{(l)}\big)
\big(u_2^{(j+1)}u_2^{(j+2)}\cdots u_2^{(l)}\big)\cdots
\big(u_m^{(j+1)}u_m^{(j+2)}\cdots u_m^{(l)}\big).
\label{eq:facwu}
\end{align}

Now, by \eqref{eq:BlockA}, the block structure conditions on $\pi_1$
in the statement of the corollary translate into the condition that the
type of $w_0^{(1)}$ is
\begin{equation} \label{eq:w0A}
A_1^{b_2}*A_2^{b_3}*\dots*A_{n-1}^{b_n}.
\end{equation}
On the other hand, 
using \eqref{eq:rkm}, we see that the rank conditions in the
statement of the corollary mean that
\begin{equation*} 
\ell_T(w_0^{(j)})=s_1+s_2+\dots+s_j,\quad j=1,2,\dots,l-1.
\end{equation*}
In combination with \eqref{eq:facwu}, this yields the conditions
\begin{equation} \label{eq:33}
\ell_T(u_1^{(j)})+\ell_T(u_2^{(j)})+\dots+\ell_T(u_m^{(j)})=s_j,\quad
j=2,3,\dots,l.
\end{equation}
Thus, we want to count the number of factorisations
\begin{equation} \label{eq:bigfact}
c=w_0^{(1)}
\big(u_1^{(2)}u_1^{(3)}\cdots u_1^{(l)}\big)
\big(u_2^{(2)}u_2^{(3)}\cdots u_2^{(l)}\big)\cdots
\big(u_m^{(2)}u_m^{(3)}\cdots u_m^{(l)}\big),
\end{equation}
where the type of $w_0^{(1)}$ is given in \eqref{eq:w0A}, and where
the ``rank conditions" \eqref{eq:33} are satisfied. So, in view of
\eqref{eq:ellrk}, we are in the situation of Theorem~\ref{thm:4} with
$n$ replaced by $n-1$,
$d=1$, $l$ replaced by $l-1$, $s_i$ replaced by $s_{i+1}$, $i=1,2,\dots,l-1$,
$T_1$ the type in \eqref{eq:w0A},
$m_1=m_2=\dots=m_{l-1}=m$, except that the factors are not exactly in the
order as in \eqref{eq:FacA}. However, by \eqref{Aa}
we know that the order of
factors is without relevance. Therefore we
just have to apply Theorem~\ref{thm:4} with the above
specialisations. If we also take into account \eqref{eq:sbA}, then we
arrive immediately at \eqref{eq:31}.
\end{proof}

This result is new even for $m=1$, that is, for the poset of
Kreweras' non-crossing partitions of $\{1,2,\dots,n\}$. It implies
all known results on Kreweras' non-crossing partitions and the
$m$-divisible non-crossing partitions of Edelman. Namely, for $l=2$
it reduces to Armstrong's result \cite[Theorem~4.4.4 with $\ell=1$]{ArmDAA}
on the number of $m$-divisible non-crossing partitions 
in $\widetilde{NC}{}^m(A_{n-1})$ with a given
block structure, which itself contains Kreweras' result
\cite[Theorem~4]{KrewAC}
on his non-crossing partitions with a given block structure as a
special case. If we sum the expression \eqref{eq:31} over all
$s_2,s_3,\dots,s_l$ with $s_2+s_3+\dots+s_l=n-1-s_1$, then we obtain
that the number of all multi-chains 
$\pi_1\le \pi_2\le \dots\le \pi_{l-1}$ in Edelman's poset 
$\widetilde{NC}{}^m(A_{n-1})$ of
$m$-divisible non-crossing partitions of $\{1,2,\dots,mn\}$ in which
$\pi_1$ has $b_i$ blocks of size $mi$ equals
\begin{multline} \label{eq:34}
\frac {1} {b_1+b_2+\dots+b_n}\binom {b_1+b_2+\dots+b_n}{b_1,b_2,\dots,b_n}
\binom {(l-1)mn} {n-s_1-1}\\=
\frac {1} {b_1+b_2+\dots+b_n}\binom {b_1+b_2+\dots+b_n}{b_1,b_2,\dots,b_n}
\binom {(l-1)mn} {b_1+b_2+\dots+b_n-1},
\end{multline}
provided that $b_1+2b_2+\dots+nb_n\le n$,
a result originally due to Armstrong \cite[Theorem~4.4.4]{ArmDAA}.
On the other hand, if we sum the expression \eqref{eq:31} over all
possible $b_1,b_2,\dots,b_n$, that is, $b_2+2b_3+\dots+(n-1)b_n=s_1$, 
use of Lemma~\ref{lem:binsum} with $M=n-s_1$ and $r=s_1$ yields that
the number of all multi-chains 
$\pi_1\le \pi_2\le \dots\le \pi_{l-1}$ in Edelman's poset 
$\widetilde{NC}{}^m(A_{n-1})\cong NC^m(A_{n-1})$
where $\pi_i$ is of rank $s_1+s_2+\dots+s_i$, $i=1,2,\dots,l-1$, equals
\begin{equation} \label{eq:35}
\frac {1} {n}\binom {n} {s_1}
\binom {mn} {s_2}\cdots
\binom {mn} {s_l},
\end{equation}
a result originally due to Edelman \cite[Theorem~4.2]{EdelAA}.
Clearly, this formula contains at the same time a formula for the
number of all $m$-divisible non-crossing partitions of
$\{1,2,\dots,mn\}$ with a given number of blocks upon setting $l=2$
(cf.\ \cite[Lemma~4.1]{EdelAA}),
as well as that it implies that the total number of multi-chains 
$\pi_1\le \pi_2\le \dots\le \pi_{l-1}$ in
the poset of these partitions is 
\begin{equation} \label{eq:36}
\frac {1} {n}\binom {(l-1)mn+n}{n-1}
\end{equation}
upon summing \eqref{eq:35} over all
non-negative integers $s_1,s_2,\dots,s_l$ with $s_1+s_2+\dots+s_l=n-1$ by means
of the multivariate Chu--Vandermonde summation, thus recovering the
formula \cite[Cor.~4.4]{EdelAA} for the zeta polynomial of the poset of
$m$-divisible non-crossing partitions of type $A_{n-1}$.
As special case $l=2$, we recover the well-known fact that the
total number of $m$-divisible non-crossing partitions of
$\{1,2,\dots,mn\}$ is $\frac {1} {n}\binom {(m+1)n}{n-1}$.

\medskip
We continue with our results for type $B_n$. We formulate the theorem
below on factorisations in $W(B_n)$ only with restrictions on the
{\it combinatorial\/} type of some factors. An analogous result with
{\it group-theoretical\/} type instead could be easily derived as
well. We omit this here because, for the combinatorial applications
that we have in mind, combinatorial type suffices.
We remark that the theorem
generalises Theorem~\ref{thm:21}, which can be obtained from the
former as the special case in which $l=1$ and $m_1=1$.

\begin{theorem} \label{thm:5}
{\em(i)} For a positive integer $d$,
let the types $T_1,T_2,\dots,T_d$ be given,
where 
$$T_i=A_1^{m_1^{(i)}}*A_2^{m_2^{(i)}}*\dots*A_n^{m_n^{(i)}},\quad 
i=1,2,\dots,j-1,j+1,\dots,d,$$ 
and
$$T_j=B_\al*A_1^{m_1^{(j)}}*A_2^{m_2^{(j)}}*\dots*A_n^{m_n^{(j)}},
$$ 
for some $\al\ge1,$
and let $l,m_1,m_2,\dots,m_l,s_1,s_2,\dots,s_l$ be given non-negative integers
with                   
\begin{equation} \label{eq:rksum}
\rk T_1+\rk T_2+\dots+\rk T_d+s_1+s_2+\dots+s_l=n.
\end{equation}
Then the number of factorisations
\begin{equation} \label{eq:FacB}
c=\si_1\si_2\cdots \si_d 
\si_1^{(1)}\si_2^{(1)}\cdots \si_{m_1}^{(1)}
\si_1^{(2)}\si_2^{(2)}\cdots \si_{m_2}^{(2)}
\cdots
\si_1^{(l)}\si_2^{(l)}\cdots \si_{m_l}^{(l)},
\end{equation}
where $c$ is a Coxeter element in $W(B_n),$
such that the {\em combinatorial}
type of $\si_i$ is $T_i,$ $i=1,2,\dots,d,$ and such
that
\begin{equation} \label{eq:rkB}
\ell_T( \si_1^{(i)})+\ell_T(\si_2^{(i)})+\cdots
+\ell_T(\si_{m_i}^{(i)})=s_i,\quad i=1,2,\dots,l,
\end{equation}
is given by
\begin{multline} \label{eq:37a}
n^{d-1}
\binom {n-\rk T_j}{m_1^{(j)},m_2^{(j)},\dots,
m_n^{(j)}}
\Bigg(
\underset{i\ne j}{\prod _{i=1} ^{d}}
\frac {1} {n-\rk T_i}\binom {n-\rk T_i}{m_1^{(i)},m_2^{(i)},\dots,
m_n^{(i)}}\Bigg)\\
\times
\binom {m_1n} {s_1}
\binom {m_2n} {s_2}\cdots
\binom {m_ln} {s_l},
\end{multline}
where the multinomial coefficient is defined as in
Lemma~{\em\ref{lem:binsum}}.

\smallskip
{\em(ii)} For a positive integer $d,$
let the types $T_1,T_2,\dots,T_d$ be given,
where 
$$T_i=A_1^{m_1^{(i)}}*A_2^{m_2^{(i)}}*\dots*A_n^{m_n^{(i)}},\quad 
i=1,2,\dots,d,$$ 
and let $l,m_1,m_2,\dots,m_l,s_1,s_2,\dots,s_l$ be given non-negative integers.
Then the number of factorisations \eqref{eq:FacB} which satisfy
\eqref{eq:rkB} plus the condition
that the {\em combinatorial} type of $\si_i$ is $T_i,$ $i=1,2,\dots,d,$ 
is given by
\begin{multline} \label{eq:37b}
n^{d-1}
(n-\rk T_1-\rk T_2-\dots-\rk T_d)
\Bigg(
{\prod _{i=1} ^{d}}
\frac {1} {n-\rk T_i}\binom {n-\rk T_i}{m_1^{(i)},m_2^{(i)},\dots,
m_n^{(i)}}\Bigg)\\
\times
\binom {m_1n} {s_1}
\binom {m_2n} {s_2}\cdots
\binom {m_ln} {s_l}.
\end{multline}
\end{theorem}

\begin{proof}
We start with the proof of item~(i).
In the factorisation \eqref{eq:FacB},
we first fix also the types of the $\si_i^{(j)}$'s. For
$i=1,2,\dots,m_j$ and $j=1,2,\dots,l$ let the type of $\si_i^{(j)}$ be
$$T_i^{(j)}=A_1^{m_1^{(i,j)}}*A_2^{m_2^{(i,j)}}*\dots*A_n^{m_n^{(i,j)}}.$$
We know that the number of these factorisations is given by
\eqref{eq:2} with $d$ replaced by $d+m_1+m_2+\dots+m_l$ and the
appropriate interpretations of the $m_i^{(j)}$'s. Next we fix
non-negative integers $r_i^{(j)}$ and sum the
expression \eqref{eq:2} over all possible types $T_i^{(j)}$ of
rank $r_i^{(j)}$, $i=1,2,\dots,m_j$, $j=1,2,\dots,l$.
The corresponding summations are completely analogous to the first
summation in the proof of Theorem~\ref{thm:21}. As a result, we obtain
\begin{multline*} 
n^{d-1}
\binom {n-\rk T_j}{m_1^{(j)},m_2^{(j)},\dots,
m_n^{(j)}}
\Bigg(
\underset{i\ne j}{\prod _{i=1} ^{d}}
\frac {1} {n-\rk T_i}\binom {n-\rk T_i}{m_1^{(i)},m_2^{(i)},\dots,
m_n^{(i)}}\Bigg)\\
\times
\binom {n} {r_1^{(1)}}
\binom {n} {r_2^{(1)}}\cdots
\binom {n} {r_{m_1}^{(1)}}
\times
\binom {n} {r_1^{(2)}}
\binom {n} {r_2^{(2)}}\cdots
\binom {n} {r_{m_2}^{(2)}}\\
\times\cdots\times
\binom {n} {r_1^{(l)}}
\binom {n} {r_2^{(l)}}\cdots
\binom {n} {r_{m_l}^{(l)}}
\end{multline*}
for the number of factorisations under consideration.
In view of \eqref{eq:rkB} and \eqref{eq:ellrk},
to obtain the final result, we must sum these expressions
over all non-negative integers $r_1^{(1)},\dots,r_{m_l}^{(l)}$ 
satisfying the equations
$$r_1^{(j)}+r_2^{(j)}+\dots+r_{m_j}^{(j)}=s_j,\quad j=1,2,\dots,l.$$
This is easily done by means of the multivariate version of the Chu--Vandermonde
summation. The formula in \eqref{eq:37a} follows.

The proof of item~(ii) is completely analogous, we must, however, cope
with the complication that the type $B$ cycle, which,
according to Theorem~\ref{thm:2}, must occur in the disjoint
cycle decomposition of {\it exactly} one of
the factors on the right-hand side of \eqref{eq:FacB}, can occur in
{\it any} of the $\si_i^{(j)}$'s. So, let us fix the types of the
$\si_i^{(j)}$'s to
$$T_i^{(j)}=A_1^{m_1^{(i,j)}}*A_2^{m_2^{(i,j)}}*\dots*A_n^{m_n^{(i,j)}},$$
$i=1,2,\dots,m_j$, $j=1,2,\dots,l$, except for 
$(i,j)=(p,q)$, where we require that the type of $\si_{p}^{(q)}$ is
$$T_{p}^{(q)}=B_\al*A_1^{\widetilde m_1}*A_2^{\widetilde
m_2}*\dots*A_n^{\widetilde m_n}.$$
Again, we know that the number of these factorisations is given by
\eqref{eq:2} with $d$ replaced by $d+m_1+m_2+\dots+m_l$ and the
appropriate interpretations of the $m_i^{(j)}$'s.
Now we fix
non-negative integers $r_i^{(j)}$ and sum the
expression \eqref{eq:2} over all possible types $T_i^{(j)}$ of
rank $r_i^{(j)}$, $i=1,2,\dots,m_j$, $j=1,2,\dots,l$. Again, 
the corresponding summations are completely analogous to the
summations in the proof of Theorem~\ref{thm:21}. In particular, 
the summation over all possible types $T_{p}^{(q)}$ of
rank $r_{p}^{(q)}$ is essentially the summation on the
right-hand side of \eqref{eq:alsum} with $d$ replaced by
$d+m_1+m_2+\dots+m_l$ and $r$ replaced by 
$r_{p}^{(q)}$. If we use what we know from the proof of
Theorem~\ref{thm:21}, then the result of the summations is found to be
\begin{multline} \label{eq:ExprB}
n^{d-1}
\Bigg(
{\prod _{i=1} ^{d}}
\frac {1} {n-\rk T_i}\binom {n-\rk T_i}{m_1^{(i)},m_2^{(i)},\dots,
m_n^{(i)}}\Bigg)\\
\times
\binom {n} {r_1^{(1)}}
\binom {n} {r_2^{(1)}}\cdots
\binom {n} {r_{m_1}^{(1)}}
\times\cdots\times
\binom {n} {r_1^{(q)}}
\cdots
{r_{p}^{(q)}}\binom {n} {r_{p}^{(q)}}
\cdots
\binom {n} {r_{m_q}^{(q)}}\\
\times\cdots\times
\binom {n} {r_1^{(l)}}
\binom {n} {r_2^{(l)}}\cdots
\binom {n} {r_{m_l}^{(l)}}.
\end{multline}
The reader should note that the term ${r_{p}^{(q)}}\binom {n}
{r_{p}^{(q)}}$ in this expression results from the summation over all
types $T_p^{(q)}$ of rank $r_p^{(q)}$ (compare \eqref{eq:rksum0} with
$r$ replaced by $r_p^{(q)}$; we have $\binom {n-1}{r_p^{(q)}-1}=
\frac {r_p^{(q)}} {n}\binom n{r_p^{(q)}}$).
Using \eqref{eq:rksum}, \eqref{eq:rkB} and \eqref{eq:ellrk}, 
we see that the sum of all
$r_p^{(q)}$ over $p=1,2,\dots,m_q$ and $q=1,2,\dots,l$ must be
$n-\rk T_1-\rk T_2-\dots-\rk T_d$.
Hence, the sum of the expressions \eqref{eq:ExprB} over all $(p,q)$
equals
\begin{multline*} 
n^{d-1}
(n-\rk T_1-\rk T_2-\dots-\rk T_d)
\Bigg(
{\prod _{i=1} ^{d}}
\frac {1} {n-\rk T_i}\binom {n-\rk T_i}{m_1^{(i)},m_2^{(i)},\dots,
m_n^{(i)}}\Bigg)\\
\times
\binom {n} {r_1^{(1)}}
\binom {n} {r_2^{(1)}}\cdots
\binom {n} {r_{m_1}^{(1)}}
\times\cdots\times
\binom {n} {r_1^{(q)}}
\binom {n} {r_2^{(q)}}
\cdots
\binom {n} {r_{m_q}^{(q)}}\\
\times\cdots\times
\binom {n} {r_1^{(l)}}
\binom {n} {r_2^{(l)}}\cdots
\binom {n} {r_{m_l}^{(l)}}.
\end{multline*}
Finally, we must sum these expressions
over all non-negative integers $r_1^{(1)},\dots,r_{m_l}^{(l)}$ 
satisfying the equations
$$r_1^{(j)}+r_2^{(j)}+\dots+r_{m_j}^{(j)}=s_j,\quad j=1,2,\dots,l.$$
Once again, 
this is easily done by means of the multivariate version of the Chu--Vandermonde
summation. As a result, we obtain the formula in \eqref{eq:37b}.
\end{proof}

In view of the combinatorial realisation of 
$m$-divisible non-crossing partitions of type $B_n$ which we
described in Section~\ref{sec:6},
the special case $d=1$ of the above theorem has the following
enumerative consequence.

\begin{corollary} \label{cor:3}
Let $l$ be a positive integer, and let $s_1,s_2,\dots,s_l$ be
non-negative integers with $s_1+s_2+\dots+s_l=n$. 
The number of multi-chains 
$\pi_1\le \pi_2\le \dots\le \pi_{l-1}$ in the poset 
$\widetilde {NC}{}^m(B_n)$ 
with the property that $\rk(\pi_i)=s_1+s_2+\cdots +s_i,$ 
$i=1,2,\dots,l-1,$ and that the
number of non-zero blocks of $\pi_1$ of size $mi$ is $2b_i,$ $i=1,2,\dots,n,$
is given by
\begin{equation} \label{eq:38}
\binom {b_1+b_2+\dots+b_n}{b_1,b_2,\dots,b_n}
\binom {mn} {s_2}\cdots
\binom {mn} {s_l},
\end{equation}
provided that $b_1+2b_2+\dots+nb_n\le n,$ and is $0$ otherwise.
\end{corollary}

\begin{remark}
The conditions in the statement of the corollary imply that
\begin{equation} \label{eq:sbB}
s_1+b_1+b_2+\dots+b_n=n.
\end{equation}
The reader should recall from Section~\ref{sec:6}, that non-zero
blocks of elements $\pi$ of $\widetilde {NC}{}^m(B_n)$ occur in pairs
since, with a block $B$ of $\pi$, also $\overline B$ is a block of
$\pi$.
\end{remark}

\begin{proof} 
The arguments are completely analogous to those of the proof of
Corollary~\ref{cor:2}. The conclusion here is that we need
Theorem~\ref{thm:5} with
$d=1$, $l$ replaced by $l-1$, $s_i$ replaced by $s_{i+1}$,
$i=1,2,\dots,l-1$, $m_1=m_2=\dots=m_{l-1}=m$, and $T_1$ of the type 
$$
B_{n-b_1-2b_2-\dots-nb_n}*A_1^{b_2}*A_2^{b_3}*\dots*A_{n-1}^{b_n}
$$
in the case that $b_1+2b_2+\dots+nb_n<n$ (which enforces the
existence of a zero block of size $2(n-b_1-2b_2-\dots-nb_n)$ in
$\pi_1$), respectively
$$
A_1^{b_2}*A_2^{b_3}*\dots*A_{n-1}^{b_n}
$$
if not. So, depending on the case in which we are, we have to apply
\eqref{eq:37a}, respectively \eqref{eq:37b}. However, for $d=1$ these
two formulae become identical. More precisely, under the above
specialisations, they reduce to
$$
\binom {n-\rk T_1}{b_2,b_3,\dots,b_n}
\binom {mn} {s_2}\cdots
\binom {mn} {s_l}.
$$
If we also take into account \eqref{eq:sbB}, then we
arrive immediately at \eqref{eq:38}.
\end{proof}

This result is new even for $m=1$, that is, for the poset of
Reiner's type $B_n$ non-crossing partitions. It implies
all known results on these non-crossing partitions and their
extension to $m$-divisible type $B_n$ non-crossing partitions due to
Armstrong. Namely, for $l=2$
it reduces to Armstrong's result \cite[Theorem~4.5.11 with $\ell=1$]{ArmDAA}
on the number of elements of $\widetilde{NC}{}^m(B_n)$
with a given block structure, which itself contains Athanasiadis' result
\cite[Theorem~2.3]{AthaAE}
on Reiner's type $B_n$ non-crossing partitions with a given block structure as a
special case. If we sum the expression \eqref{eq:38} over all
$s_2,s_3,\dots,s_l$ with $s_2+s_3+\dots+s_l=n-s_1$, then we obtain
that the number of all multi-chains 
$\pi_1\le \pi_2\le \dots\le \pi_{l-1}$ in 
$\widetilde {NC}{}^m(B_n)$ in which
$\pi_1$ has $2b_i$ non-zero blocks of size $mi$ equals
\begin{equation} \label{eq:41}
\binom {b_1+b_2+\dots+b_n}{b_1,b_2,\dots,b_n}
\binom {(l-1)mn} {n-s_1}=
\binom {b_1+b_2+\dots+b_n}{b_1,b_2,\dots,b_n}
\binom {(l-1)mn} {b_1+b_2+\dots+b_n},
\end{equation}
provided that $b_1+2b_2+\dots+nb_n\le n$, 
a result originally due to Armstrong \cite[Theorem~4.5.11]{ArmDAA}.
On the other hand, if we sum the expression \eqref{eq:38} over all
possible $b_1,b_2,\dots,b_n$, that is, over $b_2+2b_3+\dots+(n-1)b_n\le s_1$, 
use of Lemma~\ref{lem:binsum} with $M=n-s_1$ and $r=s_1-\al$
(where $\al$ stands for the difference between $s_1$ and
$b_2+2b_3+\dots+(n-1)b_n$) yields that
the number of all multi-chains 
$\pi_1\le \pi_2\le \dots\le \pi_{l-1}$ in 
$\widetilde {NC}{}^m(B_n)\cong NC^m(B_n)$ 
where $\pi_i$ is of rank $s_1+s_2+\dots+s_i$, $i=1,2,\dots,l-1$, equals
\begin{equation} \label{eq:42}
\sum _{\al=0} ^{n}\binom {n-\al-1} {s_1-\al}
\binom {mn} {s_2}\cdots
\binom {mn} {s_l}
=\binom {n} {s_1}
\binom {mn} {s_2}\cdots
\binom {mn} {s_l},
\end{equation}
another result due to Armstrong \cite[Theorem~4.5.7]{ArmDAA}.
Clearly, this formula contains at the same time a formula for the
number of all elements of
$\widetilde {NC}{}^m(B_n)\cong NC^m(B_n)$ 
with a given number of blocks (equivalently, a given rank) upon setting $l=2$
(cf.\ \cite[Theorem~4.5.8]{ArmDAA}),
as well as that it implies that the total number of multi-chains 
$\pi_1\le \pi_2\le \dots\le \pi_{l-1}$ in
$\widetilde {NC}{}^m(B_n)\cong NC^m(B_n)$ is
\begin{equation} \label{eq:43}
\binom {(l-1)mn+n}{n}
\end{equation}
upon summing \eqref{eq:42} over all 
non-negative integers $s_1,s_2,\dots,s_l$ with $s_1+s_2+\dots+s_l=n$ by means
of the multivariate Chu--Vandermonde summation, thus recovering the
formula \cite[Theorem~3.6.9]{ArmDAA} for the zeta polynomial of the poset of
generalised non-crossing partitions in the case of type $B_n$.
As special case $l=2$, we recover the fact that the
cardinality of $\widetilde {NC}{}^m(B_n)\cong NC^m(B_n)$ 
is $\binom {(m+1)n}{n}$ (cf.\ \cite[Theorem~3.5.3]{ArmDAA}).

\medskip
The final set of results in this section concerns the type $D_n$.
We start with Theorem~\ref{thm:6}, 
the result on factorisations in $W(D_n)$ which is
analogous to Theorems~\ref{thm:4} and \ref{thm:5}. Similar to
Theorem~\ref{thm:5}, we formulate the theorem
only with restrictions on the
{\it combinatorial\/} type of some factors. An analogous result with
{\it group-theoretical\/} type instead could be easily derived as
well. We refrain from doing this here because, again,
for the combinatorial applications
that we have in mind, combinatorial type suffices.
We remark that the theorem
generalises Theorem~\ref{thm:22}, which can be obtained from the
former as the special case in which $l=1$ and $m_1=1$.

\begin{theorem} \label{thm:6}
{\em(i)} For a positive integer $d,$
let the types $T_1,T_2,\dots,T_d$ be given,
where 
$$T_i=A_1^{m_1^{(i)}}*A_2^{m_2^{(i)}}*\dots*A_n^{m_n^{(i)}},\quad 
i=1,2,\dots,j-1,j+1,\dots,d,$$ 
and
$$T_j=D_\al*A_1^{m_1^{(j)}}*A_2^{m_2^{(j)}}*\dots*A_n^{m_n^{(j)}},
$$ 
for some $\al\ge2,$
and let $l,m_1,m_2,\dots,m_l,s_1,s_2,\dots,s_l$ be given non-negative integers
with                   
\begin{equation*} 
\rk T_1+\rk T_2+\dots+\rk T_d+s_1+s_2+\dots+s_l=n.
\end{equation*}
Then the number of factorisations
\begin{equation} \label{eq:FacD}
c=\si_1\si_2\cdots \si_d 
\si_1^{(1)}\si_2^{(1)}\cdots \si_{m_1}^{(1)}
\si_1^{(2)}\si_2^{(2)}\cdots \si_{m_2}^{(2)}
\cdots
\si_1^{(l)}\si_2^{(l)}\cdots \si_{m_l}^{(l)},
\end{equation}
where $c$ is a Coxeter element in $W(D_n),$
such that the {\em combinatorial}
type of $\si_i$ is $T_i,$ $i=1,2,\dots,d,$ and such
that
\begin{equation} \label{eq:rkD}
\ell_T( \si_1^{(i)})+\ell_T(\si_2^{(i)})+\cdots
+\ell_T(\si_{m_i}^{(i)})=s_i,\quad i=1,2,\dots,l,
\end{equation}
is given by
\begin{multline} \label{eq:44a}
(n-1)^{d-1}
\binom {n-\rk T_j}{m_1^{(j)},m_2^{(j)},\dots,
m_n^{(j)}}
\Bigg(
\underset{i\ne j}{\prod _{i=1} ^{d}}
\frac {1} {n-\rk T_i-1}\binom {n-\rk T_i-1}{m_1^{(i)},m_2^{(i)},\dots,
m_n^{(i)}}\Bigg)\\
\times
\binom {m_1(n-1)} {s_1}
\binom {m_2(n-1)} {s_2}\cdots
\binom {m_l(n-1)} {s_l},
\end{multline}
the multinomial coefficient being defined as in
Lemma~{\em\ref{lem:binsum}}.

\smallskip
{\em(ii)} For a positive integer $d,$
let the types $T_1,T_2,\dots,T_d$ be given,
where 
$$T_i=A_1^{m_1^{(i)}}*A_2^{m_2^{(i)}}*\dots*A_n^{m_n^{(i)}},\quad 
i=1,2,\dots,d,$$ 
and let $l,m_1,m_2,\dots,m_l,s_1,s_2,\dots,s_l$ be given non-negative integers.
Then the number of factorisations \eqref{eq:FacD} which satisfy
\eqref{eq:rkD} as well as the condition
that the {\em combinatorial} type of $\si_i$ is $T_i,$ $i=1,2,\dots,d,$ 
is given by
\begin{multline} \label{eq:44b}
2(n-1)^{d-1}\Bigg(
\sum _{j=1} ^{d}
\binom {n-\rk T_j}{m_1^{(j)},m_2^{(j)},\dots,m_n^{(j)}}
\underset{i\ne j} {\prod _{i=1} ^{d}}
\frac {1} {n-\rk T_i-1}\binom {n-\rk T_i-1}{m_1^{(i)},m_2^{(i)},\dots,
m_n^{(i)}}\Bigg)\\
\times
\binom {m_1(n-1)} {s_1}
\binom {m_2(n-1)} {s_2}\cdots
\binom {m_l(n-1)} {s_l}\\
+(n-1)^{d}
\Bigg(
{\prod _{i=1} ^{d}}
\frac {1} {n-\rk T_i-1}\binom {n-\rk T_i-1}{m_1^{(i)},m_2^{(i)},\dots,
m_n^{(i)}}\Bigg)\kern4cm\\
\times
\sum _{j=1} ^{l}m_j\binom {m_1(n-1)} {s_1}\cdots
\binom {m_j(n-1)-1} {s_j-2}\cdots
\binom {m_l(n-1)} {s_l}\\
-2(d-1)(n-1)^{d}
\Bigg(
{\prod _{i=1} ^{d}}
\frac {1} {n-\rk T_i-1}\binom {n-\rk T_i-1}{m_1^{(i)},m_2^{(i)},\dots,
m_n^{(i)}}\Bigg)\kern3cm\\
\times
\binom {m_1(n-1)} {s_1}
\binom {m_2(n-1)} {s_2}\cdots
\binom {m_l(n-1)} {s_l}.
\end{multline}
\end{theorem}

\begin{proof}
The proof of item~(i) is completely analogous to the proof of item~(i)
in Theorem~\ref{thm:5}. Making reference to that proof,
the only difference is that, instead of the
expression \eqref{eq:2}, we must use \eqref{eq:6} 
with $d$ replaced by $d+m_1+m_2+\dots+m_l$ and the
appropriate interpretations of the $m_i^{(j)}$'s. 
The summations over types $T_i^{(j)}$ with fixed rank $r_i^{(j)}$
are carried out by using \eqref{eq:binsum} with $M=n-r-1$.
Subsequently, the summations over the $r_i^{(j)}$'s satisfying
\eqref{eq:rs} are done by the multivariate version of the Chu--Vandermonde
summation. We leave it to the reader to fill in the details to
finally arrive at \eqref{eq:44a}.

Similarly, the proof of item~(ii) is analogous to the proof of item~(ii) in
Theorem~\ref{thm:5}. However, we must cope
with the complication that there may or may not be a
type $D$ cycle in the disjoint cycle decomposition of one of the
$\si_i^{(j)}$'s on the right-hand side of \eqref{eq:FacD}. 
In the case that there is no type $B$ cycle, 
we fix the types of the $\si_i^{(j)}$'s to
$$T_i^{(j)}=A_1^{m_1^{(i,j)}}*A_2^{m_2^{(i,j)}}*\dots*A_n^{m_n^{(i,j)}},$$
$i=1,2,\dots,m_j$, $j=1,2,\dots,l$, and sum the expression
\eqref{eq:7comb} with $d$ replaced by $d+m_1+m_2+\dots+m_l$ and the
appropriate interpretations of the $m_i^{(j)}$'s over all possible
types $T_i^{(j)}$ with rank $r_i^{(j)}$, $i=1,2,\dots,m_j$, $j=1,2,\dots,l$. 
This yields the expression
\begin{multline} \label{eq:ExprD1}
2(n-1)^{d-1}\Bigg(
\sum _{j=1} ^{d}
\binom {n-\rk T_j}{m_1^{(j)},m_2^{(j)},\dots,m_n^{(j)}}
\underset{i\ne j} {\prod _{i=1} ^{d}}
\frac {1} {n-\rk T_i-1}\binom {n-\rk T_i-1}{m_1^{(i)},m_2^{(i)},\dots,
m_n^{(i)}}\Bigg)\\
\times
\binom {n-1} {r_1^{(1)}}
\binom {n-1} {r_2^{(1)}}\cdots
\binom {n-1} {r_{m_1}^{(1)}}
\times\cdots\times
\binom {n-1} {r_1^{(l)}}
\binom {n-1} {r_2^{(l)}}\cdots
\binom {n-1} {r_{m_l}^{(l)}}\\
+2(n-1)^{d}\Bigg(
\sum _{j=1} ^{l}\sum _{i=1} ^{m_j}
{\prod _{i=1} ^{d}}
\frac {1} {n-\rk T_i-1}\binom {n-\rk T_i-1}{m_1^{(i)},m_2^{(i)},\dots,
m_n^{(i)}}\Bigg)\\
\times
\binom {n-1} {r_1^{(1)}}
\binom {n-1} {r_2^{(1)}}\cdots
\binom {n-1} {r_{m_1}^{(1)}}
\times\cdots\times
\binom {n-1} {r_1^{(l)}}
\binom {n-1} {r_2^{(l)}}\cdots
\binom {n-1} {r_{m_l}^{(l)}}\\
-2(d-1)(n-1)^{d-1}
\Bigg(
{\prod _{i=1} ^{d}}
\frac {1} {n-\rk T_i-1}\binom {n-\rk T_i-1}{m_1^{(i)},m_2^{(i)},\dots,
m_n^{(i)}}\Bigg)\\
\times
\binom {n-1} {r_1^{(1)}}
\binom {n-1} {r_2^{(1)}}\cdots
\binom {n-1} {r_{m_1}^{(1)}}
\times\cdots\times
\binom {n-1} {r_1^{(l)}}
\binom {n-1} {r_2^{(l)}}\cdots
\binom {n-1} {r_{m_l}^{(l)}}.
\end{multline}
In the case that there appears, however, a type $B$ cycle in $\si_p^{(q)}$,
say, we adopt the same set-up as above, except that
we restrict $\si_p^{(q)}$ to types of the form
$$T_{p}^{(q)}=D_\al*A_1^{\widetilde m_1}*A_2^{\widetilde
m_2}*\dots*A_n^{\widetilde m_n}.$$
Subsequently, we sum the expression \eqref{eq:6} 
with $d$ replaced by $d+m_1+m_2+\dots+m_l$ and the
appropriate interpretations of the $m_i^{(j)}$'s over all possible
types $T_i^{(j)}$ of rank $r_i^{(j)}$. This time, we obtain
\begin{multline} \label{eq:ExprD2}
(n-1)^{d}\Bigg(
\sum _{q=1} ^{l}\sum _{p=1} ^{m_q}\sum _{\al=2} ^{n}
{\prod _{i=1} ^{d}}
\frac {1} {n-\rk T_i-1}\binom {n-\rk T_i-1}{m_1^{(i)},m_2^{(i)},\dots,
m_n^{(i)}}\Bigg)\\
\times
\binom {n-1} {r_1^{(1)}}
\binom {n-1} {r_2^{(1)}}\cdots
\binom {n-1} {r_{m_1}^{(1)}}
\times\cdots\times
\binom {n-1} {r_1^{(q)}}
\cdots
\binom {n-\al-1} {r_{p}^{(q)}-\al}
\cdots
\binom {n-1} {r_{m_q}^{(q)}}\\
\times\cdots\times
\binom {n-1} {r_1^{(l)}}
\binom {n-1} {r_2^{(l)}}\cdots
\binom {n-1} {r_{m_l}^{(l)}}.
\end{multline}
The sum over $\al$ can be evaluated by means of the elementary
summation formula
\begin{equation*} 
\sum _{\al=2} ^{n}\binom {n-\al-1}{r-\al}=
\sum _{\al=2} ^{n}\binom {n-\al-1}{n-r-1}=\binom {n-2}{n-r}=\binom
{n-2}{r-2}.
\end{equation*}
Finally, we must sum the expressions \eqref{eq:ExprD1} and
\eqref{eq:ExprD2} over all non-negative integers $r_1^{(1)},\dots,r_{m_l}^{(l)}$ 
satisfying the equations
$$r_1^{(j)}+r_2^{(j)}+\dots+r_{m_j}^{(j)}=s_j,\quad j=1,2,\dots,l.$$
Once again, 
this is easily done by means of the multivariate version of the Chu--Vandermonde
summation. After some simplification, we obtain the formula in \eqref{eq:44b}.
\end{proof}

In view of the combinatorial realisation of 
$m$-divisible non-crossing partitions of type $D_n$ which we
described in Section~\ref{sec:6},
the special case $d=1$ of the above theorem has the following
enumerative consequence.

\begin{corollary} \label{cor:4}
Let $l$ be a positive integer, and let $s_1,s_2,\dots,s_l$ be
non-negative integers with $s_1+s_2+\dots+s_l=n$. 
The number of multi-chains 
$\pi_1\le \pi_2\le \dots\le \pi_{l-1}$ in the poset 
$\widetilde {NC}{}^m(D_n)$ 
with the property that $\rk(\pi_i)=s_1+s_2+\cdots +s_i,$ 
$i=1,2,\dots,l-1,$ and that the
number of non-zero blocks of $\pi_1$ of size $mi$ is $2b_i,$ $i=1,2,\dots,n,$
is given by
\begin{equation} \label{eq:45a}
\binom {b_1+b_2+\dots+b_n}{b_1,b_2,\dots,b_n}
\binom {m(n-1)} {s_2}\cdots
\binom {m(n-1)} {s_l},
\end{equation}
if $b_1+2b_2+\dots+nb_n<n-1,$ and
\begin{multline} \label{eq:45b}
2
\binom {b_1+b_2+\dots+b_n}{b_1,b_2,\dots,b_n}
\binom {m(n-1)} {s_2}\cdots
\binom {m(n-1)} {s_l}\\
+
\frac {m(n-1)} {b_1+b_2+\dots+b_n-1}\binom {b_1+b_2+\dots+b_n-1}
{b_1-1,b_2,\dots,b_n}\kern5cm\\
\times
\sum _{j=2} ^{l}\binom {m(n-1)} {s_2}\cdots
\binom {m(n-1)-1} {s_j-2}\cdots
\binom {m(n-1)} {s_l},
\end{multline}
if $b_1+2b_2+\dots+nb_n=n$.
\end{corollary}

\begin{remark}
The conditions in the statement of the corollary imply that
\begin{equation} \label{eq:sbD}
s_1+b_1+b_2+\dots+b_n=n.
\end{equation}
The reader should recall from Section~\ref{sec:6}, that non-zero
blocks of elements $\pi$ of $\widetilde {NC}{}^m(D_n)$ occur in pairs
since, with a block $B$ of $\pi$, also $\overline B$ is a block of
$\pi$. The condition $b_1+2b_2+\dots+nb_n<n-1$, which is required for
Formula~\eqref{eq:45a} to hold, implies that $\pi_1$ must contain a zero block
of size $2(n-b_1-2b_2-\dots-nb_n)$, while the equality
$b_1+2b_2+\dots+nb_n=n$, which is required for Formula~\eqref{eq:45b}
to hold, implies that $\pi_1$ contains no zero block. The extra
condition on zero blocks that are imposed on elements of
$\widetilde{NC}{}^m(D_n)$ implies that $b_1+2b_2+\dots+nb_n$ cannot
be equal to $n-1$.
\end{remark}

\begin{proof} 
Again, the arguments are completely analogous to those of the proof of
Corollary~\ref{cor:2}. Here we need
Theorem~\ref{thm:6} with
$d=1$, $l$ replaced by $l-1$, $s_i$ replaced by $s_{i+1}$,
$i=1,2,\dots,l-1$, $m_1=m_2=\dots=m_{l-1}=m$,
and $T_1$ of the type 
$$
D_{n-b_1-2b_2-\dots-nb_n}*A_1^{b_2}*A_2^{b_3}*\dots*A_{n-1}^{b_n}
$$
in the case that $b_1+2b_2+\dots+nb_n<n-1$, respectively
$$
A_1^{b_2}*A_2^{b_3}*\dots*A_{n-1}^{b_n}
$$
if not. So, depending on the case in which we are, we have to apply
\eqref{eq:44a}, respectively \eqref{eq:44b}. 
If we also take into account \eqref{eq:sbD}, then we
arrive at the claimed result after little manipulation. Since we have
done similar calculations already several times, the details are left
to the reader.
\end{proof}

This result is new even for $m=1$, that is, for the poset of
type $D_n$ non-crossing partitions of 
Athanasiadis and Reiner \cite{AtReAA}, and of Bessis and Corran
\cite{BeCoAA}. 
Not only does it imply all known results on these non-crossing partitions 
and their extension to $m$-divisible type $D_n$ non-crossing partitions 
due to Armstrong, it allows us as well to solve several open
enumeration problems on the $m$-divisible type $D_n$ non-crossing
partitions. We state these new results separately in the corollaries
below.

To begin with, if we set $l=2$ in Corollary~\ref{cor:4}, then we obtain
the following extension to $\widetilde {NC}{}^m(D_n)$ 
of Athanasiadis and Reiner's result \cite[Theorem~1.3]{AtReAA}
on the number of type $D_n$ non-crossing partitions with a given block
structure.

\begin{corollary} \label{cor:5A}
The number of all elements of
$\widetilde {NC}{}^m(D_n)$ which have
$2b_i$ non-zero blocks of size $mi$ equals
\begin{equation} \label{eq:48A}
\binom {b_1+b_2+\dots+b_n}{b_1,b_2,\dots,b_n}
\binom {m(n-1)} {b_1+b_2+\dots+b_n}
\end{equation}
if $b_1+2b_2+\dots+nb_n<n-1,$ and
\begin{multline} \label{eq:48B}
2
\binom {b_1+b_2+\dots+b_n}{b_1,b_2,\dots,b_n}
\binom {m(n-1)} {b_1+b_2+\dots+b_n}\\
+
\binom {b_1+b_2+\dots+b_n-1}
{b_1-1,b_2,\dots,b_n}
\binom {m(n-1)} {b_1+b_2+\dots+b_n-1}
\end{multline}
if $b_1+2b_2+\dots+nb_n=n$.
\end{corollary}

On the other hand, if we sum the expression \eqref{eq:45a},
respectively \eqref{eq:45b}, over all
$s_2,s_3,\dots,s_l$ with $s_2+s_3+\dots+s_l=n-s_1$, then we obtain
the following generalisation.

\begin{corollary} \label{cor:5}
The number of all multi-chains 
$\pi_1\le \pi_2\le \dots\le \pi_{l-1}$ in 
$\widetilde {NC}{}^m(D_n)$ in which
$\pi_1$ has $2b_i$ non-zero blocks of size $mi$ equals
\begin{equation} \label{eq:48a}
\binom {b_1+b_2+\dots+b_n}{b_1,b_2,\dots,b_n}
\binom {(l-1)m(n-1)} {b_1+b_2+\dots+b_n},
\end{equation}
if $b_1+2b_2+\dots+nb_n<n-1,$ and
\begin{multline} \label{eq:48b}
2
\binom {b_1+b_2+\dots+b_n}{b_1,b_2,\dots,b_n}
\binom {(l-1)m(n-1)} {b_1+b_2+\dots+b_n}\\
+
\binom {b_1+b_2+\dots+b_n-1}
{b_1-1,b_2,\dots,b_n}
\binom {(l-1)m(n-1)} {b_1+b_2+\dots+b_n-1}
\end{multline}
if $b_1+2b_2+\dots+nb_n=n$.
\end{corollary}

Next we sum the expressions \eqref{eq:45a} and \eqref{eq:45b} over all
possible $b_1,b_2,\dots,b_n$, that is, 
we sum \eqref{eq:45a} over $b_2+2b_3+\dots+(n-1)b_n< s_1-1$, 
and we sum the expression \eqref{eq:45b} over
$b_2+2b_3+\dots+(n-1)b_n=s_1$.
With the help of Lemma~\ref{lem:binsum} and the simple binomial
summation \eqref{eq:rksum1}, these sums can indeed be evaluated.
In this manner, we obtain the following result on rank-selected chain
enumeration in $\widetilde {NC}{}^m(D_n)$.

\begin{corollary} \label{cor:6}
The number of all multi-chains 
$\pi_1\le \pi_2\le \dots\le \pi_{l-1}$ in 
$\widetilde {NC}{}^m(D_n)\cong NC^m(D_n)$ 
where $\pi_i$ is of rank $s_1+s_2+\dots+s_i,$ $i=1,2,\dots,l-1,$ equals
\begin{multline} \label{eq:49}
2\binom {n-1} {s_1}
\binom {m(n-1)} {s_2}\cdots
\binom {m(n-1)} {s_l}\\
+
m\sum _{j=2} ^{l}\binom {n-1} {s_1}\binom {m(n-1)} {s_2}\cdots
\binom {m(n-1)-1} {s_j-2}\cdots
\binom {m(n-1)} {s_l}\\
+\binom {n-2} {s_1-2}
\binom {m(n-1)} {s_2}\cdots
\binom {m(n-1)} {s_l}.
\end{multline}
\end{corollary}

This formula extends Athanasiadis and Reiner's formula
\cite[Theorem~1.2(ii)]{AtReAA} from $NC(D_n)$ to
$\widetilde{NC}{}^m(D_n)$.
Setting $l=2$, we obtain a formula for the
number of all elements in
$\widetilde {NC}{}^m(D_n)\cong NC^m(D_n)$ 
with a given number of blocks (equivalently, of given rank);
cf.\ \cite[Theorem~4.6.3]{ArmDAA}.
Next, summing \eqref{eq:49} over all 
non-negative integers $s_1,s_2,\dots,s_l$ with $s_1+s_2+\dots+s_l=n$ by means
of the multivariate Chu--Vandermonde summation, we find 
that the total number of multi-chains 
$\pi_1\le \pi_2\le \dots\le \pi_{l-1}$ in
$\widetilde {NC}{}^m(D_n)\cong NC^m(D_n)$ is given by
\begin{multline} \label{eq:50}
2\binom {((l-1)m+1)(n-1)} {n}
+
\binom {((l-1)m+1)(n-1)} {n-1}\\=
\frac {2(l-1)m(n-1)+n} {n}
\binom {((l-1)m+1)(n-1)} {n-1},
\end{multline}
thus recovering the
formula \cite[Theorem~3.6.9]{ArmDAA} for the zeta polynomial of the poset of
generalised non-crossing partitions for the type $D_n$.
The special case $l=2$ of \eqref{eq:50} gives the well-known fact that the
cardinality of $\widetilde {NC}{}^m(D_n)\cong NC^m(D_n)$ 
is $\frac {2m(n-1)+n} {n}
\binom {(m+1)(n-1)} {n-1}$ (cf.\ \cite[Theorem~3.5.3]{ArmDAA}).

In the following section,
Corollary~\ref{cor:6} will enable us to provide a new proof of Armstrong's $F=M$
(Ex-)Conjecture in type $D_n$.


\section{Proof of the $F=M$ Conjecture for type $D$}
\label{sec:8} 

Armstrong's $F=M$ (Ex-)Conjecture \cite[Conjecture~5.3.2]{ArmDAA}, 
which extends an earlier
conjecture of Chapoton \cite{ChaFAA}, relates the ``$F$-triangle" of
the generalised cluster complex of Fomin and Reading \cite{FoReAA} to
the ``$M$-triangle" of Armstrong's generalised non-crossing
partitions. The $F$-triangle is a certain refined face count in the
generalised cluster complex. We do not give the definition here and,
instead, refer the reader to \cite{ArmDAA,KratCB}, because it will
not be important in what follows. It suffices to know that, again
fixing a finite root system $\Phi$ of rank $n$ and a positive integer
$m$, the $F$-triangle 
$F_\Phi^m(x,y)$ for the generalised cluster complex $\De^m(\Phi)$ is a 
polynomial in $x$ and $y$, and that it was computed in
\cite{KratCB} for all types. What we need here is that it was
shown in \cite[Sec.~11, Prop.~D]{KratCB} that
\begin{multline} \label{eq:F}
(1-xy)^n F^m_{D_n}\left(\frac {x(1+y)} {1-xy},\frac {xy}
{1-xy}\right)\\=
\sum _{r,s\ge0} ^{}x^sy^r\Bigg(
2\binom {n-1}{s-1}\binom {m(n-1)}r\binom {m(n-1)+s-r-1}{s-r}\\+
\binom {n-2}{s}\binom {m(n-1)}r\binom
{m(n-1)+s-r-1}{s-r}\\
+m\binom {n-1}{s-1}
\binom {m(n-1)-1}{r-2}\binom {m(n-1)+s-r-1}{s-r}\\-
m\binom {n-1}{s-1}\binom {m(n-1)}r\binom {m(n-1)+s-r-2}{s-r-2}
\Bigg).
\end{multline}

The ``$M$-triangle" of $NC^m(\Phi)$ is the polynomial defined by
$$
M^m_\Phi(x,y)=\sum _{u,w\in NC^m(\Phi)} ^{}\mu(u,w)\,x^{\rk u}y^{\rk w},$$
where $\mu(u,w)$ is the M\"obius function in $NC^m(\Phi)$.
It is called ``triangle" because the M\"obius function
$\mu(u,w)$ vanishes unless $u\le w$, and, thus, the only coefficients
in the polynomial which may be non-zero are the coefficients of $x^ky^l$
with $0\le k\le l\le n$.

An equivalent object is the {\it dual $M$-triangle}, which
is defined by
$$
(M^m_\Phi)^*(x,y)=
\sum _{u,w\in(NC^m(\Phi))^*}
^{}\mu^*(u,w)\,x^{\rk^*w}y^{\rk^*u},
$$
where $(NC^m(\Phi))^*$ denotes the poset {\rm dual} to $NC^m(\Phi)$
{\rm(}i.e., the poset which arises from $NC^m(\Phi)$ by reversing all
order relations{\rm)}, where $\mu^*$ denotes the M\"obius
function in $(NC^m(\Phi))^*$, and where $\rk^*$ denotes the rank
function in $(NC^m(\Phi))^*$. It is equivalent since, obviously, we
have
\begin{equation} \label{eq:02}
(M^m_\Phi)^*(x,y)=(xy)^nM^m_\Phi(1/x,1/y).
\end{equation}

Given this notation, Armstrong's $F=M$ (Ex-)Conjecture 
\cite[Conjecture~5.3.2]{ArmDAA} reads as
follows.

\begin{conjectureFM}
For any finite root system $\Phi$ of rank $n,$ we have
$$F^m_\Phi(x,y)=y^n\,M^m_\Phi\(\frac {1+y} {y-x},\frac {y-x} {y}\).$$
Equivalently, 
\begin{equation} \label{eq:FM}
(1-xy)^n F^m_\Phi\left(\frac {x(1+y)} {1-xy},\frac {xy}
{1-xy}\right)=
\sum _{u,w\in(NC^m(\Phi))^*}
^{}\mu^*(u,w)\,(-x)^{\rk^*w}(-y)^{\rk^*u}.
\end{equation}
\end{conjectureFM}

So, Equation~\eqref{eq:F} provides an expression for the left-hand
side of \eqref{eq:FM} for $\Phi=D_n$. 
With our result on rank-selected chain enumeration in $NC^m(D_n)$
given in Corollary~\ref{cor:6}, we are now able to calculate the
right-hand side of \eqref{eq:FM} directly. 
As we mentioned already in the Introduction, together with the
results from \cite{KratCB,KratCF},
this completes a computational case-by-case proof of
Conjecture~FM. A case-free proof had been found earlier by Tzanaki in
\cite{TzanAB}.

The only ingredient that we need for the proof is the well-known 
link between chain enumeration and the M\"obius function. 
(The reader should consult \cite[Sec.~3.11]{StanAP} for
more information on this topic.) 
Given a poset $P$ and two elements $u$ and $w$, $u\le w$,
in the poset, the {\it zeta polynomial\/} of the interval $[u,w]$,
denoted by $Z(u,w;z)$, is the number of (multi)chains from $u$ to
$w$ of length $z$. (It can be shown that this is indeed a
polynomial in $z$.) Then the M\"obius function of $u$ and $w$ is equal
to $\mu(u,w)=Z(u,w;-1)$.

\begin{proof}[Proof of Conjecture~FM in type $D_n$]
We now compute the right-hand side of \eqref{eq:FM}, that is,
$$
\sum _{u,w\in (NC^m(D_n))^*} ^{}\mu^*(u,w)
(-x)^{\rk^* w}(-y)^{\rk^* u}.$$
In order to compute the coefficient of $x^sy^r$ in this expression,
$$
(-1)^{r+s}\underset{\text{ with }\rk^* u=r\text{ and }\rk^* w=s}
{\sum _{u,w\in (NC^m(D_n))^*} ^{}}\mu^*(u,w),$$
we compute the sum of all corresponding zeta polynomials (in the
variable $z$), multiplied by $(-1)^{r+s}$,
$$
(-1)^{r+s}\underset{\text{ with }\rk^* u=r\text{ and }\rk^* w=s}
{\sum _{u,w\in (NC^m(D_n))^*} ^{}}Z(u,w;z),
$$
and then put $z=-1$. 

For computing this sum of zeta polynomials,
we must set
$l=z+2$, $n-s_1=s$, $s_l=r$, $s_2+s_3+\dots+s_{l-1}=s-r$ in 
\eqref{eq:49}, and
then sum the resulting expression over all possible
$s_2,s_3,\dots,s_{l-1}$. (The reader should keep in mind that
the roles of $s_1,s_2,\dots,s_l$ in Corollary~\ref{cor:6} have to be reversed,
since we are aiming at computing zeta polynomials in the poset {\it
dual\/} to $NC^m(D_n)$.)
By using the Chu--Vandermonde summation,
one obtains
\begin{multline*} 
2\binom {m(n-1)}{r}\binom {zm(n-1)}{s-r}\binom {n-1}{s-1}+
m\binom {m(n-1)-1}{r-2}
\binom {zm(n-1)}{s-r}\binom {n-1}{s-1}\\+
zm\binom {m(n-1)}{r}
\binom {zm(n-1)-1} {s-r-2}\binom {n-1}{s-1}+
\binom {m(n-1)}{r}\binom {zm(n-1)}{s-r}\binom {n-2}{s}.
\end{multline*}
If we put $z=-1$ in this expression and multiply it by $(-1)^{r+s}$,
then we obtain exactly the coefficient of $x^sy^r$ in \eqref{eq:F}.
\end{proof}

\section{A conjecture of Armstrong on maximal intervals containing a
random multichain}
\label{sec:8a} 

Given a finite root system of rank $n$,
Conjecture~3.5.13 in \cite{ArmDAA} says the following:

\medskip
{\it
If we choose an $l$-multichain uniformly at random from the set
\begin{equation} \label{eq:chains} 
\big\{\pi_1\le\pi_2\le\dots\le\pi_l:\pi_i\in NC^m(\Phi),i=1,\dots,l,
\text{ and }\rk(\pi_1)=i\big\},
\end{equation}
then the expected number of maximal intervals in $NC^m(\Phi)$
containing this multichain is
\begin{equation} \label{eq:Nar}
\frac {\Nar^m(\Phi,n-i)} {\Nar^1(\Phi,n-i)}, 
\end{equation}
where $\Nar^m(\Phi,i)$ is the $i$-th Fu\ss--Narayana number associated
to $NC^m(\Phi)$, that is, the number of elements of $NC^m(\Phi)$ of
rank $i$. In particular, this expected value is independent of $l$.}

We show in this section that, for types $A_n$ and
$B_n$, the conjecture follows easily from Edelman's \eqref{eq:35}
respectively Armstrong's \eqref{eq:42} (presumably, this fact constituted
the evidence for setting up the conjecture), while
an analogous computation using our new result \eqref{eq:49} 
demonstrates that 
it fails for type $D_n$. At the end of this section,
we comment on what we think happens for the exceptional types.

The computation of the expected value in the above conjecture can be
approached in the following way. One first observes that a maximal
interval in $NC^m(\Phi)$ is an interval between an element $\pi_0$ of
rank $0$ and the global maximum $(c;\ep,\dots,\ep)$. Therefore, to
compute the proposed expected value, we may count the number of chains
\begin{equation} \label{eq:Chains} 
\pi_0\le\pi_1\le\pi_2\le\dots\le\pi_l,\quad \rk(\pi_0)=0\text{ and }
\rk(\pi_1)=i,
\end{equation}
and divide this number by the total number of all chains in
\eqref{eq:chains}. Clearly, in types $A_n$, $B_n$, and $D_n$, this kind of
chain enumeration can be easily accessed by \eqref{eq:35},
\eqref{eq:42}, and \eqref{eq:49}, respectively.

\medskip
We begin with type $A_n$. By \eqref{eq:35}, the number of chains
\eqref{eq:Chains} equals
\begin{multline*}
\sum _{s_2+\dots+s_{l+1}=n-i} ^{}
\frac {1} {n+1}\binom {n+1} {0}
\binom {m(n+1)} {i}
\binom {m(n+1)} {s_2}\cdots
\binom {m(n+1)} {s_{l+1}}\\=
\frac {1} {n+1}\binom {m(n+1)} {i}
\binom {ml(n+1)} {n-i},
\end{multline*}
while the number of chains in \eqref{eq:chains} equals
\begin{multline*}
\sum _{s_2+\dots+s_{l+1}=n-i} ^{}
\frac {1} {n+1}\binom {n+1} {i}
\binom {m(n+1)} {s_2}\cdots
\binom {m(n+1)} {s_{l+1}}\\=
\frac {1} {n+1}\binom {n+1} {i}
\binom {ml(n+1)} {n-i}.
\end{multline*}
In both cases, we used the multivariate Chu--Vandermonde summation to
evaluate the sums over $s_2,\dots,s_{l+1}$.
The quotient of the two numbers is
$$\frac {\displaystyle\frac {1} {n+1}\binom {m(n+1)} {i}} 
{\displaystyle\frac {1} {n+1}\binom
{n+1} {i}}=
\frac {\displaystyle\frac {1} {n+1}\binom {n+1}{n-i}\binom {m(n+1)} {i}} 
{\displaystyle\frac {1} {n+1}\binom {n+1}{n-i}\binom {n+1} {i}},
$$
which by \eqref{eq:35} with $n$ replaced by $n+1$, $l=2$, $s_1=n-i$, and
$s_2=i$ agrees indeed with \eqref{eq:Nar} for $\Phi=A_n$. 

\medskip
For type $B_n$, there is an analogous computation using \eqref{eq:42},
the details of which we leave to the reader. The result is that the
desired expected value equals
$$\frac {\displaystyle\binom {mn} {i}} 
{\displaystyle\binom {n} {i}}=
\frac {\displaystyle\binom {n}{n-i}\binom {mn} {i}} 
{\displaystyle\binom {n}{n-i}\binom {n} {i}},
$$
which by \eqref{eq:42} with $l=2$, $s_1=n-i$, and
$s_2=i$ agrees indeed with \eqref{eq:Nar} for $\Phi=B_n$. 

\medskip
The analogous computation for type $D_n$ uses \eqref{eq:49}. 
The number of chains \eqref{eq:Chains} equals
\begin{align}
2&\sum _{s_2+\dots+s_{l+1}=n-i} ^{}
\binom {n-1} {0}
\binom {m(n-1)} {i}
\binom {m(n-1)} {s_2}\cdots
\binom {m(n-1)} {s_{l+1}}\notag\\
&\quad +
m\sum _{s_2+\dots+s_{l+1}=n-i} ^{}
\binom {n-1} {0}\binom {m(n-1)-1} {i-2}\binom {m(n-1)} {s_2}\cdots
\binom {m(n-1)} {s_{l+1}}\notag\\
&\quad +
m\sum _{j=3} ^{l+1}\sum _{s_2+\dots+s_{l+1}=n-i} ^{}
\binom {n-1} {0}\binom {m(n-1)} {i}\binom {m(n-1)} {s_2}\notag\\
&\kern7cm
\cdots
\binom {m(n-1)-1} {s_j-2}\cdots
\binom {m(n-1)} {s_{l+1}}\notag\\
&=
2\binom {m(n-1)} {i}
\binom {ml(n-1)} {n-i}+
m\binom {m(n-1)-1} {i-2}
\binom {ml(n-1)} {n-i}\notag\\
&\quad +
m(l-1)\binom {m(n-1)} {i}
\binom {ml(n-1)-1} {n-i-2},
\label{eq:1chain}
\end{align}
while the number of chains in \eqref{eq:chains} equals
\begin{align}
2&\sum _{s_2+\dots+s_{l+1}=n-i} ^{}
\binom {n-1} {i}
\binom {m(n-1)} {s_2}\cdots
\binom {m(n-1)} {s_{l+1}}\notag\\
&\quad +
m\sum _{j=2} ^{l+1}\sum _{s_2+\dots+s_{l+1}=n-i} ^{}
\binom {n-1} {i}\binom {m(n-1)} {s_2}\cdots
\binom {m(n-1)-1} {s_j-2}\cdots
\binom {m(n-1)} {s_{l+1}}\notag\\
&\quad +
\sum _{s_2+\dots+s_{l+1}=n-i} ^{}
\binom {n-2} {i-2}\binom {m(n-1)} {s_2}\cdots
\binom {m(n-1)} {s_{l+1}}\notag\\
&=
2\binom {n-1} {i}
\binom {ml(n-1)} {n-i}+
ml\binom {n-1} {i}
\binom {ml(n-1)-1} {n-i-2}+
\binom {n-2} {i-2}
\binom {ml(n-1)} {n-i}.
\label{eq:2chain}
\end{align}
The quotient of \eqref{eq:1chain} and \eqref{eq:2chain} gives the
desired expected value. It is, however, not independent of $l$, and
therefore Armstrong's conjecture does not hold for $\Phi=D_n$.

In the case that $\Phi$ is of exceptional type, then, as we outline in
the next section, the knowledge of the corresponding decomposition
numbers (see the Appendix) allows one to access the rank selected
chain enumeration. Using this, the approach for computing the expected
value proposed by Armstrong that we used above for the classical
types can be carried through as well for the exceptional types. We
have not done this, but we expect that, similarly to the case of $D_n$, 
for most exceptional types 
the expected value {\it will\/} depend on $l$, so that
Armstrong's conjecture will probably also fail in these cases.

\section{Chain enumeration in the poset of generalised non-crossing
partitions for the exceptional types}
\label{sec:9} 

Although it is not the main topic of our paper, we want to briefly demonstrate in this
section that the knowledge of the decomposition numbers 
also enables one to do refined enumeration in the generalised 
non-crossing partition posets $NC^m(\Phi)$ for {\it exceptional\/} root
systems $\Phi$ (of rank $n$). 
We restrict the following considerations to the 
rank-selected chain enumeration. This means that we want to count the
number of all multi-chains $\pi_1\le \pi_2\le \dots\le \pi_{l-1}$ in 
$NC^m(\Phi)$, where $\pi_i$ is of rank $s_1+s_2+\dots+s_i$,
$i=1,2,\dots,l-1$. Let us denote this number by
$R_\Phi(s_1,s_2,\dots,s_l)$, with $s_l=n-s_1-s_2-\dots-s_l$. 
Now, the considerations at the beginning of the
proof of Corollary~\ref{cor:2}, leading to the factorisation
\eqref{eq:bigfact} with rank constraints on the factors, are
also valid for $NC^m(\Phi)$ instead of $NC^m(A_{n-1})$, that is, they
are independent of the underlying root system. Hence, to determine
the number $R_\Phi(s_1,s_2,\dots,s_l)$, we have to count all possible
factorisations 
\begin{equation*} 
c=w_0^{(1)}
\big(u_1^{(2)}u_1^{(3)}\cdots u_1^{(l)}\big)
\big(u_2^{(2)}u_2^{(3)}\cdots u_2^{(l)}\big)\cdots
\big(u_m^{(2)}u_m^{(3)}\cdots u_m^{(l)}\big),
\end{equation*}
under the rank constraints
\eqref{eq:33} and $\ell_T(w_0^{(1)})=s_1$,
where $c$ is a Coxeter element in $W(\Phi)$. As we remarked in the
proof of Corollary~\ref{cor:2}, equivalently we may count all
factorisations
\begin{equation} \label{eq:bigfact1}
c=w_0^{(1)}
\big(u_1^{(2)}u_2^{(2)}\cdots u_m^{(2)}\big)
\big(u_1^{(3)}u_2^{(3)}\cdots u_m^{(3)}\big)\cdots
\big(u_1^{(l)}u_2^{(l)}\cdots u_m^{(l)}\big)
\end{equation}
which satisfy \eqref{eq:33} and $\ell_T(w_0^{(1)})=s_1$.
We can now obtain an explicit expression by fixing first the types of
$w_0^{(1)}$ and all the $u_i^{(j)}$'s. Under these constraints, the
number of factorisations \eqref{eq:bigfact1} is just the
corresponding decomposition number. Subsequently, we sum the
resulting expressions over all possible types. 

Before we are able to state the formula which we obtain in this way,
we need to recall some standard {\it integer} partition notation
(cf.\ e.g.\ \cite[Sec.~7.2]{StanBI}).
An {\it integer partition} $\la$ (with $n$ parts) is an $n$-tuple
$\la=(\la_1,\la_2,\dots,\la_n)$ of integers satisfying
$\la_1\ge\la_2\ge\dots\ge\la_n\ge0$. It is called an integer partition
{\it of $N$}, written in symbolic notation as $\la\vdash N$, if
$\la_1+\la_2+\dots+\la_n=N$. The number of parts (components) of $\la$ 
of size $i$ is denoted by $m_i(\la)$.

Then, making again use of the notation for the 
multinomial coefficient introduced in Lemma~\ref{lem:binsum}, 
the expression for 
$R_\Phi(s_1,s_2,\dots,s_l)$ which we obtain in the way described
above is
\begin{equation} \label{eq:ranksel}
{\sum _{} ^{}}{}^{\displaystyle\prime}\kern3pt
N_\Phi(T_0^{(1)},T_1^{(2)},T_2^{(2)},\dots,T_n^{(l)})
\prod _{j=2} ^{l}\binom m{m_1(\la^{(j)}),m_2(\la^{(j)}),\dots,
m_n(\la^{(j)})},
\end{equation}
where $\sum{}^{\textstyle\prime}$ is taken over all integer
partitions $\la^{(2)},\la^{(3)},\dots,\la^{(l)}$ satisfying
$\la^{(2)}\vdash s_2$, $\la^{(3)}\vdash s_3$,
\dots, $\la^{(l)}\vdash s_l$,
over all types $T_0^{(1)}$ with $\rk(T_0^{(1)})=s_1$,
and over all types $T_i^{(j)}$ with $\rk(T_i^{(j)})=\la_i^{(j)}$, 
$i=1,2,\dots,n$, $j=2,3,\dots,l$.

By way of example, using this formula and the values of the
decomposition numbers $N_{E_8}(\dots)$ given in Appendix~\ref{NE8} 
(and a computer), we obtain that the number 
$R_{E_8}(4,2,1,1)$ of all chains $\pi_1\le \pi_2\le \pi_3$ in $NC^m(E_8)$, 
where $\pi_1$ is of rank $4$, $\pi_2$ is of rank $6$,
and $\pi_3$ is of rank $7$, is given by
$$
 {\frac {75 {m^3} \left( 8055 m -1141 \right) } 2}
$$
(which, by the independence \eqref{Aa} of decomposition numbers from
the order of the types, is also equal to 
$R_{E_8}(4,1,2,1)$ and $R_{E_8}(4,1,1,2)$),
while the number $R_{E_8}(2,4,1,1)$ 
of all chains $\pi_1\le \pi_2\le \pi_3$ in $NC^m(E_8)$, 
where $\pi_1$ is of rank $2$, $\pi_2$ is of rank $6$,
and $\pi_3$ is of rank $7$, is given by
$$
 {\frac{75 {m^3} \left( 73125 {m^3} - 58950 {m^2}  + 15635 m-2154\right)}
      8} 
$$
(which is also equal to $R_{E_8}(2,1,4,1)$ and $R_{E_8}(2,1,1,4)$).

\section*{Acknowledgements}
The authors thank the anonymous referee for a very careful 
reading of the original manuscript.

\bigskip

\appendix

\section{The decomposition numbers for the exceptional types}

\subsection{The decomposition numbers for type $I_2(a)$ \cite[Sec.~13]{KratCB}}

We have\linebreak $N_{I_2(a)}(I_2(a))=1$, $N_{I_2(a)}(A_1,A_1)=a$, 
$N_{I_2(a)}(A_1)=a$, $N_{I_2(a)}(\emptyset)=1$,
all other numbers $N_{I_2(a)}(T_1,T_2,\dots,T_d)$ being zero.

\subsection{The decomposition numbers for type $H_3$ \cite[Sec.~14]{KratCB}}

We have $N_{H_3}(H_3)=1$, $N_{H_3}(A_1^2,A_1)=5$,
$N_{H_3}(A_2,A_1)=5$, $N_{H_3}(I_2(5),A_1)=5$,
$N_{H_3}(A_1,A_1,A_1)=50$, plus the assignments implied by
\eqref{Aa} and \eqref{Ab},
all other numbers $N_{H_3}(T_1,T_2,\dots,T_d)$ being zero.

\subsection{The decomposition numbers for type $H_4$ \cite[Sec.~15]{KratCB}}

We have $N_{H_4}(H_4)=1$, $N_{H_4}(A_1*A_2,A_1)=15$,
$N_{H_4}(A_3,A_1)=15$, 
$N_{H_4}(H_3,A_1)=15$,
$N_{H_4}(A_1*I_2(5),A_1)=15$,
$N_{H_4}(A_1^2,A_1^2)=30$,
$N_{H_4}(A_1^2,A_2)=30$,
$N_{H_4}(A_1^2,I_2(5))=15$,
$N_{H_4}(A_2,A_2)=5$,
$N_{H_4}(A_2,I_2(5))=15$,
$N_{H_4}(I_2(5),I_2(5))=3$,
$N_{H_4}(A_1^2,A_1,A_1)=225$,
$N_{H_4}(A_2,A_1,A_1)=150$,
$N_{H_4}(I_2(5),A_1,A_1)=90$,
$N_{H_4}(A_1,A_1,A_1,A_1)=1350$, 
plus the assignments implied by
\eqref{Aa} and \eqref{Ab},
all other numbers $N_{H_4}(T_1,T_2,\dots,T_d)$ being zero.

\subsection{The decomposition numbers for type $F_4$ \cite[Sec.~16]{KratCB}}

We have $N_{F_4}(F_4)=1$, $N_{F_4}(A_1*A_2,A_1)=12$,
$N_{F_4}(B_3,A_1)=12$, 
$N_{F_4}(A_1^2,A_1^2)=12$,
$N_{F_4}(A_1^2,B_2)=12$,
$N_{F_4}(A_2,A_2)=16$,
$N_{F_4}(B_2,B_2)=3$,
$N_{F_4}(A_1^2,A_1,A_1)=72$,
$N_{F_4}(A_2,A_1,A_1)=48$,
$N_{F_4}(B_2,A_1,A_1)=36$,
$N_{F_4}(A_1,A_1,A_1,A_1)=432$, 
plus the assignments implied by
\eqref{Aa} and \eqref{Ab},
all other numbers $N_{F_4}(T_1,T_2,\dots,T_d)$ being zero.

\subsection{The decomposition numbers for type $E_6$ \cite[Sec.~17]{KratCB}}

We have $N_{E_6}(E_6)=1$, 
 $N_{E_6}(A_1*A_2^2, A_1) = 6$,
 $N_{E_6}(A_1*A_4, A_1) = 12$, $N_{E_6}(A_5, A_1) = 6$,  
 $N_{E_6}(D_5, A_1) = 12$, 
 $N_{E_6}(A_1^2*A_2, A_2) = 36$, $N_{E_6}(A_2^2, A_2) = 8$, 
 $N_{E_6}(A_1*A_3, A_2) = 24$, $N_{E_6}(A_4, A_2) = 24$, $N_{E_6}(D_4, A_2) = 4$, 
 $N_{E_6}(A_1^2*A_2, A_1^2) = 18$, 
 $N_{E_6}(A_1*A_3, A_1^2) = 36$, $N_{E_6}(A_4, A_1^2) = 36$, $N_{E_6}(D_4, A_1^2) = 18$, 
 $N_{E_6}(A_1^3, A_1^3) = 12$, $N_{E_6}(A_1*A_2, A_1^3) = 24$, $N_{E_6}(A_1*A_2, A_1*A_2) = 48$, 
 $N_{E_6}(A_3, A_1^3) = 36$, $N_{E_6}(A_3, A_1*A_2) = 72$, $N_{E_6}(A_3, A_3) = 27$, 
 $N_{E_6}(A_1^2*A_2, A_1, A_1) = 144$, $N_{E_6}(A_2^2, A_1, A_1) = 24$, 
 $N_{E_6}(A_1*A_3, A_1, A_1) = 144$, $N_{E_6}(A_4, A_1, A_1) = 144$, $N_{E_6}(D_4, A_1, A_1) = 48$, 
 $N_{E_6}(A_1^3, A_1^2, A_1) = 180$, $N_{E_6}(A_1^3, A_2, A_1) = 168$, 
 $N_{E_6}(A_1*A_2, A_1^2, A_1) = 360$,
 $N_{E_6}(A_1*A_2, A_2, A_1) = 336$, 
 $N_{E_6}(A_3, A_1^2, A_1) = 378$, $N_{E_6}(A_3, A_2,\break A_1) = 180$, 
 $N_{E_6}(A_1^2, A_1^2, A_1^2) = 432$, $N_{E_6}(A_2, A_1^2, A_1^2) = 504$, 
 $N_{E_6}(A_2, A_2, A_1^2) = 288$, $N_{E_6}(A_2, A_2, A_2) = 160$,
 $N_{E_6}(A_1^2, A_1^2, A_1, A_1) = 2376$, 
 $N_{E_6}(A_2, A_1^2, A_1, A_1) = 1872$,\break $N_{E_6}(A_2, A_2, A_1, A_1) = 1056$,
 $N_{E_6}(A_1^3, A_1, A_1, A_1) = 864$, $N_{E_6}(A_1*A_2, A_1, A_1, A_1) = 1728$, 
 $N_{E_6}(A_3, A_1, A_1, A_1) = 1296$, 
 $N_{E_6}(A_1^2, A_1, A_1, A_1, A_1) = 10368$, $N_{E_6}(A_2, A_1, A_1,
A_1,\break A_1) = 6912$,
 $N_{E_6}(A_1, A_1, A_1, A_1, A_1, A_1) = 41472$,
plus the assignments implied by
\eqref{Aa} and \eqref{Ab},
all other numbers $N_{E_6}(T_1,T_2,\dots,T_d)$ being zero.

\subsection{The decomposition numbers for type $E_7$
\cite[Sec.~6]{KratCF}}

We have $N_{E_7}(E_7)=1$, 
$N_{E_7}(E_6, A_1) = 9$,
$N_{E_7}(D_6, A_1) = 9$,
$N_{E_7}(A_6, A_1) = 9$,
$N_{E_7}(A_1*D_5, A_1) = 9$,
$N_{E_7}(A_1*A_5, A_1) = 9$,
$N_{E_7}(A_2*D_4, A_1) = 0$,
$N_{E_7}(A_2*A_4, A_1) = 9$,
$N_{E_7}(A_1^2*D_4, A_1) = 0$,
$N_{E_7}(A_1^2*A_4, A_1) = 0$,
$N_{E_7}(A_3^2, A_1) = 0$,
$N_{E_7}(A_1*A_2*A_3, A_1) = 9$,
$N_{E_7}(A_1^3*A_3, A_1) = 0$, 
$N_{E_7}(A_2^3, A_1) = 0$,
$N_{E_7}(A_1^2*A_2^2, A_1) = 0$,
$N_{E_7}(A_1^4*A_2, A_1) = 0$, 
$N_{E_7}(A_1^6, A_1) = 0$, 
$N_{E_7}(D_5, A_2) = 18$,
$N_{E_7}(A_5, A_2) = 30$,
$N_{E_7}(A_1*A_4, A_2) = 54$,
$N_{E_7}(A_1*D_4, A_2) = 9$,
$N_{E_7}(A_2*A_3, A_2) = 36$,
$N_{E_7}(A_1^2*A_3, A_2) = 36$,
$N_{E_7}(A_1*A_2^2, A_2) = 36$,
$N_{E_7}(A_1^3*A_2, A_2) = 12$,
$N_{E_7}(A_1^5, A_2) = 0$,
$N_{E_7}(D_5, A_1^2) = 54$,
$N_{E_7}(A_5, A_1^2) = 63$,
$N_{E_7}(A_1*D_4, A_1^2) = 27$,
$N_{E_7}(A_1*A_4, A_1^2) = 81$,
$N_{E_7}(A_2*A_3, A_1^2) = 27$,
$N_{E_7}(A_1^2*A_3, A_1^2) = 27$,
$N_{E_7}(A_1*A_2^2, A_1^2) = 27$, 
$N_{E_7}(A_1^3*A_2, A_1^2) = 9$,
$N_{E_7}(A_1^5, A_1^2) = 0$, 
$N_{E_7}(D_5, A_1, A_1) = 162$,
$N_{E_7}(A_5, A_1, A_1) = 216$,
$N_{E_7}(A_1*D_4, A_1, A_1) = 81$,
$N_{E_7}(A_1*A_4, A_1, A_1) = 324$,
$N_{E_7}(A_2*A_3, A_1, A_1) = 162$,
$N_{E_7}(A_1^2*A_3, A_1, A_1) = 162$,
$N_{E_7}(A_1*A_2^2, A_1, A_1) = 162$,
$N_{E_7}(A_1^3*A_2, A_1, A_1) = 54$,
$N_{E_7}(A_1^5, A_1, A_1) = 0$, 
$N_{E_7}(D_4,A_3)=9$, 
$N_{E_7}(A_4,A_3)=54$,
$N_{E_7}(A_1*A_3, A_3) = 135$,
$N_{E_7}(A_2^2, A_3) = 54$,
$N_{E_7}(A_1^2*A_2, A_3) = 162$,
$N_{E_7}(A_1^4, A_3) = 27$,
$N_{E_7}(D_4, A_1*A_2) = 45$,
$N_{E_7}(A_4, A_1*A_2) = 162$,
$N_{E_7}(A_1*A_3, A_1*A_2) = 243$,
$N_{E_7}(A_2^2, A_1*A_2) = 54$,
$N_{E_7}(A_1^2*A_2, A_1*A_2) = 162$,
$N_{E_7}(A_1^4, A_1*A_2) = 27$,
$N_{E_7}(D_4, A_1^3) = 30$,
$N_{E_7}(A_4, A_1^3) = 99$,
$N_{E_7}(A_1*A_3, A_1^3) = 126$,
$N_{E_7}(A_2^2, A_1^3) = 18$,
$N_{E_7}(A_1^2*A_2, A_1^3) = 54$,
$N_{E_7}(A_1^4, A_1^3) = 9$, 
$N_{E_7}(D_4, A_2, A_1) = 81$,
$N_{E_7}(A_4, A_2, A_1) = 378$,
$N_{E_7}(A_1*A_3, A_2, A_1) = 783$,
$N_{E_7}(A_2^2, A_2,  A_1) = 270$,
$N_{E_7}(A_1^2*A_2, A_2, A_1) = 810$,
$N_{E_7}(A_1^4, A_2, A_1) = 135$,
$N_{E_7}(D_4, A_1^2, A_1) = 243$,
$N_{E_7}(A_4, A_1^2, A_1) = 891$,
$N_{E_7}(A_1*A_3, A_1^2, A_1) = 1377$,
$N_{E_7}(A_2^2,\break A_1^2, A_1) = 324$,
$N_{E_7}(A_1^2*A_2, A_1^2, A_1) = 972$,
$N_{E_7}(A_1^4, A_1^2, A_1) = 162$,
$N_{E_7}(D_4, A_1, A_1,\break A_1) = 729$,
$N_{E_7}(A_4,  A_1, A_1, A_1) = 2916$,
$N_{E_7}(A_1*A_3, A_1, A_1, A_1) = 5103$,
$N_{E_7}(A_2^2, A_1,\break A_1, A_1) = 1458$,
$N_{E_7}(A_1^2*A_2, A_1, A_1, A_1) = 4374$,
$N_{E_7}(A_1^4, A_1, A_1, A_1) = 729$,
$N_{E_7}(A_3,\break A_3, A_1) = 486$,
$N_{E_7}(A_3, A_1*A_2, A_1) = 1458$,
$N_{E_7}(A_3, A_1^3, A_1) = 891$,
$N_{E_7}(A_1*A_2, A_1*A_2, A_1) = 2430$,
$N_{E_7}(A_1*A_2, A_1^3, A_1) = 1215$,
$N_{E_7}(A_1^3, A_1^3, A_1) = 540$,
$N_{E_7}(A_3, A_2, A_2) = 432$,
$N_{E_7}(A_1*A_2, A_2, A_2) = 1188$,
$N_{E_7}(A_1^3, A_2, A_2) = 711$,
$N_{E_7}(A_3, A_2, A_1^2) = 1053$,
$N_{E_7}(A_1*A_2, A_2, A_1^2) = 2349$,
$N_{E_7}(A_1^3, A_2, A_1^2) = 1323$,
$N_{E_7}(A_3, A_1^2, A_1^2) = 2430$,
$N_{E_7}(A_1*A_2, A_1^2, A_1^2) = 3402$,
$N_{E_7}(A_1^3, A_1^2, A_1^2) = 1539$, 
$N_{E_7}(A_3, A_2, A_1, A_1) = 3402$,
$N_{E_7}(A_1*A_2, A_2, A_1, A_1) = 8262$,
$N_{E_7}(A_1^3, A_2, A_1, A_1) = 4779$,
$N_{E_7}(A_3, A_1^2, A_1, A_1) = 8019$,\break
$N_{E_7}(A_1*A_2, A_1^2, A_1, A_1) = 13851$,
$N_{E_7}(A_1^3, A_1^2, A_1, A_1) = 7047$,
$N_{E_7}(A_3, A_1, A_1, A_1,\break A_1) = 26244$,
$N_{E_7}(A_1*A_2, A_1, A_1, A_1, A_1) = 52488$,
$N_{E_7}(A_1^3, A_1, A_1, A_1, A_1) = 28431$,
$N_{E_7}(A_2, A_2, A_2, A_1) = 2916$,
$N_{E_7}(A_2, A_2, A_1^2, A_1) = 6561$,
$N_{E_7}(A_2, A_1^2, A_1^2, A_1) = 13122$,
$N_{E_7}(A_1^2, A_1^2, A_1^2, A_1) = 19683$, 
$N_{E_7}(A_2, A_2, A_1, A_1, A_1) = 21870$,
$N_{E_7}(A_2, A_1^2, A_1, A_1,\break A_1) = 45927$,
$N_{E_7}(A_1^2, A_1^2, A_1, A_1, A_1) = 78732$,
$N_{E_7}(A_2, A_1, A_1, A_1, A_1, A_1) = 157464$,
$N_{E_7}(A_1^2, A_1, A_1, A_1, A_1, A_1) = 295245$,
$N_{E_7}(A_1,A_1,A_1,A_1,A_1,A_1,A_1)=1062882$,
plus the assignments implied by
\eqref{Aa} and \eqref{Ab},
all other numbers $N_{E_7}(T_1,T_2,\dots,T_d)$ being zero.

\subsection{The decomposition numbers for type $E_8$
\cite[Sec.~7]{KratCF}} \label{NE8}

We have $N_{E_8}(E_8)= 1$, 
$N_{E_8}(E_7, A_1)= 15$,
$N_{E_8}(D_7, A_1)= 15$,
$N_{E_8}(A_7, A_1)= 15$,
$N_{E_8}(A_1*E_6, A_1)= 15$,
$N_{E_8}(A_1*D_6, A_1)= 0$,
$N_{E_8}(A_1*A_6, A_1)= 15$,
$N_{E_8}(A_2*D_5, A_1)= 15$,
$N_{E_8}(A_2*A_5, A_1)= 0$,
$N_{E_8}(A_1^2*D_5, A_1)= 0$,
$N_{E_8}(A_1^2*A_5, A_1)= 0$,
$N_{E_8}(A_3*D_4, A_1)= 0$,
$N_{E_8}(A_3*A_4, A_1)= 15$,
$N_{E_8}(A_1*A_2*D_4, A_1)= 0$,
$N_{E_8}(A_1*A_2*A_4, A_1)= 15$,
$N_{E_8}(A_1^3*D_4, A_1)= 0$,
$N_{E_8}(A_1^3*A_4, A_1)= 0$,
$N_{E_8}(A_1*A_3^2, A_1)= 0$,
$N_{E_8}(A_2^2*A_3, A_1)= 0$,
$N_{E_8}(A_1^2*A_2*A_3, A_1)= 0$,
$N_{E_8}(A_1^4*A_3, A_1)= 0$,
$N_{E_8}(A_1*A_2^3, A_1)= 0$,
$N_{E_8}(A_1^3*A_2^2, A_1)= 0$,
$N_{E_8}(A_1^5*A_2, A_1)= 0$,
$N_{E_8}(A_1^7, A_1)= 0$,
$N_{E_8}(E_6, A_2)= 20$,
$N_{E_8}(D_6, A_2)= 15$,
$N_{E_8}(A_6, A_2)= 60$,
$N_{E_8}(A_1*D_5, A_2)= 60$,
$N_{E_8}(A_1*A_5, A_2)= 60$,
$N_{E_8}(A_2*D_4, A_2)= 20$,
$N_{E_8}(A_2*A_4, A_2)= 90$,
$N_{E_8}(A_3^2, A_2)= 45$,
$N_{E_8}(A_1^2*D_4, A_2)= 0$,
$N_{E_8}(A_1^2*A_4, A_2)= 90$,
$N_{E_8}(A_1*A_2*A_3, A_2)= 90$,
$N_{E_8}(A_1^3*A_3, A_2)= 0$,
$N_{E_8}(A_2^3, A_2)= 0$,
$N_{E_8}(A_1^2*A_2^2, A_2)= 45$,
$N_{E_8}(A_1^4*A_2, A_2)= 0$,
$N_{E_8}(A_1^6, A_2)= 0$,
$N_{E_8}(E_6, A_1^2)= 45$,
$N_{E_8}(D_6, A_1^2)= 90$,
$N_{E_8}(A_6, A_1^2)= 135$,
$N_{E_8}(A_1*D_5, A_1^2)= 135$,
$N_{E_8}(A_1*A_5, A_1^2)= 135$,
$N_{E_8}(A_2*D_4, A_1^2)= 45$,
$N_{E_8}(A_2*A_4, A_1^2)= 90$,
$N_{E_8}(A_3^2, A_1^2)= 45$,
$N_{E_8}(A_1^2*D_4, A_1^2)= 0$,
$N_{E_8}(A_1^2*A_4, A_1^2)= 90$,
$N_{E_8}(A_1*A_2*A_3, A_1^2)= 90$,
$N_{E_8}(A_1^3*A_3, A_1^2)= 0$,
$N_{E_8}(A_2^3, A_1^2)= 0$,
$N_{E_8}(A_1^2*A_2^2, A_1^2)= 45$,
$N_{E_8}(A_1^4*A_2, A_1^2)= 0$,
$N_{E_8}(A_1^6, A_1^2)= 0$,
$N_{E_8}(E_6, A_1, A_1)= 150$,
$N_{E_8}(D_6, A_1, A_1)= 225$,
$N_{E_8}(A_6, A_1, A_1)= 450$,
$N_{E_8}(A_1*D_5, A_1, A_1)= 450$,
$N_{E_8}(A_1*A_5, A_1, A_1)= 450$,
$N_{E_8}(A_2*D_4, A_1, A_1)= 150$,
$N_{E_8}(A_2*A_4, A_1, A_1)= 450$,
$N_{E_8}(A_3^2, A_1, A_1)= 225$,
$N_{E_8}(A_1^2*D_4, A_1, A_1)= 0$,
$N_{E_8}(A_1^2*A_4, A_1, A_1)= 450$,
$N_{E_8}(A_1*A_2*A_3, A_1, A_1)= 450$,
$N_{E_8}(A_1^3*A_3, A_1, A_1)= 0$,
$N_{E_8}(A_2^3, A_1, A_1)= 0$,
$N_{E_8}(A_1^2*A_2^2, A_1, A_1)= 225$,
$N_{E_8}(A_1^4*A_2, A_1, A_1)= 0$,
$N_{E_8}(A_1^6, A_1, A_1)= 0$,
$N_{E_8}(D_5, A_3)= 45$,
$N_{E_8}(A_5, A_3)= 90$,
$N_{E_8}(A_1*A_4, A_3)= 315$,
$N_{E_8}(A_1*D_4, A_3)= 45$,
$N_{E_8}(A_2*A_3, A_3)= 270$,
$N_{E_8}(A_1^2*A_3, A_3)= 270$,
$N_{E_8}(A_1*A_2^2, A_3)= 225$,
$N_{E_8}(A_1^3*A_2, A_3)= 225$,
$N_{E_8}(A_1^5, A_3)= 0$,
$N_{E_8}(D_5, A_1*A_2)= 195$,
$N_{E_8}(A_5, A_1*A_2)= 390$,
$N_{E_8}(A_1*A_4, A_1*A_2)= 690$,
$N_{E_8}(A_1*D_4, A_1*A_2)= 195$,
$N_{E_8}(A_2*A_3, A_1*A_2)= 495$,
$N_{E_8}(A_1^2*A_3, A_1*A_2)= 495$,
$N_{E_8}(A_1*A_2^2, A_1*A_2)= 300$,
$N_{E_8}(A_1^3*A_2, A_1*A_2)= 300$,
$N_{E_8}(A_1^5, A_1*A_2)= 0$,
$N_{E_8}(D_5, A_1^3)= 150$,
$N_{E_8}(A_5, A_1^3)= 300$,
$N_{E_8}(A_1*A_4, A_1^3)= 375$,
$N_{E_8}(A_1*D_4, A_1^3)= 150$,
$N_{E_8}(A_2*A_3, A_1^3)= 225$,
$N_{E_8}(A_1^2*A_3, A_1^3)= 225$,
$N_{E_8}(A_1*A_2^2, A_1^3)= 75$,
$N_{E_8}(A_1^3*A_2, A_1^3)= 75$,
$N_{E_8}(A_1^5, A_1^3)= 0$,
$N_{E_8}(D_5, A_2, A_1)= 375$,
$N_{E_8}(A_5, A_2, A_1)= 750$,
$N_{E_8}(A_1*A_4, A_2, A_1)= 1950$,
$N_{E_8}(A_1*D_4, A_2, A_1)= 375$,
$N_{E_8}(A_2*A_3, A_2, A_1)= 1575$,
$N_{E_8}(A_1^2*A_3, A_2, A_1)= 1575$,
$N_{E_8}(A_1*A_2^2, A_2, A_1)= 1200$,
$N_{E_8}(A_1^3*A_2, A_2, A_1)= 1200$,
$N_{E_8}(A_1^5, A_2, A_1)= 0$,
$N_{E_8}(D_5, A_1^2, A_1)= 1125$,
$N_{E_8}(A_5, A_1^2, A_1)= 2250$,
$N_{E_8}(A_1*A_4, A_1^2, A_1)= 3825$,
$N_{E_8}(A_1*D_4, A_1^2, A_1)= 1125$,
$N_{E_8}(A_2*A_3, A_1^2, A_1)= 2700$,
$N_{E_8}(A_1^2*A_3, A_1^2, A_1)= 2700$,
$N_{E_8}(A_1*A_2^2, A_1^2, A_1)= 1575$,
$N_{E_8}(A_1^3*A_2, A_1^2, A_1)= 1575$,
$N_{E_8}(A_1^5, A_1^2, A_1)= 0$,
$N_{E_8}(D_5, A_1, A_1, A_1)= 3375$,
$N_{E_8}(A_5, A_1, A_1, A_1)= 6750$,
$N_{E_8}(A_1*A_4, A_1, A_1, A_1)= 13500$,
$N_{E_8}(A_1*D_4, A_1, A_1, A_1)= 3375$,
$N_{E_8}(A_2*A_3, A_1, A_1, A_1)= 10125$,
$N_{E_8}(A_1^2*A_3, A_1, A_1, A_1)= 10125$,
$N_{E_8}(A_1*A_2^2, A_1, A_1, A_1)= 6750$,
$N_{E_8}(A_1^3*A_2, A_1, A_1, A_1)= 6750$,
$N_{E_8}(A_1^5, A_1, A_1, A_1)= 0$,
$N_{E_8}(D_4, D_4)= 5$,
$N_{E_8}(D_4, A_4)= 15$,
$N_{E_8}(A_4, A_4)= 138$,
$N_{E_8}(D_4, A_1*A_3)= 105$,
$N_{E_8}(A_4, A_1*A_3)= 390$,
$N_{E_8}(A_1*A_3, A_1*A_3)= 1155$,
$N_{E_8}(D_4, A_2^2)= 35$,
$N_{E_8}(A_4, A_2^2)= 180$,
$N_{E_8}(A_1*A_3, A_2^2)= 360$,
$N_{E_8}(A_2^2, A_2^2)= 95$,
$N_{E_8}(D_4, A_1^2*A_2)= 135$,
$N_{E_8}(A_4, A_1^2*A_2)= 630$,
$N_{E_8}(A_1*A_3, A_1^2*A_2)= 1035$,
$N_{E_8}(A_2^2, A_1^2*A_2)= 270$,
$N_{E_8}(A_1^2*A_2, A_1^2*A_2)= 495$,
$N_{E_8}(D_4, A_1^4)= 30$,
$N_{E_8}(A_4, A_1^4)= 165$,
$N_{E_8}(A_1*A_3, A_1^4)= 255$,
$N_{E_8}(A_2^2, A_1^4)= 60$,
$N_{E_8}(A_1^2*A_2, A_1^4)= 135$,
$N_{E_8}(A_1^4, A_1^4)= 30$,
$N_{E_8}(D_4, A_3, A_1)= 225$,
$N_{E_8}(A_4, A_3, A_1)= 1215$,
$N_{E_8}(A_1*A_3, A_3, A_1)= 4050$,
$N_{E_8}(A_2^2, A_3, A_1)= 1575$,
$N_{E_8}(A_1^2*A_2, A_3, A_1)= 5400$,
$N_{E_8}(A_1^4, A_3, A_1)= 1350$,
$N_{E_8}(D_4,\break A_1*A_2, A_1)= 975$,
$N_{E_8}(A_4, A_1*A_2, A_1)= 4590$,
$N_{E_8}(A_1*A_3, A_1*A_2, A_1)= 10800$,
$N_{E_8}(A_2^2, A_1*A_2, A_1)= 3450$,
$N_{E_8}(A_1^2*A_2, A_1*A_2, A_1)= 9900$,
$N_{E_8}(A_1^4, A_1*A_2, A_1)= 2475$,
$N_{E_8}(D_4, A_1^3, A_1)= 750$,
$N_{E_8}(A_4, A_1^3, A_1)= 3375$,
$N_{E_8}(A_1*A_3, A_1^3, A_1)= 6750$,
$N_{E_8}(A_2^2, A_1^3, A_1)= 1875$,
$N_{E_8}(A_1^2*A_2, A_1^3, A_1)= 4500$,
$N_{E_8}(A_1^4, A_1^3, A_1)= 1125$,
$N_{E_8}(D_4,\break A_2, A_2)= 175$,
$N_{E_8}(A_4, A_2, A_2)= 1140$,
$N_{E_8}(A_1*A_3,  A_2, A_2)= 3300$,
$N_{E_8}(A_2^2, A_2, A_2)= 1300$,
$N_{E_8}(A_1^2*A_2, A_2, A_2)= 4500$,
$N_{E_8}(A_1^4, A_2, A_2)= 1125$,
$N_{E_8}(D_4, A_2, A_1^2)= 675$,
$N_{E_8}(A_4, A_2, A_1^2)= 3015$,
$N_{E_8}(A_1*A_3, A_2, A_1^2)= 8550$,
$N_{E_8}(A_2^2, A_2, A_1^2)= 2925$,
$N_{E_8}(A_1^2*A_2, A_2, A_1^2)= 9000$,
$N_{E_8}(A_1^4, A_2, A_1^2)= 2250$,
$N_{E_8}(D_4,  A_1^2, A_1^2)= 1800$,
$N_{E_8}(A_4, A_1^2,\break A_1^2)= 8640$,
$N_{E_8}(A_1*A_3, A_1^2, A_1^2)= 17550$,
$N_{E_8}(A_2^2, A_1^2,  A_1^2)= 5175$,
$N_{E_8}(A_1^2*A_2, A_1^2,\break A_1^2)= 13500$,
$N_{E_8}(A_1^4, A_1^2, A_1^2)= 3375$,
$N_{E_8}(D_4, A_2, A_1,  A_1)= 1875$,
$N_{E_8}(A_4, A_2, A_1,\break A_1)= 9450$,
$N_{E_8}(A_1*A_3, A_2, A_1, A_1)= 27000$,
$N_{E_8}(A_2^2, A_2,  A_1, A_1)= 9750$,
$N_{E_8}(A_1^2*A_2, A_2, A_1, A_1)= 31500$,
$N_{E_8}(A_1^4, A_2, A_1, A_1)= 7875$,
$N_{E_8}(D_4, A_1^2, A_1, A_1)= 5625$,
$N_{E_8}(A_4, A_1^2, A_1, A_1)= 26325$,
$N_{E_8}(A_1*A_3, A_1^2, A_1, A_1)= 60750$,
$N_{E_8}(A_2^2, A_1^2, A_1, A_1)= 19125$,
$N_{E_8}(A_1^2*A_2, A_1^2, A_1, A_1)= 54000$,
$N_{E_8}(A_1^4, A_1^2, A_1, A_1)= 13500$,
$N_{E_8}(D_4, A_1, A_1,\break A_1, A_1)= 16875$,
$N_{E_8}(A_4, A_1, A_1, A_1, A_1)= 81000$,
$N_{E_8}(A_1*A_3,  A_1, A_1, A_1, A_1)= 202500$,
$N_{E_8}(A_2^2, A_1, A_1, A_1, A_1)= 67500$,
$N_{E_8}(A_1^2*A_2, A_1, A_1, A_1,  A_1)= 202500$,\break
$N_{E_8}(A_1^4, A_1, A_1, A_1, A_1)= 50625$,
$N_{E_8}(A_3, A_3, A_2)= 1350$,
$N_{E_8}(A_3,  A_1*A_2, A_2)= 5175$,
$N_{E_8}(A_3, A_1^3, A_2)= 3825$,
$N_{E_8}(A_1*A_2, A_1*A_2, A_2)= 15000$,
$N_{E_8}(A_1*A_2, A_1^3, A_2)= 9825$,
$N_{E_8}(A_1^3, A_1^3, A_2)= 6000$,
$N_{E_8}(A_3, A_3, A_1^2)= 4050$,
$N_{E_8}(A_3,  A_1*A_2, A_1^2)= 13500$,
$N_{E_8}(A_3, A_1^3, A_1^2)= 9450$,
$N_{E_8}(A_1*A_2, A_1*A_2, A_1^2)= 30825$,
$N_{E_8}(A_1*A_2, A_1^3, A_1^2)= 17325$,
$N_{E_8}(A_1^3, A_1^3, A_1^2)= 7875$,
$N_{E_8}(A_3, A_3, A_1, A_1)= 12150$,
$N_{E_8}(A_3, A_1*A_2, A_1, A_1)= 42525$,
$N_{E_8}(A_3, A_1^3, A_1, A_1)= 30375$,
$N_{E_8}(A_1*A_2, A_1*A_2, A_1,  A_1)= 106650$,
$N_{E_8}(A_1*A_2, A_1^3, A_1, A_1)= 64125$,
$N_{E_8}(A_1^3, A_1^3, A_1, A_1)= 33750$,
$N_{E_8}(A_3,  A_2, A_2, A_1)= 10575$,
$N_{E_8}(A_3, A_2, A_1^2, A_1)= 29700$,
$N_{E_8}(A_3, A_1^2, A_1^2, A_1)= 76950$,
$N_{E_8}(A_1*A_2, A_2, A_2, A_1)= 35700$,
$N_{E_8}(A_1*A_2, A_2, A_1^2, A_1)= 84825$,
$N_{E_8}(A_1*A_2, A_1^2, A_1^2,  A_1)= 171450$,
$N_{E_8}(A_1^3, A_2,\break A_2, A_1)= 25125$,
$N_{E_8}(A_1^3, A_2, A_1^2, A_1)= 55125$,
$N_{E_8}(A_1^3, A_1^2,  A_1^2, A_1)= 94500$,
$N_{E_8}(A_3,\break A_2, A_1, A_1, A_1)= 91125$,
$N_{E_8}(A_3, A_1^2, A_1, A_1, A_1)= 243000$,
$N_{E_8}(A_1*A_2, A_2, A_1, A_1,\break A_1)= 276750$,
$N_{E_8}(A_1*A_2, A_1^2, A_1, A_1, A_1)= 597375$,
$N_{E_8}(A_1^3,  A_2, A_1, A_1, A_1)= 185625$,
$N_{E_8}(A_1^3, A_1^2, A_1, A_1, A_1)= 354375$,
$N_{E_8}(A_3, A_1, A_1, A_1, A_1,  A_1)= 759375$,
$N_{E_8}(A_1*A_2,\break A_1, A_1, A_1, A_1, A_1)= 2025000$,
$N_{E_8}(A_1^3, A_1, A_1, A_1, A_1, A_1)= 1265625$,
$N_{E_8}(A_2, A_2, A_2,\break A_2)= 9350$,
$N_{E_8}(A_2, A_2, A_2, A_1^2)= 24975$,
$N_{E_8}(A_2, A_2, A_1^2,  A_1^2)= 64350$,
$N_{E_8}(A_2, A_1^2,\break A_1^2, A_1^2)= 143100$,
$N_{E_8}(A_1^2, A_1^2, A_1^2, A_1^2)= 261225$,
$N_{E_8}(A_2, A_2,  A_2, A_1, A_1)= 78000$,
$N_{E_8}(A_2, A_2, A_1^2, A_1, A_1)= 203625$,
$N_{E_8}(A_2, A_1^2, A_1^2, A_1, A_1)=  479250$,
$N_{E_8}(A_1^2, A_1^2, A_1^2,\break A_1, A_1)= 951750$,
$N_{E_8}(A_2, A_2, A_1, A_1, A_1, A_1)= 641250$,
$N_{E_8}(A_2, A_1^2, A_1, A_1, A_1, A_1)= 1569375$,
$N_{E_8}(A_1^2, A_1^2, A_1, A_1, A_1, A_1)= 3341250$,
$N_{E_8}(A_2,  A_1, A_1, A_1, A_1, A_1, A_1)=\break 5062500$,
$N_{E_8}(A_1^2, A_1, A_1, A_1, A_1, A_1, A_1)= 11390625$,
$N_{E_8}(A_1,  A_1, A_1, A_1, A_1, A_1, A_1,\break A_1)= 37968750$,
plus the assignments implied by
\eqref{Aa} and \eqref{Ab},
all other numbers $N_{E_8}(T_1,T_2,\dots,T_d)$ being zero.

\end{document}